\documentclass[reqno]{amsart}
\usepackage{amsmath,  amssymb, latexsym}
\newcommand{\be}{\begin{equation}}

\newcommand{\mc}{\mathcal}
      \newcommand{\ee}{\end{equation}}
      \newcommand{\ba}{\begin{eqnarray}}
       \newcommand{\ea}{\end{eqnarray}}
\newcommand{\ban}{\begin{eqnarray*}}
       \newcommand{\ean}{\end{eqnarray*}}
\newcommand{\f}{\frac}
\newcommand{\lr}{\lrcorner}
\newcommand{\p}{\partial}

\newcommand{\no}{\nonumber\\ &&}

 \newcommand{\Pf}{\noindent {\em Proof.} }

\newcommand{\sect}[1]{\section{#1} \setcounter{equation}{0}}

\newtheorem{theo}{Theorem}[section]

 \newtheorem{lem}[theo]{Lemma}

\begin{document}
\title{Existence, Convergence and Limit Map of the Laplacian Flow}
\author{Feng Xu \& Rugang Ye} 
\address{MSI, The Australian National University, ACT 0200, Australia \& Department of Mathematics, University of California, Santa Barbara, CA93106, USA}
\email{feng.xu@anu.edu.au, yer@math.ucsb.edu}
\date{}
\begin{abstract}
We prove short time existence and uniqueness of the Laplacian flow starting at an arbitrary closed $G_2$-structure. We establish long time existence and convergence of the Laplacian flow starting near a torsion-free $G_2$-structure. We analyze the limit map of the Laplacian flow in relation to the moduli space of 
torsion-free $G_2$-structures. We also present a number of   
results which constitute a fairly complete 
algebraic and analytic basis for studying the Laplacian flow.
\end{abstract}
\maketitle

\noindent {Table of Contents} \\
1. Introduction \\
2. $G_2$-structures\\
3. An identity for the Hodge Laplacian on 3-forms\\
4. Additional differential identities\\
5. Linear parabolic theory \\
6. Short time solutions of the gauge fixed Laplacian flow\\
7. Short time solutions of the Laplacian flow \\
8. Long time existence and convergence of the gauge fixed Laplacian flow \\
9. Long time existence and convergence of the Laplacian flow\\
10. The limit map of the Laplacian flow\\
Appendix: Space-time function spaces\\
 
\sect{Introduction}
\vspace{2mm}

The Riemannian holonomy group of a 7-dimensional manifold $M$  
equipped with a torsion-free $G_2$-structure is contained in the Lie group $G_2$. As a consequence, 
$M$ is Ricci flat. If the fundamental group of $M$ is finite, then the holonomy group of $M$ actually equals $G_2$ and the spinor bundle of $M$ splits off a parallel $\mathbb{R}$ summand w.r.t.~the Levi-Civita connection. These 
properties are the main reason for the importance of 
torsion-free $G_2$-structures and in general,  
$G_2$-structures, in differential geometry. 
In particular, $G_2$ holonomy appears as one important case of the Berger classification of holonomy groups of Riemannian manifolds.  Note that manifolds with 
$G_2$ holonomy play an important role in $M$-theory.
Namely the compactification of $M$-theory on a manifold with $G_2$ holonomy leads to an ${\mc N}=1$ (3+1)-dimensional quantum field theory, which is similar to the compactification of 
heterotic string theory on Calabi-Yau manifolds.
(The parallel $\mathbb{R}$ summand of the spinor bundle provides the ${\mc N}=1$ supersymmetry.)

A fundamental problem here is how to deform a 
given $G_2$-structure on a manifold to a torsion-free
$G_2$-structure.  R.~Bryant proposed the following 
Laplacian flow for closed $G_2$-structures
\ba \label{lpflow}
\f{\p \sigma}{\p t}=\Delta_{\sigma} \sigma,
\ea
where $\Delta_{\sigma}$ denotes the Hodge Laplacian 
of the Riemannian metric induced by the $G_2$-structure 
$\sigma$, cf.~[B2]. In [BX], this flow is interpreted as the gradient flow of Hitchin's volume functional 
w.r.t.~an unusual metric. 

It turns out that the structures of the Laplacian flow are rather complicated. On the other hand, the Laplacian flow shares some features with the Ricci flow, which is worth noting.  
Under the Laplacian flow, the induced metric $g=g(t)$ evolves as follows [B2]
\ba
\f{\p g}{\p t}=-2Ric+\f{8}{21}|\tau|^2 g + \f{1}{4} 
{\bf j}(*(\tau \wedge \tau)),
\ea
where $\tau$ denotes the adjoint torsion of $\sigma$ (cf.~Section 2), $*$ the Hodge star, and $\bf j$ a certain linear operator associated with $\sigma$. Thus we see the leading part 
$-2Ric$ which appears in the Ricci flow. The perturbation part is given by the adjoint torsion $\tau$, which is a key quantity because its vanishing is equivalent to the torsion-free condition.  One may wonder here whether it is possible to use the Ricci flow to deform $G_2$-structures. However, the induced metric 
of a $G_2$-structure $\sigma$ does not determine
$\sigma$ completely, and in general a metric may not be induced from a $G_2$-structure.  Hence a suitable additional coupled equation would be needed in order to use the Ricci flow to deform 
$G_2$-structures. Based on some calculations we are convinced that such coupled equations do not exist in general.

There are four main parts in this paper.  First, 
we prove the short time existence and uniqueness of solutions of the Laplacian flow with given initial data. Second, we establish the long time existence and convergence of the Laplacian flow starting near torsion-free $G_2$-structures. Third, we establish the smoothness of 
the limit map of the Laplacian flow around torsion-free $G_2$-structures,  and determine the limit projection of the Laplacian flow into the moduli space of torsion-free $G_2$-structures. This reveals the deep relation of the Laplacian flow
with the moduli space of torsion-free $G_2$-structures.
Fourth, we present a number of   
results which constitute a fairly complete 
algebraic and analytic basis for studying the Laplacian flow. These include algebraic formulas, differential identities, a linear parabolic theory, and a detailed analysis of the basic analytic structure of the Laplacian flow and the gauge fixed Laplacian flow. 

The highlights of our main results in the first three parts are formulated in the following three main theorems. (We refer to the subsequent sections for the results in the last part.)  A {\it closed solution} means a solution given by closed 
$G_2$-structures. 

\begin{theo} \label{existence} Let $M$ be a compact 7-dimensional manifold. Let $\sigma_1$ be a closed   
$G_2$ structure of class $C^{4+\mu}$ on $M$ for 
some $0<\mu<1$.  Then there is a closed
${\mc C}^{2+\mu, (3+\mu)/2}$ solution 
$\sigma=\sigma(t)$ of the Laplacian flow on 
a time interval $[0, T]$ with $T>0$,
such that 
$\sigma(0)=\sigma_1$. This solution
is unique among all ${\mc C}^{2+\mu,(2+\mu)/2}$
functions (with closed $G_2$-structures as values) with the initial value $\sigma_1$.
For each $0<\epsilon<T$, there is a family of diffeomorphisms $\phi(\cdot, t)$ of class ${\mc C}^{3+\mu, (4+\mu)/2}$ on $[\epsilon, T]$, such that $\tilde \sigma(t)=\phi(\cdot, t)^*\sigma(t)$ is
a $C^{\infty}$ solution of the Laplacian flow.  Moreover, there holds $\sigma \in {\mc C}^{l-2, (l-1)/2}$ on 
$M \times [0, T]$, provided that $\sigma_1 \in C^l$ for a noninteger $l>4+\mu$. In particular, $\sigma$ is smooth 
if $\sigma_1$ is smooth. 
\end{theo}

For an estimate for $T$ from below and other a priori estimates, we refer to   Theorem \ref{gauge1}, Theorem \ref{gauge2} and the proof of Theorem \ref{existence}.  The definitions of the involved function spaces are given in Appendix. In particular, the space $C^{4+\mu}$ means the H\"{o}lder space $C^{4, \mu}$ in the conventional notation. The parabolic H\"{o}lder spaces, i.e.~the ${\mc C}^{l, l/2}$ spaces, and their generalizations ${\mc C}^{l, l'/2}$ spaces, involve spacial and time derivatives 
in patterns which are particularly suitable for handling 
second order partial differential equations of parabolic type or related types. Several statements of this theorem actually hold true under more general or weaker assumptions. 
On the other hand, short time existence and uniqueness of solutions of Sobolev classes can also be obtained for Sobolev initial data.

\begin{theo}\label{stability} 
Let $\sigma_0$ be a smooth torsion-free $G_2$-structure on a compact  manifold $M$ of dimension 7. Let $0<\mu<1$. Then there exists a strong $C^{2+\mu}$-neighborhood $\mathcal{U}_{\sigma_0}$ of $\sigma_0$ in the space of closed smooth  $G_2$-structures on $M$ such that whenever $\sigma_1\in \mathcal{U}_{\sigma_0},$
 the Laplacian flow (\ref{lpflow}) starting at $\sigma_1$ has a unique closed smooth solution $\sigma=\sigma(t)$ on $M \times [0, \infty)$  which converges exponentially to a smooth torsion free $G_2$-structure $\sigma_{\infty}$ as $t\rightarrow \infty.$  Thus torsion-free $G_2$-structures are 
stable  in the space of closed $G_2$-structures with respect to the Laplacian flow.
\end{theo}

\begin{theo} \label{limitmap} Let $\mc{F}$ denote the limit map of the Laplacian flow in the situation of Theorem \ref{stability}, i.e.~$\mc{F}(\sigma_1)=\sigma_{\infty}$. Then 
$\mc{F}: \mc{U}_{\sigma_0} \rightarrow 
\mc{T}$ is a smooth map, where $\mc{T}$ denotes the space of smooth torsion-free $G_2$-structures on $M$.   Moreover, there holds
\ba
\pi \circ \mc{F}=\Pi,
\ea
where $\Pi$ is a canonical 
projection into the moduli space $\mc{T}/\text{Diff}\,_0(M)$ of smooth torsion-free 
$G_2$-structures on $M$ and $\pi$ is the quotient 
projection from $\mc{T}$ onto $\mc{T}/\text{Diff}\,_0(M)$.
\end{theo}

For relevant definitions (such as strong 
$C^{2+\mu}$ neighborhood and the projection $\Pi$) we refer to Sections 9 and 10. Next we explain the backgrounds and main ideas of the above results. \\ 

\noindent {\bf Existence and uniqueness of short time solutions}

 As it turns out, existence and uniqueness of short time solutions of the Laplacian flow are a rather delicate problem. Indeed, the Laplacian flow is not a parabolic equation, and there  seems to be no way to restore full parabolicity for it by a transformation such as gauge fixing as employed in the DeTurck trick for the Ricci flow. Previously, it was proved in [BX] via rather complicated computations  that a partial parabolicity, namely the parabolicity in the direction of closed forms, can be restored for the Laplacian flow by a certain gauge fixing, i.e. the gauge fixed Laplacian flow is parabolic in the said direction.
However, the gauge fixed Laplacian flow fails to be 
parabolic in the complementary directions. To cope with this situation of lack of full parabolicity, the set-up of Fr\'{e}chet space of smooth forms
and Nash-Moser implicit function theorem were employed in [BX].  This way, short time existence and uniqueness of closed smooth solutions of the Laplacian flow starting at a smooth closed $G_2$-structure were obtained in [BX].

In this paper, we first introduce a new gauge fixing for the Laplacian flow, which restores the partial parabolicity, i.e. the parabolicity in the direction of closed forms.  This new gauge fixing is 
simpler and more transparent than the one used in [BX],
and is based on a new identity for the Hodge Laplacian, 
which in turn is based on some delicate differential identities involving splittings of the exterior differential via irreducible $G_2$-representations. 
Next we develop a new linear parabolic theory for closed forms which is tailored to handle operators which are 
only parabolic in the direction of closed forms.  Using this theory we are then able to establish the short time 
existence and uniqueness of the Laplacian flow starting  at $C^{4+\mu}$ initial data and obtain estimates depending only on the $C^{4+\mu}$ properties of the initial data.  (Note that we avoid using  the 
Nash-Moser implicit function theorem.) 

The improvement to $C^{l}$ initial data with $l\ge 4+\mu$ provided by Theorem \ref{existence} in comparsion with the result in [BX] is an obvious analytic aspect of it. More important is the complete understanding and resolution of the problem of short time solutions of the Laplacian flow.  The existence for $C^{l}$ initial data, the associated estimates, as well as the basic 
$\mc{C}^{l, l/2}$ set-up also play an important role for establishing the long time existence and convergence of the Laplacian flow and the smoothness of its limit map 
as presented in Theorem \ref{stability} and Theorem \ref{limitmap}, as will further be explained below. 
Moreover, the framework and strategy for Theorem \ref{existence} also allow to handle e.g.~the Laplacian flow on complete noncompact manifolds. This will be presented in a subsequent paper. \\

\noindent {\bf Long time existence and convergence}

The second main theorem of this paper, Theorem \ref{stability}, provides the first result on long time behavior of the Laplacian flow. From the dynamical point of view, 
this result can be viewed as stability of torsion-free $G_2$-structures in regard to the Laplacian flow and the Hitchin volume functional. As the Laplacian flow is very natural 
geometrically, this dynamical stability is also very natural from a geometric point of view. We also believe that it is significant for 
the $M$-theory. Previously, the dynamical stability of Einstein metrics w.r.t.~the volume-normalized Ricci flow and that of 
Ricci flat metrics w.r.t.~the Ricci flow were proved by the second named author under the condition of positive first eigenvalue of the Lichnerowitz Laplacian, as a consequence of a general convergence result for the Ricci flow [Y1]. We would like to mention that we have also obtained a general long time convergence result for the Laplacian flow under the assumption of small 
torsion of the initial $G_2$-structure [XY1].  

The basic scheme of the proof for Theorem \ref{stability} 
is to derive exponential decay estimates for the solution under the assumption of certain smallness and boundedness. The said smallness and boundedness on a small time interval follow from our results on short time solutions, but are 
not known a priori for all time. Hence it is crucial to obtain strong feedback via exponential decay, such that 
they can be shown to always hold true.  The key starting point of the exponential decay is the exponential $L^2$-decay, which is based on the spectral property of the Hodge Laplacian. Such an exponential $L^2$-decay scheme was first 
implemented successfully in [Y1] for proving long time 
convergence of the Ricci flow. The situation in this paper is more delicate for the following reason. The 
involved PDE has a second order perturbation term besides
the leading Laplacian term, which makes it more difficult to apply the maximum principle to convert $L^2$ estimates into $C^0$ estimates.
Moreover, for the purpose of establishing the smoothness of the limit map of the Laplacian flow, we need to derive linear power  
decay estimates rather than estimates with fractional powers. Here the conventional 
$L^2$ version of Moser type maximum principle is not suitable. We derive an $L^1$ version instead and apply it to overcome the trouble. 

Another tool employed here is a result on the local smooth structure of the moduli space of torsion-free $G_2$-structures. 
It is used to locate the target torsion-free $G_2$-structure for the Laplacian flow to converge to.
(The actual limit differs from this target by a 
diffeomorphism.) The said result  is a refinement of D.~Joyce's well-known result [J] on the same topic, and its proof is presented in [XY2].
Note that this result can be viewed as the stationary version of Theorem \ref{existence}. Indeed, it is in part based 
on our understanding of some features of the Laplacian flow, see [XY2] for details. 

Note that the analytic set-up for the above scheme of 
exponential decay has to be carefully chosen. Indeed, the ${\mc C}^{l, l/2}$ spaces and the parabolic estimates in Section 5 play a crucial role here. The main reason for this is that the estimates in these spaces require minimal amount of bounds while 
providing strong control directly, in contrast to 
e.g.~Sobolev space estimates which leave a large gap 
because of Sobolev embeddings. \\

\noindent {\bf The limit map and its projection into the moduli space}

Our third main theorem, Theorem \ref{limitmap} (and additional results in Section 10),  is the first result of its kind regarding the limit map of a nonlinear geometric evolution equation, the space of 
its stationary solutions, and 
the corresponding moduli space. Besides individual 
torsion-free $G_2$-structures, their moduli space 
is an important geometric object. In particular, it plays an important role in M-theory. It is therefore 
very desirable to understand the relation between 
the Laplacian flow and the space of torsion-free
$G_2$-structures and the associated moduli space.
(The application of the result in [XY2] mentioned above 
is only one aspect of this relation.) 
Theorem \ref{limitmap} provides a complete understanding of this  deep relation.

The smoothness of the limit map of the Laplacian flow 
is rather intricate. Indeed, one encounters anlaytic troubles if one deals with the Laplacian flow directly. 
Indeed, the equation satisfied by the difference of two solutions of the Laplacian flow (with two different initial values) fails to be parabolic, and hence it is not clear how to derive estimates for this difference directly. 
(This goes back to the lack of parabolicity of the Laplacian flow itself. We are able to handle it in the context of short time solutions by a suitable gauge fixing as explained 
before.) Our basic strategy for proving the said smoothness is to go through the gauge fixed Laplacian flow. The proof requires a number of additional ingredients, and involves various exponential decay estimates. Indeed, the linear theory in Section 5,  Theorem \ref{stability} and the techniques in its proof have to be applied in various fashions. In particular, as mentioned above, the 
linear power nature of the decay estimates is crucial here. 

Finally, the identification of the limit projection of the Laplacian flow in terms of a canonical projection is achieved via the detailed convergence analysis of the Laplacian flow.  \\

Now some additional brief descriptions of the main content of the subsequent sections. In Section 2, we present a short introduction to the basics of $G_2$-structures. We explain the basic concepts, present some useful facts and algebraic formulas, such as the important $G_2$-irreducible decompostions of forms and associated formulas, and also provide some basic set-ups of this paper. 
In Section 3, we derive the new identity for the Hodge Laplacian mentioned above. Along the way, we present a detailed treatment of a typical one of Bryant's differential identities, and also derive a new one.  In Section 4, we present some additional differential identities. Note that the differential identities in these two sections are tied to the irreducible decompositions of forms and are only available when a $G_2$-structure is present. Obviously, the applications of these differential identities in the study of the Laplacian flow as presented in [BX] and this paper  offer a unique 
new perspective in geometric analysis and nonlinear analysis. In Section 5, we develop the new linear parabolic theory 
for closed forms described above.  A subtle point here is that the corresponding linear parabolic problem for exact forms is ill-behaved due to the lack of completeness of some involved function spaces.
(This phanomenon is uncovered for the first time 
in this paper.) In Section 6, we construct our new 
gauge, which is motivated by the Hodge Laplacian identity in Section 3, and apply the theory in Section 5 and the classic inverse function theorem to prove existence and uniqueness of short time solutions of the 
gauge fixed Laplacian flow.     Note that the inverse function theorem immediately implies a local uniqueness. To obtain global uniqueness, we utilize the special quadratic structure 
in the equation to handle its nonlinear second order perturbation part, and appeal to  the Bochner-Weitzenb\"{o}ck formula for the Hodge Laplacian. In Section 7, we apply the results of Section 6 to derive existence and uniqueness of short time solutions of the Laplacian flow. Here the differential identities 
in Section 4 play an important role.

In Section 8, we prove long time 
convergence at exponential rate of the gauge fixed Laplacian flow 
starting near a torsion-free $G_2$-structure. In Section 9 we combine the result in the previous section and 
results on the local smooth structure of the moduli space of torsion-free $G_2$-structures to derive long time existence and convergence of the Laplacian flow.
In the last section, we prove the smoothness of the limit map of the Laplacian flow and identify its projection into the moduli space of torsion-free 
$G_2$-structures.\\

The first named author would like to thank Prof.~Robert Bryant, Prof.~Mike Eastwood and Prof.~Mark Haskins 
for stimulating discussions on the subjects of this paper. He would also like to thank MSRI for its hospitality when he was a postdoctoral fellow there.\\

\sect{$G_2$-structures}
\vspace{2mm}

In Subsection 2.1 we present some basics of $G_2$-structures. In Subsection 2.2, we describe 
the decompositions of forms into irreducible components, which play a crucial role for various computations in this paper.  

\subsection{Basics}
\hspace{1cm} 
\vspace{1mm}

Let $e_i, i=1,2...,7$ denote the standard orthonormal 
basis of ${\mathbb R}^7$ and $e^i=dx^i$ its dual basis. 
The standard $G_2$-structure on ${\mathbb R}^7$  
is 
\ba \label{sigma}
&&\sigma_{{\mathbb R}^7}=e^1 \wedge (e^2 \wedge e^3+
e^4 \wedge e^5+e^6 \wedge e^7)+e^2 \wedge (e^4 \wedge
e^6-e^5 \wedge e^7) \nonumber \\
&&\hspace{1.5cm}-e^3 \wedge (e^4 \wedge e^7+e^5 \wedge e^6) \nonumber \\
&&\hspace{.8cm}=e^1 \wedge \omega_{{\mathbb R}^6}+\text{Re}\,\Omega_{{\mathbb C}^3},
\ea
where $\omega_{{\mathbb R}^6}=e^2 \wedge e^3+e^4 \wedge e^5 + e^6 \wedge e^7$ is the standard symplectic form on 
${\mathbb R}^6$ and $\Omega_{{\mathbb C}^3}=dz^1\wedge dz^2 \wedge dz^3$ 
is the standard holomorphic volume form on ${\mathbb C}^3={\mathbb R}^6$, 
w.r.t.~the decomposition ${\mathbb R}^7={\mathbb R}\oplus 
{\mathbb R}^6$. (Thus $z^1=x^2+\sqrt{-1}x^3, z^2=x^4+\sqrt{-1}x^5$ and 
$z^3=x^6+\sqrt{-1}x^7$.) 
 The group $G_2$ can be defined as follows 
\ba
G_2=\{A \in GL(7, {\mathbb R}): A^*\sigma_{{\mathbb R}^7}=
\sigma_{{\mathbb R}^7}\}.
\ea
It is a 14-dimensional compact, connected, simply-connected and simple Lie subgroup of $SO(7)$, 
cf.~[B1][B2].

Set $\Lambda^3_+({\mathbb R}^7)^*=\{
L^*\sigma_{{\mathbb R}^7}: L \in GL({\mathbb R}, 7)\}$. This is the set of constant $G_2$-structures on ${\mathbb R}^7$. It is open in $\Lambda^3 ({\mathbb R}^7)^*$, cf.~[B2]. Let $M$ be a smooth 7-dimensional manifold. 
For each $p\in M$, set $\Lambda^3_+T_p^*M=
\{\sigma\in \Lambda^3T^*_pM: 
\sigma=L^* \sigma_{{\mathbb R}^n} \mbox{ for an isomorphism }L: T_pM \rightarrow {\mathbb R}^7\}$.  Then 
$\Lambda^3_+(T^*M)=\cup_{p\in M} \Lambda^3_+T^*_pM$ is 
an open subbundle of $\Lambda^3T^*M$ as a fiber bundle.
This is the bundle of {\it positive} 3-forms. 
An $A$ in $\sigma=A^* \sigma_{{\mathbb R}^n}$ is called 
an {\it inducing map} of $\sigma$. An {\it induced orthonormal basis} for $\sigma$ is $A^{-1}e_1,...,A^{-1}e_7$, where $A$ is an inducing map of $\sigma$. \\

\noindent {\bf Definition 2.1} Let $l\ge 0$.  
$G_2$-structures of 
class $C^l$, or $C^l$ $G_2$-structures are defined to be 3-forms of class $C^l$ with values  in $\Lambda^3_+T^*M$.  In other words, they are 
$C^l$ sections of $\Lambda^3_+T^*M$.
(It is easy to show that they are in one-to-one correspondance with $C^l$ principal $G_2$ subbundles of the principal frame bundle of $M$.)\\

Note that the existence of a $G_2$-structure 
(of class $C^l, l\ge 0$) is equivalent to the vanishing 
of the first two Stiefel-Whitney classes, i.e. 
equivalent to $M$ being orientable and spinnable, 
cf.~[B2]. 

Since $G_2 \subset SO(7)$, a $C^l$ $G_2$-structure $\sigma$ induces a $C^l$ Riemannian metric $g_{\sigma}$
on $M$
and an $C^l$ orientation of $M$, namely 
$g_{\sigma}=L^*g_{{\mathbb R}^7}$ and $dvol_{\sigma}
=L^*(e^1\wedge \cdot \cdot \cdot \wedge e^7)$, if $\sigma=L^* \sigma_{{\mathbb R}^7}$.  All quantities associated with $g_{\sigma}$ will often be indicated by the subscript 
$\sigma$. For example, $*_{\sigma}$ denotes the Hodge 
$*$ of $g_{\sigma}$. Note that $g_{\sigma}$ can
be given by an explicit algebraic formula in terms of $\sigma$. Indeed there holds, as is easy to verify
\ba \label{metric}
g_{\sigma}(u, v)=\f{(u \lr \sigma) \wedge (v \lr \sigma) \wedge \sigma}{6dvol_{\sigma}}
\ea
at each $p\in M$ and for all $u, v \in T_pM$. Moreover, there holds  
\ba \label{volume}
dvol_{\sigma}=6^{-\f{7}{9}} (det_{\Omega} \sigma)^{\f{1}{9}}\Omega,
\ea
where $\Omega$ denotes an arbitrary volume form at any given $p$ 
(i.e.~a nonzero element of $\Lambda^7_pT^*M$), and 
the determinant $det_{\Omega} \sigma$ is defined to be 
the determinant of the quadratic form 
$(u \lr \sigma) \wedge (v \lr \sigma) \wedge \sigma)
/\Omega$
on a basis $u_1,...,u_7$ such that $\Omega(u_1,...,u_7)=1$. Hence the formula (\ref{metric}) gives the metric $g_{\sigma}$ explicitly in terms of $\sigma$.

Next we note the following 
simple, but important fact. 

\begin{lem} \label{simple} There are universal positive numbers $\epsilon_0\le 1$, $\mu_0$ and $C_0$ with the following property. Let $p\in M$. If $\sigma \in \Lambda^3_+(T^*_pM), \gamma 
\in \Lambda^3T^*_pM$ and $|\gamma-\sigma|_{\sigma}\le
\epsilon_0$, then $\gamma \in \Lambda^3_+T^*_pM$.
Moreover, there holds $|g_{\gamma}|_{\sigma}\le C_0$ and the eigenvalues of $g_{\gamma}$ w.r.t.~
$g_{\sigma}$ are bounded below by $\mu_0$.
\end{lem}
\Pf Since $\Lambda^3_+({\mathbb R}^7)^*$ is open in 
$\Lambda^3 ({\mathbb R}^7)^*$, there is a positive number $\epsilon_0\le 1$ 
such that $\gamma \in \Lambda^3_+ ({\mathbb R}^7)^*$ 
whenever $\gamma \in  \Lambda^3 ({\mathbb R}^7)^*$ and 
$|\gamma-\sigma_{{\mathbb R}^7}|\le \epsilon_0$. By continuity and compactness, there is a positive number 
$\mu_0$ such that the eigenvalues of $g_{\gamma}$ w.r.t.~the Euclidean metric are bounded from below 
by $\mu_0$.  
The claims of the lemma then follow from the induced nature of $g_{\sigma}$. \qed\\

\noindent {\bf Definition 2.2} The {\it total torsion} of 
a $G_2$-structure $\sigma$ is defined to be 
$\nabla_{\sigma} \sigma$. Its {\it adjoint torsion} 
$\tau=\tau_{\sigma}$ is defined to be 
\ba
\tau=d^*_{\sigma} \sigma=-*_{\sigma} d *_{\sigma} \sigma.
\ea
Note that $d\tau=\Delta_{\sigma} \sigma$, if $\sigma$ is a {\it closed} $G_2$-structure, i.e.~a $G_2$-structue which is a closed form.
A $G_2$-structure is said to be {\it torsion-free}, provided that 
its total torsion vanishes everywhere. 
(If we do not specify the $C^l$ class of $\sigma$ in a discussion, then $\sigma$ is assumed to be in $C^l$ for the minimal $l$ as required in the discussion.)  \\

A fundamental fact [B1][B2][FG][S] is that a $G_2$-structure is torsion-free precisely when 
the induced metric has a subgroup of $G_2$ as its holonomy group and hence is Ricci-flat. 
On the other hand, it is well-known [B2][FG][S] that a  $G_2$-structure $\sigma$ is torsion-free
precisely when it is closed and its adjoint torsion vanishes, i.e.~when $\sigma$ is a harmonic form (w.r.t.~$g_{\sigma}$) in the case of a compact $M$.  Indeed, the full torsion $\nabla_{\sigma}  \sigma$ can be expressed in terms of $d\sigma$ and $\tau_{\sigma}$, which follows from the arguments in [Proof of Proposition 2, B2], see also [Theorem 2.27, K] for an explicit formula.  This explicit formula leads to the following lemma regarding 
closed $G_2$-structures.

\begin{lem} \label{Klemma}  Let $\sigma$ be a closed $G_2$-structure. Then there holds at each $p\in M$
\ba
\nabla_{\sigma} \sigma=-\f{1}{3} <\tau_{\sigma}, *_{\sigma} \sigma>_{2,1},
\ea
where $<\cdot, \cdot>_{2,1}$ denotes the contraction
$<\cdot, \cdot>_{2,1}: \otimes^2 T_p^*M \times 
(T_p^*M \otimes \Lambda^3 T^*_pM) \rightarrow 
T_p^*M \otimes \Lambda^3 T_p^*M$ given by 
\ba
<\alpha_1 \otimes \alpha_2, 
\alpha_3 \otimes \gamma>_{2,1}=(\alpha_2 \cdot \alpha_3)
\alpha_1 \otimes \gamma
\ea
for $\alpha_1,\alpha_2, \alpha_3 \in T^*_pM$ and $\gamma 
\in \Lambda^3 T^*_pM$. (Note that $\Lambda^4 T_p^*M 
\subset T_p^*M \otimes \Lambda^3 T^*_pM$. For relevant discussions of a similar contraction see Lemma \ref{adjoint} below.)
\end{lem} 
\Pf This is a reformulation of [Theorem 2.27, K] in the special case of a closed $G_2$-structure. \qed \\

\subsection{Irreducible Decomposition of Forms}
\hspace{1cm}
\vspace{1mm}

Let $\sigma$ be a $G_2$-structure on $M$.  For each $p\in M$ and $1\le j\le 7$, the exterior space $\Lambda^j T^*_pM$ decomposes orthogonally into irreducible representations of $G_2$,
which then leads to the corresponding decompositions 
of the bundles $\Lambda^j T^*M$, and hence of differential $j$-forms. We have [B2]
\ba \label{charac}
&&\Lambda^3 T^*M=
\Lambda^3_1 (T^*M) \oplus \Lambda^3_7 (T^*M)
\oplus \Lambda^3_{27} (T^*M), \nonumber \\
&&\Lambda^2 T^*M
=\Lambda^2_7(T^*M) \oplus \Lambda^2_{14} (T^*M),
\Lambda^1T^*M= \Lambda^1_7 (T^*M),
\ea
and the corresponding ones $\Omega^3(M)=
\Omega^3_1(M)\oplus \Omega^3_7(M)\oplus
\Omega^3_{27}(M)$ etc.~(as well as for forms of various 
$C^l$ classes), where the subscript indicates the dimension of representation.
We have the characterizations 
\ba \label{decomp-ch}
&&\Lambda^3_1(T^*M)
=\{c\sigma_p: c \in \mathbb R, p\in M\}, 
\Lambda^3_7(T^*M)=\{*_{\sigma}(\alpha \wedge \sigma): \alpha \in T^*M\}, \nonumber \\ &&\Lambda_{27}^3(T^*M)=\{\gamma \in \Lambda^3(T^*M): \gamma \wedge \sigma=0, \gamma\wedge *_{\sigma} \sigma=0\}, \nonumber \\
&&\Lambda^2_7(T^*M)=\{*_{\sigma} (\alpha \wedge *_{\sigma}\sigma): \alpha \in \Lambda^1T^*M\}=
\{\alpha\in \Lambda^2T^*M: \alpha \wedge \sigma=2*_{\sigma}\alpha\}, \nonumber \\ 
&&\Lambda^2_{14}(T^*M)=\{\alpha \in \Lambda^2T^*M: 
\alpha \wedge \sigma=-*_{\sigma} \alpha\},
\ea
 cf.~[B2].  (Obviously, e.g.~$\gamma\wedge \sigma$ means $\gamma \wedge \sigma_p$ for $\gamma \in \Lambda^3(T_p^*M)$.) 
It follows that
\ba
\pi^2_7 \alpha=\f{1}{3} \alpha+\f{1}{3} *_{\sigma}(\alpha \wedge \sigma), \,\, \pi^2_{14} \alpha=\f{2}{3}\alpha-\f{1}{3}*_{\sigma}(\alpha \wedge \sigma),
\ea
where $\pi^i_j$ denotes the orthogonal projection from 
$\Lambda^iT_p^*M$ to $\Lambda^i_{j} (T^*_pM), p\in M$.
On the other hand, by (\ref{decomp-ch}), the formula for the decomposition of 
$\gamma \in \Lambda^3T_p^*M$ for 
$p\in M$ can be written as follows
\ba \label{decomp}
\gamma=f^0 \sigma+ *_{\sigma} (f^1 \wedge \sigma)
+f^3,
\ea
with $f^0 \in \mathbb{R}, f^1 \in T_p^*(M)$ and 
$f^3=\pi^3_{27}\gamma$. ($\sigma$ stands for $\sigma_p$.)
We present a formula for computing $f^1$, which will be needed later. For this purpose, we first present two lemmas, which will also be useful for other purposes. 

\begin{lem} \label{productlemma} Let $p\in M$ and $\alpha_1, \alpha_2 \in T_p^*M$. Then there hold
\ba \label{4product}
(\alpha_1 \wedge \sigma) \cdot (\alpha_2 \wedge \sigma)=
4 \alpha_1 \cdot \alpha_2
\ea
and 
\ba \label{3product}
(\alpha_1 \wedge *_{\sigma}\sigma) \cdot (\alpha_2 \wedge 
*_{\sigma}\sigma)=3 \alpha_1 \cdot \alpha_2,
\ea
where $\sigma$ means $\sigma_p$.
\end{lem}
\Pf By the induced nature of the metric $g_{\sigma}$,
it suffices to consider the Eulidean space. So we can assume $\sigma_p=\sigma_{\mathbb R^7}$. By linearity, it 
suffices to verify (\ref{4product}) and 
(\ref{3product}) for $\alpha_1=e_i$ and $\alpha_2=e_j$.
Since $G_2\subset SO(7)$ and it acts transitively on unit vectors and 
on orthonormal pairs [B1][B2], we can assume $(i, j)=(1,1)$ or $(1,2)$. Now it is straightforward to 
verify 
\ba
(e^1 \wedge \sigma_{\mathbb{R}^7}) \cdot (e^1 \wedge \sigma_{\mathbb{R}^7})=4, \,\,
(e^1 \wedge \sigma_{\mathbb{R}^7}) \cdot (e^2 \wedge \sigma_{\mathbb{R}^7})=0.
\ea 
On the other hand, using the formula 
\ba
*\sigma_{\mathbb{R}^7}&=&e^1 \wedge e^5 \wedge e^6 \wedge e^7
+e^2 \wedge e^3 \wedge e^6 \wedge e^7
+e^2 \wedge e^3 \wedge e^4 \wedge e^5
+e^1 \wedge e^3 \wedge e^5 \wedge e^7 \nonumber \\
&&-e^1 \wedge e^3 \wedge e^4 \wedge e^6-
e^1 \wedge e^2 \wedge e^5 \wedge e^6-e^1 \wedge e^2 
\wedge e^4 \wedge e^7
\ea
it is also straightforward to verify
\ba
(e^1 \wedge *\sigma_{\mathbb{R}^7}) 
\cdot (e^1 \wedge *\sigma_{\mathbb{R}^7})=3, (e^1 \wedge *\sigma_{\mathbb{R}^7}) 
\cdot (e^2 \wedge *\sigma_{\mathbb{R}^7})=0.
\ea
\qed \\

\begin{lem} \label{adjoint} Let $p\in M$. Consider the linear map $\sigma_p \wedge:T_pM^* \rightarrow \Lambda^4T_p^*M$ and its adjoint
$(\sigma_p \wedge)^*: \Lambda^4T_p^*M \rightarrow 
T_pM^*$. We have the following formula
\ba
(\sigma_p \wedge)^*=\sigma_p \neg_{\sigma}|_{\Lambda^4 T_p^*M}, 
\ea
where $\neg_{\sigma}:\Lambda^3T_p^*M \times (\Lambda^3 T_pM^*\otimes T_p^*M) \rightarrow 
T_p^*M$ denotes the front contraction 
w.r.t.~$g_{\sigma}$, i.e.
\ba
\gamma_1 \neg_{\sigma} (\gamma_2 \otimes \alpha)=
(\gamma_1 \cdot \gamma_2) \alpha
\ea
for $\gamma_1, \gamma_2 \in \Lambda^3 T^*_pM$ and  $\alpha \in T_p^*M$. (Note that it equals $\f{1}{6}$ times the restriction of the front contraction between $\otimes^3 T_p^*M$
and $\otimes^4 T_p^*M$ which is given by
\ba
(\gamma_1, \gamma_2 \otimes \alpha) \mapsto
(\gamma_1 \cdot \gamma_2) \alpha
\ea
for $\gamma_1, \gamma_2\in \otimes^3T^*_pM$ and $\alpha
\in T_pM^*$. The factor $\f{1}{6}$ is due to the fact that the inner product between $\gamma_1, \gamma_2
\in \Lambda^3T_p^*M$ equals $\f{1}{6}$ times their inner product as elements of $\otimes^3T_p^*M$.) 
\end{lem}
\Pf It suffices to consider the Euclidean space. We need to verify 
\ba
(\sigma_{\mathbb{R}^7} \wedge \alpha) \cdot \gamma=
\alpha \cdot (\sigma_{\mathbb{R}^7} \neg \gamma)
\ea
for all $\alpha \in (\mathbb{R}^7)^*$ and $\gamma \in \Lambda^4 (\mathbb{R}^7)^*$. Since $G_2\subset SO(7)$ and it acts transitively on unit vectors, 
we can assume $\alpha=e_1$. There holds
\ba
\sigma_{\mathbb{R}^7} \wedge e^1=-e^1 \wedge e^2 \wedge 
e^4 \wedge e^6+e^1 \wedge e^2 \wedge e^5 \wedge e^7
+e^1 \wedge e^3 \wedge e^4 \wedge e^7+ e^1 \wedge e^3 \wedge e^5 \wedge e^6.
\ea
Hence we have for $\gamma=\sum_{i<j<k<l} a_{ijkl} e^i \wedge e^j \wedge e^k \wedge e^l$
\ba
(\sigma_{\mathbb{R}^7} \wedge e^1) \cdot \gamma
=-a_{1246}+a_{1257}+a_{1347}+a_{1356}.
\ea
On the other hand, there holds $e^1\wedge e^2 \wedge e^3 \neg \gamma=
a_{1234} e^4+a_{1235}e^5+a_{1236} e^6
+a_{1237} e^7$ etc.~and hence a straightforward calculation yields
\ba
e^1 \cdot (\sigma_{\mathbb{R}^7} \neg \gamma)=
-a_{1246}+a_{1257}+a_{1347}+a_{1356}.
\ea
\qed \\

Now we present the formulas for computing $f^0, f^1$ 
 and $f^3$ in (\ref{decomp}).

\begin{lem} \label{3-decomp} The forms $f^0, f^1$ and $f^3$ 
in (\ref{decomp}) can be computed from 
$\gamma$ as follows
\ba \label{formula-1}
f^0=\f{1}{7}\gamma \cdot \sigma,
f^1=-\f{1}{4} \sigma \neg_{\sigma}(*_{\sigma} \gamma),
f^3=\gamma-f^0-*_{\sigma}(f^1 \wedge \sigma).
\ea 
 In other words, we have
\ba \label{3proj}
\pi^3_1 \gamma=\f{1}{7} (\gamma \cdot \sigma) \sigma, 
\pi^3_7 \gamma=-\f{1}{4} *_{\sigma}((\sigma \neg_{\sigma}(*_{\sigma} \gamma))\wedge \sigma).
\ea
\end{lem}
\Pf Taking the inner product of (\ref{decomp}) with an 
arbitrary element $*_{\sigma}(f \wedge \sigma)$ of 
$\Lambda^3_7T_p^*M$ (with $f \in T_p^*M$) we deduce,
on account of Lemma \ref{productlemma}
\ba
\gamma \cdot *_{\sigma} (f \wedge \sigma)
&=& *_{\sigma}(f^1 \wedge \sigma) \cdot *_{\sigma}
(f \wedge \sigma) \nonumber \\
&=& 4 f^1 \cdot f.
\ea
Since $\gamma \cdot *_{\sigma} (f \wedge \sigma)=
-*_{\sigma} \gamma \cdot (\sigma \wedge f)$, we can apply Lemma \ref{adjoint} to arrive at the 
formula for $f^1$ in (\ref{formula-1}). The formula 
for $f^0$ in (\ref{formula-1}) is obtained by taking 
the inner product of (\ref{decomp}) with $\sigma$. 
\qed \\  

Note that a general $G_2$ structure $\sigma$
has four torsion forms $\tau_0, \tau_1, \tau_2$
and $\tau_3$, with e.g.~$\tau_2$ having values in  $\Lambda^2_{14}(T^*M)$, see [B2]. If $\sigma$ is closed, then its adjoint torsion $\tau$ is precisely $\tau_2$.
Indeed, in that case, we have by [Proposition 1, B2]
the equation $\tau_2 \wedge \sigma=d*_{\sigma} \sigma$.
Hence we have by the above characterizations $
\tau_2=-*_{\sigma} (\tau_2 \wedge \sigma)=-*_{\sigma}
d *_{\sigma} \sigma=\tau$.\\

\sect{An identity for the Hodge Laplacian on 3-forms}
\vspace{2mm} 

The main purpose of this section is to present a new
identity for the Hodge Laplacian, which will play a crucial role in Section 6 for constructing 
a suitable gauge fixing for the Laplacian flow.

\subsection{A New Differential Identity of First Order} 
\hspace{1cm}\vspace{1mm}

In [B2] Bryant introduced differentials (exterior 
derivatives) which are adapted to the above decompositions of forms. These adapted differentials are very natural and indeed unique up to zeroth order perturbations. 
He also found remarkable (but very natural) identities for those differentials [B2]. Here we present a typical one of them and obtain a new one.  

As before, a $G_2$-structure $\sigma$ on $M$ is given.\\

\noindent {\bf Definition 3.1} The differential $d^7_7:
\Omega^1(M)=\Omega^1_7(M) \rightarrow 
\Omega^1_7(M)$ is defined to be 
\ba \label{d77}
d^7_7 \alpha =*_{\sigma} d(\alpha \wedge *_{\sigma}\sigma)=*_{\sigma}(d\alpha \wedge *_{\sigma}\sigma-\alpha\wedge *\tau_{\sigma}).
\ea
The differential $d^7_{14}: \Omega^1_7(M)
\rightarrow \Omega^2_{14}(M)$ is defined to be 
\ba
d^7_{14} \alpha=\pi^2_{14} d \alpha.
\ea
(We can also define $d^7_7$ on $\Omega^2_7(M)$ and  $\Omega^3_7(M)$. But the formulas for $d^7_7$ on different spaces are different.  The situations with $d^7_{14}$ and other adapted differentials are similar.) \\

The differential identity (\ref{id}) below
without the lower order term can be found in [B2] for  the special case of a torsion-free 
$\sigma$. 
For the purpose of computations in dealing with the Laplacian flow, we need to understand the precise nature of the additional lower order term which appears in the identity in the general case.

\begin{lem} \label{idlemma} We drop the subscript 
$\sigma$ in the notations. There holds for all $\alpha\in \Omega^1(M)$
(or $\alpha \in C^1(T^*M)$)
\ba \label{id}
d\alpha=\f{1}{3}*(d^7_7 \alpha \wedge *\sigma)
+d^7_{14}\alpha+\f{1}{3}*\left(*\sigma \wedge *(\alpha \wedge *\tau)\right).
\ea
\end{lem}
\Pf  The identity (\ref{id}) is equivalent to  the identity
\ba \label{idd}
\pi^2_7d\alpha=\f{1}{3}*(d^7_7 \alpha \wedge *\sigma)+\f{1}{3}*(*\sigma\wedge*(\alpha \wedge *\tau)).
\ea
To prove (\ref{idd}) we first observe 
\ba \label{observ}
\pi^2_7 d\alpha=\pi^2_7(\sum_i e^i \wedge (e_i\lr \nabla 
\alpha))=F(\nabla \alpha),
\ea
where $e_i$ denotes a local orthonormal basis and
\ba
F_p(\Theta)=\pi^2_7(\sum_i e^i \wedge (e_i\lr \Theta))
\ea
for all $p \in M$ and $\Theta \in T^*_pM \otimes T_p^*M$.
On the other hand, we have 
\ba 
*(d^7_7 \alpha \wedge *\sigma)+*\left(*\sigma\wedge*(\alpha \wedge *\tau)\right)
&=&*(*\sigma \wedge (d^7_7\alpha+
*(\alpha \wedge *\tau))) \nonumber \\
&=&*\left( *\sigma \wedge *(d\alpha \wedge *\sigma)\right) \nonumber \\
&=&\tilde F(\nabla \alpha),
\ea
where
\ba
\tilde F_p(\Theta)=\sum_i*\left(*\sigma \wedge*((e^i \wedge e_i \lr \Theta) \wedge *\sigma)\right).
\ea 
By (\ref{charac}), $\tilde F_p$ has 
values in $\Lambda^2_7(T_p^*M)$.

It is easy to verify that $F$ and $\tilde F$ are independent of the choice of the basis. Let $F$ and $\tilde F$ stand for $F_p$ and $\tilde F_p$ respectively for an arbitrary $p\in M$.  They are 
linear maps from $T^*_pM \otimes
T^*_pM$ to $\Lambda^2_7(T_p^*M)$.  One readily verifies that they are $G_2$ equivariant.  Now we have the 
orthogonal decomposition into irreducible 
$G_2$ representations
\ba
T^*_pM \otimes T^*_pM=span (g_{\sigma}|_p)
\oplus S^2_0(T^*_pM)  \oplus \Lambda^2_7(T^*_pM)
\oplus \Lambda^2_{14}(T_p^*M),
\ea
where $S^2(T^*_pM)$ consists of traceless symmetric 
2-tensors. The dimensions of these representations are
obviously different from each other.  By Schur lemma, the restrictions of $F$ and $\tilde F$ to the complement of $\Lambda^2_7(T^*_pM)$ are trivial. On the other hand, it is easy to see that their restrictions $F_{\Lambda^2_7}$ and 
$\tilde F_{\Lambda^2_7}$ to $\Lambda^2_7(T^*_pM)$ are nontrivial. Indeed, we can choose the basis $e_i$ to be 
induced from the standard basis on ${\mathbb R}^7$ via  an inducing map of $\sigma$. Then the formula 
(\ref{sigma}) holds true for $\sigma$. Using it we easily deduce for $\Theta=e^1 \otimes e^2$ 
\ba \label{P}
&&F(e^1 \otimes e^2)=\pi^2_7(e^1 \wedge e^2)=\f{1}{3}
e^1 \wedge e^2+\f{1}{3}*(e^1 \wedge e^2 \wedge \sigma)
\nonumber \\&& \hspace{2cm}=\f{1}{3}(e^1\wedge e^2-e^4 \wedge e^7-e^5 \wedge e^6)
\ea
and 
\ba \label{P1}
&&\tilde F(e^1 \otimes e^2)=*(*\sigma \wedge *(e^1 \wedge e^2 \wedge *\sigma))=*(*\sigma \wedge e^3)\nonumber \\&& \hspace{2cm}=e^1\wedge e^2-e^4\wedge e^7-e^5 \wedge e^6.
\ea
 
By Schur's lemma, $F_{\Lambda^2_7}$ and 
$\tilde F_{\Lambda^2_7}$ are isomorphisms. Since 
$\Lambda^2_7(T^*_pM)$ is odd dimensional, the isomorphism $F_{\Lambda^2_7} 
\tilde F_{\Lambda^2_7}^{-1}$ has at least one nontrivial 
eigenspace. By the irreducibility we then conclude that 
it is a scalar multiple of the identity. By (\ref{P}) and (\ref{P1}) the scalar is $\f{1}{3}$. 
Hence we conclude that $F=\f{1}{3}\tilde F$, which leads to (\ref{idd}). 

An alternative proof of (\ref{observ}) is in terms of the characterization (\ref{decomp-ch}), the orthogonality relations and integration by parts, analogous to the proof of Lemma \ref{d**lemma} below. \qed\\  

Next we present the said new differential identity. 

\begin{lem} \label{newdiff} There holds for all $\alpha\in \Omega^1(M)$
\ba \label{id-new}
&&d\alpha=*(d^7_7 \alpha \wedge *\sigma)-*(d\alpha \wedge \sigma)+\f{1}{3}\xi\left((*(*\sigma \wedge *(\alpha \wedge *\tau_{\sigma}))\right),
\ea
where $\xi=\xi_{\sigma}$ is defined as follows
\ba \label{T}
\xi(\gamma)=\gamma+*(\sigma \wedge \gamma).
\ea
\end{lem}
\Pf 
By the identity (\ref{id}) and the formulas in 
(\ref{decomp-ch}) we deduce
\ba \label{iddd}
d\alpha \wedge \sigma=\f{2}{3}d^7_7 \alpha \wedge *\sigma
-*d^7_{14} \alpha+\f{1}{3}\sigma\wedge*\left(*\sigma\wedge*(\alpha \wedge *\tau_{\sigma})\right),
\ea
which leads to 
\ba \label{idddd}
*(d\alpha \wedge \sigma)=\f{2}{3}*(d^7_7 \alpha \wedge *\sigma)
-d^7_{14} \alpha+\f{1}{3}*\left[\sigma\wedge*\left(*\sigma\wedge*(\alpha \wedge *\tau_{\sigma})\right) \right].
\ea
Adding (\ref{id}) and (\ref{idddd}) we then arrive 
at (\ref{id-new}). \qed \\

\subsection{A New Identity for the Hodge Laplacian on 3-Forms}
\hspace{1cm}\vspace{1mm}

Let $\sigma$ be a given closed $G_2$-structure on $M$. In the ensuing compuations in this subsection, we'll drop the subscript $\sigma$. Thus $*=*_{\sigma}, \Delta=\Delta_{\sigma}$ and 
$\tau=\tau_{\sigma}$. For a closed form $\theta\in C^2(\Lambda^3 T^*M)$ we apply the decomposition (\ref{decomp}) and compute
\ba
&&-\Delta \theta =-*d*d\theta+d*d*\theta=d*d*\theta,
\ea
\ba \label{dtheta}
&&*d*\theta=*d*(f^0\sigma+*(f^1\wedge\sigma)+f^3) \nonumber \\&&\hspace{1.15cm}=*(df^0\wedge*\sigma+df^1\wedge\sigma+d*f^3)-f^0\tau.
\ea

Next we consider the differential operator $H=H_{\sigma}$:
\ba \label{H}
&&H(\theta)=*d*(\frac{4}{3}f^0\sigma+*(f^1\wedge\sigma)-f^3)
\nonumber \\
&& \hspace{1cm}=*(\frac{4}{3}df^0\wedge *\sigma
 +df^1\wedge \sigma-d*f^3)-\frac{4}{3}f^0\tau.
\ea
This is an important operator because of its role in the 
linearization of the Laplacian flow, as will be 
shown in Section 6 below. We would like to compute the difference $d\circ H-\Delta=d(H+*d*)$.
By (\ref{H}) and (\ref{dtheta}) there holds
\ba  \label{HH}
(H+*d*)\theta=\frac{7}{3}*(df^0\wedge *\sigma)
 +2*(df^1\wedge \sigma)-\f{7}{3}f^0\tau.
\ea 
We would like to convert the term 
$2*(df^1\wedge \sigma)$ involving the 2-form $df^1$ into 
an expression involving a 1-form. This is achieved by the following lemma.
The 2-form $df^1$ still appears in the new formula
(\ref{Hequation}), but 
is separated  from other quantities. Hence it disappears 
in (\ref{hodge}) because of differentiation. 

\begin{lem} There holds 
\ba \label{Hequation}
&&(H+*d*)\theta=\frac{7}{3}*(df^0\wedge *\sigma)
+2*(d^7_7 f^1\wedge *\sigma)-2df^1\no
+\f{2}{3}\xi\left(*(*\sigma \wedge *(f^1 \wedge *\tau))\right)-\f{7}{3}f^0\tau \no
=\frac{7}{3} (df^0)_{\sharp} \lr \sigma
+2 (d^7_7 f^1)_{\sharp}\lr \sigma-2df^1\no
+\f{2}{3}\xi\left(*(*\sigma \wedge *(f^1 \wedge *\tau))\right)-\f{7}{3}f^0\tau. \nonumber \\
\ea
Consequently, 
\ba \label{hodge}
&&\Delta \theta =d(H(\theta))-d(\frac{7}{3}(df^0)_{\sharp} \lr  \sigma
+2(d^7_7f^1)_{\sharp} \lr \sigma)\nonumber \\&& \hspace{1.2cm}+d\left[\f{7}{3}f^0\tau-\f{2}{3}\xi\left(*(*\sigma \wedge *(f^1 \wedge *\tau))\right)
\right].\nonumber \\
\ea
\end{lem}
\Pf Applying (\ref{id-new}) with $\alpha=f^1$ 
we obtain
\ba
*(df^1 \wedge \sigma)=-df^1+*(d^7_7f^1 \wedge *\sigma)
+\f{1}{3}\xi\left(*(*\sigma \wedge *(\alpha \wedge *\tau_{\sigma}))\right).
\ea
Combinig this with (\ref{HH}) we then arrive at 
(\ref{Hequation}). \qed \\

\sect{Additional differential identities}

In this section we establish several differential identities  which will be used in Section 7 for proving  the uniqueness  of 
the solution of the Laplacian flow with given initial data.  As in the last section, the proofs of these identities determine the precise forms of the additional lower order terms which arise in the situation of a general closed 
 $G_2$-structure in comparison with a torsion-free 
 $G_2$-structure.  
 
 Let a closed 
$G_2$-structure $\sigma$ on $M$ be given. As in the last section, we drop the subscript $\sigma$ in the notations. 

\subsection{Two First Order Identities} 
\hspace{1cm}
\vspace{2mm}

\noindent {\bf Definition 4.1} The differential 
$d^1_7: \Omega^0(M) \rightarrow \Omega^1(M)$ is defined to be 
\ba
d^1_7f=df.
\ea
The differential $d^7_1: \Omega^1(M) \rightarrow \Omega^0(M)$ is defined to be the former $L^2$-adjoint of $d^1_7$, thus 
\ba
d^7_1 \alpha=d^*\alpha=-*d*\alpha.
\ea
The differential $d^7_{27}: \Omega^1(M) \rightarrow 
\Omega^3_{27}(M)$ is defined to be
\ba
d^7_{27} \alpha=\pi^3_{27} d*(\alpha \wedge *\sigma).
\ea
The differential $d^{14}_7: \Omega^2_{14}(M) \rightarrow 
\Omega^1(M)$ is defined to be the formal 
$L^2$ adjoint of $d^7_{14}$, whose definition is given in the last section. Thus $d^{14}_7=(d^7_{14})^*$.
Finally, we define $d^{14}_{27}: \Omega^2_{14}(M) 
\rightarrow \Omega^3_{27}(M)$ by the formula 
$d^{14}_{27} \beta=\pi^3_{27} d\beta.$ 
\vspace{2mm}

First we present a new differential identity. It has the remarkable feature of expressing the co-differential of a special kind of 2-form in terms of its differential. This is impossible for general 2-forms.

\begin{lem} \label{2-7} There holds for all $\beta \in \Omega^2_7(M)$
\ba \label{2form}
*d*\beta=\f{1}{2}\sigma \neg*d\beta-\frac{1}{2} \sigma \neg *(e^i \wedge *(\alpha 
\wedge *\nabla_{e_i} \sigma))-*(\alpha \wedge *\tau),
\ea
where $\alpha\in \Omega^1(M)$ is uniquely determined by the equation 
$\beta=*(\alpha \wedge *\sigma)$ (according to 
(\ref{decomp-ch})). 
In other words, there holds 
\ba \label{*d}
d^7_7 \alpha \equiv *d (\alpha \wedge *\sigma)=\f{1}{2} 
\sigma \neg *d*(\alpha \wedge *\sigma)
-\frac{1}{2} \sigma \neg *(e^i \wedge *(\alpha 
\wedge *\nabla_{e_i} \sigma))-*(\alpha \wedge *\tau)
\ea
for all 1-forms $\alpha$.
\end{lem}
\Pf   There holds
\ba
*d (\alpha \wedge *\sigma)&=&*(d\alpha \wedge 
*\sigma)-*(\alpha \wedge *\tau)\nonumber \\
&=&\Phi(\nabla \alpha)-*(\alpha \wedge *\tau),
\ea
where for $p\in M$ and $\Theta\in T^*_pM \otimes T^*_pM$
\ba
\Phi_p(\Theta)=*(e^i \wedge (e_i \lr \Theta) \wedge 
*\sigma).
\ea 
On the other hand, there holds
\ba
\sigma \neg *d*(\alpha \wedge *\sigma)
&=& \sigma \neg *(e^i \wedge 
*(\nabla_{e_i} \alpha \wedge *\sigma+
\alpha \wedge *\nabla_{e_i} \sigma)) 
\nonumber \\
&=& \Psi(\nabla \alpha)+\sigma \neg *(e^i \wedge 
*(\alpha \wedge *\nabla_{e_i} \sigma)),
\ea
where 
\ba
\Psi_p(\Theta)=\sigma \neg *(e^i \wedge 
*(e_i \lr \Theta \wedge *\sigma)) 
\ea 
with the above $\Theta$. For a fixed $p$, $\Phi_p$ and 
$\Psi_p$ are $G_2$-equivariant linear maps from 
$T^*_pM \otimes T^*_pM$  into 
$T_p^*M$. By the arguments in the proof of Lemma, 
$\Phi_p=\lambda \Psi_p$ for a scalar $\lambda$. 
Consider an induced orthonormal basis $e_i$ and its dual 
$e^i$. There hold 
\ba
\Phi_p(e^1 \otimes e^2)=*(e^1 \wedge e^2 \wedge *\sigma)
=*e^1 \wedge e^2 \wedge e^4 \wedge e^5 \wedge e^6 \wedge e^7=e^3
\ea
and 
\ba
\Psi_p(e^1 \otimes e^2)&=&\sigma \neg *(e^1 \wedge *(
e^2 \wedge *\sigma))\nonumber \\
&=&\sigma 
\neg *(e^1 \wedge *(e^2 \wedge e^4 \wedge e^5 \wedge e^6
\wedge e^7-e^1 \wedge e^2 \wedge e^3 \wedge e^5 \wedge e^7
+e^1 \wedge e^2 \wedge e^3 \wedge e^4 \wedge e^6))
\nonumber \\
&=&\sigma \neg (e^2 \wedge e^3 \wedge e^4 \wedge e^6
-e^2 \wedge e^3 \wedge e^5 \wedge e^7)
\nonumber \\
&=& 2 e^3. 
\ea
It follows that $\Phi_p=\f{1}{2} \Psi_p$, which leads 
to (\ref{*d}). \qed \\

The differential identity in the next lemma without the lower order terms can be found
in [B2] for the special case of a torsion-free $G_2$-structure.

\begin{lem} \label{d**lemma} There holds for all 
$\alpha\in \Omega^1(M)$
\ba \label{d**}
d*(\alpha \wedge *\sigma)=-\f{3}{7} (d^7_1 \alpha) 
\sigma-\f{1}{2} *(d^7_7 \alpha \wedge \sigma)
+d^7_{27} \alpha+\zeta(\alpha)
\ea
with
\ba
\zeta(\alpha)=\zeta_{\sigma}(\alpha)=-\f{1}{7} ((\alpha \wedge *\sigma) \cdot
*\tau) \sigma-\frac{1}{4} \sigma \neg *(e^i \wedge *(\alpha 
\wedge *\nabla_{e_i} \sigma))-\f{1}{2}*(\alpha \wedge *\tau).
\ea
\end{lem}
\Pf   First we decompose $d*(\alpha \wedge *\sigma)$ into irreducible parts
\ba \label{d**1}
d*(\alpha \wedge *\sigma)=\pi^3_1 d*(\alpha \wedge *\sigma)+\pi^3_7 d*(\alpha \wedge *\sigma)
+\pi^3_{27}d*(\alpha \wedge *\sigma).
\ea
The first part can be determined by employing orthogonal relations and integration by parts. (One can also argue as in the proof of Lemma \ref{idlemma} and Lemma.  By (\ref{decomp-ch}) it can be written as 
$f\sigma$ for a scalar function $f$.  Taking the $L^2$ inner product of (\ref{d**1}) with $\tilde f \sigma$ for an arbitrary scalar function $\tilde f$ we infer
\ba
7\int_M f \tilde f&=& \int_M d*(\alpha \wedge *\sigma)
\cdot 
\tilde f \sigma\nonumber \\
&=&-\int_M \alpha \wedge *\sigma \cdot   
d (\tilde f *\sigma \nonumber) \\
&=& -\int_M  (\alpha \wedge *\sigma) \cdot (d \tilde f) 
\wedge *\sigma-\int_M \tilde f (\alpha \wedge *\sigma) \cdot
*\tau. 
\ea 
Appealing to Lemma \ref{productlemma} we then deduce
\ba
7\int_M f \tilde f&=& 
-3\int_M  \alpha  \cdot d \tilde f 
-\int_M \tilde f (\alpha \wedge *\sigma) \cdot
*\tau \nonumber \\
&=& -3 \int_M \tilde f d^7_1 \alpha-
\int_M \tilde f (\alpha \wedge *\sigma) \cdot
*\tau.
\ea 
We conclude that $f=-\f{3}{7}d^7_1 \alpha
-\f{1}{7} (\alpha \wedge *\sigma) \cdot
*\tau$ and hence 
\ba \label{3-1}
\pi^3_1d*(\alpha \wedge *\sigma)=-\f{3}{7}(d^7_1 \alpha)\sigma-\f{1}{7} ((\alpha \wedge *\sigma) \cdot
*\tau) \sigma.
\ea

Next we determine the second term in the decomposition (\ref{d**1}).  By Lemma \ref{3-decomp} and Lemma we 
deduce 
\ba \label{3-7}
\pi^3_7 d*(\alpha \wedge *\sigma)&=&-\f{1}{4} *((\sigma \neg *d*(\alpha \wedge *\sigma))\wedge \sigma)))
\nonumber \\
&=& -\f{1}{2}(d^7_7 \alpha \wedge \sigma)-\frac{1}{4} \sigma \neg *(e^i \wedge *(\alpha 
\wedge *\nabla_{e_i} \sigma))-\f{1}{2}*(\alpha \wedge *\tau).
\ea
Combining (\ref{3-1}) and (\ref{3-7}) we arrive at 
(\ref{d**}). \qed

\subsection{An Identity for the Hodge Laplacian on 1-Forms} 
\hspace{1cm}
\vspace{1mm}

The following second order differential identity 
without the lower order term can be found in [B2] 
for the special case of a torison-free $G_2$-structure.

\begin{lem} \label{hodge-1form} There holds for all 
$\alpha \in \Omega^1(M)$
\ba \label{Delta}
\Delta \alpha=(d^1_7d^7_1+d^7_7d^7_7) \alpha
+\f{1}{3}*d*\xi\left(*(*\sigma \wedge *(\alpha \wedge *\tau))\right), 
\ea
where $\xi$ is given in (\ref{T}).
\end{lem}
\Pf There holds $\Delta \alpha=dd^* \alpha+d^* d\alpha
=dd^* \alpha+*d*\alpha$.
By the definitions of $d^1_7$ and $d^7_1$ we have
\ba
dd^* \alpha=d^1_7d^7_1 \alpha.
\ea
 On the other hand, 
we have by Lemma \ref{newdiff}
\ba
*d*d\alpha=*d(d^7_7 \alpha \wedge *\sigma)
-*d(d\alpha \wedge \sigma)+\f{1}{3}*d*\xi\left(*(*\sigma \wedge *(\alpha \wedge *\tau_{\sigma}))\right).
\ea
By the definition of $d^7_7$ the first term on the 
above right hand side is precisely $d^7_7d^7_7\alpha$. 
The second term vanishes because $\sigma$ is closed. 
Combining the above calculations we arrive at 
(\ref{Delta}). \qed \\

\sect{Linear Parabolic Theory}
\vspace{2mm}

As mentioned in Introduction, the gauge fixed Laplacian flow is parabolic only in the direction of closed forms. Hence there are troubles with applying the conventional theory of parabolic equations. For this reason, an approach in terms of  Nash-Moser implicit function theorem was adopted in [BX]. In this section, we develop a new linear parabolic theory for closed forms, which will enable us to construct short time solutions of  the gauge fixed Laplacian flow via the classical implicit function theorem. 

For the sake of completeness, we also include the corresponding theory for exact forms. (We also have the corresponding theories for co-closed and co-exact forms.) There is a subtlety here as mentioned in  Introduction. The structure of the Laplacian equation or the gauge-fixed Laplacian equation for closed forms allows one to treat them as equations for exact forms, namely one can assume that $\sigma-\sigma_0$ is exact with
$\sigma_0$ denoting the initial $G_2$-structure. 
However, the linear parabolic theory for exact forms is not suitable 
for treating the issue of short time solutions due to the lack of completeness of the involved function spaces
$C^{l,l/2}_d(\pi^*(\Lambda^3_d T^*M))$, cf. the discussions below. 

The parabolic H\"{o}lder spaces, namely the 
${\mc C}^{l, l/2}$ spaces, play an important role in this paper both for 
handling  short time solutions and the convergence of 
the Laplacian flow. These function spaces are used 
e.g. in the classical text [LSU].  They were first 
introduced in a geometric set-up in [Y3]. Alternatively, we can also use 
parabolic Sobolev spaces to handle short time solutions.
But the use of the ${\mc C}^{l,l/2}$ spaces is crucial for 
proving convergence, see Section 5.  The definition of the ${\mc C}^{l,l/2}$ spaces is given in Appendix.

In this section, $M$ stands for a compact smooth manifold of dimension $n\ge 2$.  Fix a background  metric $g_*$ on $M$, which is used to make various measurements. It is required to have enough smoothness 
in each context. Throughout this section, {\it all norms are measured w.r.t.~to $g_*$}, unless otherwise indicated. 
Note that we can choose $g_*$ according to our needs in 
each situation.  
For example, we can choose $g_*$ to be the induced metric of a given torsion-free $G_2$-structure in the context of Theorem \ref{stability}. 
We can translate easily the measurements w.r.t.~one background metric  into measurements w.r.t.~another backgroup metric. \\

\subsection{Linear parabolic theory for general forms}
\hspace{1cm}\vspace{1mm}

Consider the vector bundle  $E=\Lambda^jT^*M,
1\le j \le 7$.  We'll fix $j$ in the discussions below. Let $C^{l}(E)$ denote the space of $C^{l}$ sections of 
$E$, equipped with the $C^l$-norm, which is defined w.r.t. $g_*$.
Fix $T>0$.  let $\pi=\pi_{[0, T]}: M \times [0, T] \rightarrow M$ denote the projection $\pi_{[0, T]}(p, t)=p$.  Let ${\mc C}^{l, l/2}(\pi^*E)$
denote the ${\mc C}^{l, l/2}$ sections of $\pi^*E$, equipped 
with the ${\mc C}^{l,l/2}$-norm, which is defined w.r.t. 
$g_*$. Note that a section $\gamma$ of $\pi^* E$ has arguments $(p, t) \in M \times [0, T]$ and 
satisfies $\gamma(p, t)\in E_p$. 

Let $l>2$ and $\mathcal U$ an open subset of ${\mc C}^{l, l/2}(\pi^*E)$.  To each operator
$F: {\mathcal U} \rightarrow {\mc C}^{l-2, (l-2)/2}(\pi^*E)$ we associate its  
$P$-operator 
\ba
P_F: {\mathcal U} \rightarrow 
{\mc C}^{l-2, (l-2)/2}(\pi^*E)
\ea 
given by $P_F=\f{\p }{\p t}+F$ and its $P$-map
\ba
{\mathcal P}_F: {\mathcal U} \rightarrow 
{\mc C}^{l-2, (l-2)/2}(\pi^*E) \times C^{l}(E) 
\ea
given by ${\mathcal P}_F(\gamma)=(\f{\p \gamma}{\p t}+F(\gamma), 
\gamma(0)).$

\begin{theo} \label{full} Let $g_0$ be a $C^l$ metric on $M$ for a given noninteger $l>2$, and 
$\Delta$ its Hodge Laplacian. Let $\pi=\pi_{[0,T]}$ for a given $T>0$. There is a positive constant $\delta_0=\delta_0(\|g_0\|_{C^0}, \|g_0^{-1}\|_{C^0}, l, g_*)$ with the following properties. Let $\Phi_0 \in {\mc C}^{l-1, (l-1)/2}(Hom(\pi^*(\Lambda^{j} T^*M, \Lambda^{j-1} T^*M))$ and 
\ba
\Phi_1 \in {\mc C}^{l-1, (l-1)/2}(T^*M \otimes \Lambda^j T^*M, \Lambda^{j-1} T^*M). \nonumber
\ea
 Set 
$\Phi(\gamma)=\Phi_0(\gamma)+\Phi_1(\nabla \gamma)$. 
Assume 
\ba \label{phi-small}
\|\Phi_1\|_{C^0} \le \delta_0.
\ea
Then the $P$-map of the operator $\Delta+d \circ \Phi$
\ba
{\mathcal P}_{\Delta+d \circ \Phi}: 
{\mc C}^{l,l/2}(\pi^*\Lambda^j T^*M) \rightarrow  
{\mc C}^{l-2,(l-2)/2}(\pi^*\Lambda^j T^*M) \times C^{l}(\Lambda^j T^*M))
\ea
is an isomorphism. Moreover, there hold
\ba \label{norm}
\|{\mathcal P}_{\Delta+d \circ \Phi}\| \le C
\text{ and } \|{\mathcal P}^{-1}_{\Delta+d \circ \Phi}\| \le C
\ea
for a positive constant $C=C(n,l,T,\|g_0\|_{C^{l-1}},
\|g^{-1}_0\|_{C^0},  \|\Phi_0\|_{{\mc C}^{l-1, (l-1)/2}}, 
\|\Phi_1\|_{{\mc C}^{l-1, (l-1)/2}},$
$ g_*)$. 

The number $\delta_0$ depends on each involved scalar quantity decreasingly, while the number $C$ has increasing dependences. The dependences of $\delta_0$ and $C$ on $g_*$ are in terms of its Riemannian norm $|g_*|_{C^{l-1}}$ (see [Y3] for the definition of this norm).  The  dependences of constants  on $g_*$ below are all of the same nature. 
\end{theo}
 
\noindent {\it Proof of Theorem \ref{full}}  We have the following Bochner-Weitzenb\"{o}ck  formula
\ba \label{bochner}
\Delta=\nabla^* \nabla+{\mathcal R},
\ea
where ${\mathcal R}={\mathcal R}_j$ is a linear 
action of the curvature operator of $g_0$ on $j$-forms. In a local 
chart, the leading term of the operator $\nabla^* \nabla$
takes the form $-\sum_{ij} g^{ij} \partial_i \partial_j$.
Hence the parabolic theory in [LSU] can be applied, and 
the desired isomorphism property and estimates follow,
see [Y2] or [Y3] for details. Note that the smallness condition (\ref{phi-small}) is for the purpose of obtaining uniform strong ellipticity of the operator $-\Delta-d \circ \Phi$.\qed \\

Theorem \ref{full} will be applied below to establish 
a linear parabolic theory for closed forms. On the other hand, we have the following time-interior version of 
Theorem \ref{full}, which will be used in Section 7 for 
handling long time existence and convergence of the Laplacian flow. 

\begin{lem} \label{interiorlemma} Assume the same set-up as in Theorem \ref{full}. Moreover, assume (\ref{phi-small}).
Let $\gamma\in C^{l, l/2}(\pi^*E)$ and 
$\alpha \in C^{l-2, (l-2)/2}(\pi^*E)$ satisfy
\ba \label{psiequation}
\f{\p \gamma}{\p t}+\Delta \gamma+d(\Phi(\gamma))=\alpha
\ea
on $M \times [0, T]$. Let $0<\epsilon_1<\epsilon_2<T$. Then there is a positive constant $C=C(l, T, (\epsilon_2-\epsilon_1)^{-1})$ depending only on $l, T$ and $(\epsilon_2-\epsilon_1)^{-1}$ such that 
\ba \label{reduction}
\|\gamma\|_{{\mathcal{C}}^{l, l/2}(M\times [\epsilon_2, T])} &\le& C \cdot C(l, T, (\epsilon_2-\epsilon_1)^{-1})
[\|\alpha\|_{{\mathcal{C}}^{l-2, (l-2)/2}(M\times [\epsilon_1, T])}\nonumber \\
&&+(\epsilon_2-\epsilon_1)^{-1}\|\gamma\|_{{\mathcal{C}}^{l-2, (l-2)/2}(M\times [\epsilon_1, T])}],
\ea 
where $C$ is the constant in (\ref{norm}). 
\end{lem}
\Pf Fix a nonnegative smooth function $\eta$ on 
$\mathbb{R}$ such that $\eta(t)=0$ for $t\le 0$ and 
$\eta(t)=1$ for $t\ge 1$. Then we set $\eta_{\epsilon_1, 
\epsilon_2}(t)=\eta((\epsilon_2-\epsilon_1)^{-1} (t-\epsilon_1))$ and $\tilde \gamma=\eta_{\epsilon_1, 
\epsilon_2}(t) \gamma$. There holds
\ba
\f{\p \tilde \gamma}{\p t}+\Delta \tilde \gamma+d(\Phi(\tilde \gamma))=\eta_{\epsilon_1, \epsilon_2}  \alpha-\eta_{\epsilon_1, \epsilon_2}'
\gamma
\ea
on $M \times [0, T]$.
Let $\delta=l-[l]$. Then there hold
\ba
\|\eta_{\epsilon_1, \epsilon_2}\|_{{\mathcal{C}}^{l-2, (l-2)/2}(M\times [0, T])}
\le C(l, T, (\epsilon_2-\epsilon_1)^{-1})
\ea
and 
\ba
\|\eta_{\epsilon}'\|_{{\mathcal{C}}^{l-2, (l-2)/2}(M\times [0, T])}\le 
(\epsilon_2-\epsilon_1)^{-1}C(l, T, (\epsilon_2-\epsilon_1)^{-1}),
\ea
where
\ba
C(l, T, x)=C_{[l]}(\max\{T^{1-\delta}, T^{1-\delta/2}, T^{(1-\delta)/2}\}
x^{[\f{l}{2}]+1}+ \max_{0\le j\le [\f{l}{2}]} x^{j})
\ea
for a positive constant $C_{[l]}$ depending only on $[l]$. 
Then it follows that
\ba
\|\eta_{\epsilon_1, \epsilon_2} \alpha\|_{{\mathcal{C}}^{l-2, (l-2)/2}(M\times [0, T])}
\le C(l, T, (\epsilon_2-\epsilon_1)^{-1}) \|\alpha\|_{{\mathcal{C}}^{l-2, (l-2)/2}(M\times [\epsilon_1, T])}
\ea
and
\ba
\|\eta_{\epsilon_1, \epsilon_2}'\gamma\|_{{\mathcal{C}}^{l-2, (l-2)/2}(M\times [0, T])} \le
(\epsilon_2-\epsilon_1)^{-1}C(l, T, (\epsilon_2-\epsilon_1)^{-1})\|\gamma\|_{{\mathcal{C}}^{l-2, (l-2)/2}(M\times [\epsilon_1, T])}.
\ea
Applying Theorem \ref{full} we then arrive at
\ba
\|\tilde \gamma\|_{{\mathcal{C}}^{l, l/2}(M\times [0, T])} &\le& C\cdot C(l, T, (\epsilon_2-\epsilon_1)^{-1})
(\|\alpha\|_{{\mathcal{C}}^{l-2, (l-2)/2}(M\times [\epsilon_1, T])} \nonumber \\
&&+(\epsilon_2-\epsilon_1)^{-1}\|\gamma\|_{{\mathcal{C}}^{l-2, (l-2)/2}(M\times [\epsilon_1, T])}),
\ea  
which implies (\ref{reduction}). \qed 

\begin{theo} \label{interiortheorem} Assume the same set-up as in Theorem \ref{full}. Moreover, assume (\ref{phi-small}).
Let $\gamma\in C^{l, l/2}(\pi^*E)$ and 
$\alpha \in C^{l-2, (l-2)/2}(\pi^*E)$ satisfy
(\ref{psiequation}) on $M \times [0, T]$. Let $0<\epsilon<T$. Then there is a positive constant $C(l, T, \epsilon^{-1},C)$ depending only on $l, T, \epsilon^{-1}$ and the $C$ in (\ref{norm}) such that 
\ba \label{interiorestimate}
\|\gamma\|_{{\mathcal{C}}^{l, l/2}(M\times [\epsilon, T])} \le C(l, T, \epsilon^{-1},C)
(\|\alpha\|_{{\mathcal{C}}^{l-2, (l-2)/2}(M\times [0, T])}+\|\gamma\|_{{\mathcal{C}}^{m, m/2}(M\times [0, T])}),
\ea 
where $m=l-2k\ge 0$ for the largest nonnegative integer 
$k$. (If $l=2k+\mu$ for $0<\mu<1$, then $m=\mu$.
If $l=2k+1+\mu$ for $0<\mu<1$, then $m=1+\mu$.) \end{theo}
\Pf Apply Lemma \ref{interiorlemma} successively to 
the sequence of pairs $(\epsilon/2, \epsilon), (\epsilon/4, 
\epsilon/2),...$ (playing the role of $(\epsilon_1, \epsilon_2)$), 
with a sequence of decreasing $l$, i.e. $l,l-2,...$.
After finitely many steps we then arrive at 
(\ref{interiorestimate}). \qed \\

\subsection{Linear parabolic theory for closed forms}
\hspace{1cm}\vspace{1mm}

We set for $l\ge 0$
\ba
&&C^{l}_{o}(\Lambda^j T^*M)=\{
\gamma \in C^l(\Lambda^j T^*M): d\gamma=0\}.
\ea
Here the equation $d\gamma=0$ is in the sense of 
distribution in the case $0\le l <1$, i.e. 
\ba
<\gamma, d^*_{g_*}\theta>_{L^2_{g_*}}=0
\ea
for all $\theta\in \Omega^{j+1}(M)$, where $d^*_{g_*}$ is the co-differential associated with $g_*$, and 
$<,>_{L^2_{g_*}}$ denotes the $L^2$ inner product w.r.t.~$g_*$.
(We can also replace $g_*$ by a given $g_0$ as in 
Theorem \ref{full}.)   Obviously, $C^{l}_{o}(\Lambda^j T^*M)$ is a closed subspace, and hence  a Banach subspace of $C^l(\Lambda^j T^*M)$. 
For a noninteger $l>0$ we set
\ba
&&C^{l,l/2}_{o}(\pi^*\Lambda^j T^*M)
=\{\gamma \in C^{l,l/2}(\pi^*\Lambda^j T^*M):
d\gamma(\cdot, t)=0 \mbox{ for each } 
t \in [0, T]\} \nonumber \\
\ea
which is obviously a closed subspace, and hence a Banach subspace of 
$C^{l,l/2}(\pi^*\Lambda^j T^*M)$.
(Again, the equation $d \gamma(\cdot, t)=0$ is in the 
sense of distribution in the case $0<l<1$.) \\

\begin{theo} \label{closed} Let $l>2$ be a noninteger. Let $g_0, \Delta$, $\Phi_0, \Phi_1, \Phi$ and  $\delta_0$ be the same as in Theoerem \ref{full}.  
Assume (\ref{phi-small}).
Then the $P$-map of the operator $\Delta+d \circ \Phi$
\ba
{\mathcal P}_{\Delta+d \circ \Phi}: 
{\mc C}^{l,l/2}_o(\pi^*\Lambda^j T^*M) \rightarrow  
{\mc C}^{l-2,(l-2)/2}_o(\pi^*\Lambda^j T^*M) \times C^{l}_o(\Lambda^j T^*M)
\ea
is an isomorphism. Moreover, there hold
\ba \label{closednorm}
\|{\mathcal P}_{\Delta+d \circ \Phi}\| \le C
\text{ and } \|{\mathcal P}^{-1}_{\Delta+d \circ \Phi}\| \le C
\ea
for a positive constant $C=C(n,l,T,\|g_0\|_{C^{l-1}},
\|g^{-1}_0\|_{C^0},  \|\Phi_0\|_{{\mc C}^{l-1, (l-1)/2}}, 
\|\Phi_1\|_{{\mc C}^{l-1, (l-1)/2}},$
$ g_*)$. 
 Thus, for each $\beta \in  C^{l}_o(\Lambda^j T^*M)$ and 
$\alpha \in  
{\mc C}^{l-2,(l-2)/2}_o(\pi^*\Lambda^j T^*M)$ there is a unique solution $\gamma \in {\mc C}^{l,l/2}_o(\pi^*\Lambda^j$
$ T^*M)$
of the initial value problem
\ba \label{closedheat}
\f{\p \gamma}{\p t}+\Delta \gamma+d(\Phi(\gamma))=\alpha \text{ (on } [0, T]),
\ea
\ba    
\gamma(\cdot, 0)=\beta,
\ea
such that 
\ba \label{backforward}
\|\gamma\|_{{\mc C}^{l,l/2}}\le C(\|\alpha\|_{{\mc C}^{l-2, (l-2)/2}}+\|\beta\|_{C^l}).
\ea
\end{theo}
\Pf  We have 
\ba
d \circ \Phi(\gamma)=\sum_i e^i \wedge [(\nabla_{e_i} 
\Phi_0)(\gamma)+\Phi_0(\nabla_{e_i} \gamma)
+(\nabla_{e_i} \Phi_1)(\nabla \gamma)
+\Phi_1(\nabla_{e_i} \nabla \gamma)].
\ea
Applying Theorem \ref{full}, we infer that the extended $P$-map
\ba
{\mathcal P}_{\Delta+d \circ \Phi}: 
{\mc C}^{l,l/2}(\pi^*\Lambda^j T^*M) \rightarrow 
{\mc C}^{l-2, (l-2)/2}(\pi^*\Lambda^j M) \times 
C^{l}(\Lambda^j T^*M)
\ea
is an isomorphism and satisfies the estimate (\ref{norm}). 
Hence it suffices to 
show 
\ba  \label{inclusion}
{\mc C}^{l,l/2}_o(\pi^*\Lambda^j T^*M)={\mathcal P}_{\Delta+d \circ \Phi}^{-1}({\mc C}^{l-2, (l-2)/2}_o(\pi^*\Lambda^j M)) \times 
C^{l}_o(\Lambda^j T^*M)).
\ea
First we show that the LHS of (\ref{inclusion}) is contained in the RHS of (\ref{inclusion}). It suffices to show that $\f{\p \gamma}{\p t}+\Delta \gamma
+d\Phi(\gamma)$ is closed in the sense of distribution
for each $\gamma \in {\mc C}^{l,l/2}_o(\pi^*\Lambda^j T^*M)$. For such a $\gamma$ and an arbitrary 
$\theta \in \Omega^{j+1}(M)$ we indeed have 
\ba
&&<\f{\p \gamma}{\p t}+\Delta \gamma
+d\Phi(\gamma), d^*_{g_*}\theta>_{L^2_{g_*}}=\f{\p }{\p t}
<\gamma, d^*_{g_*}\theta>_{L^2_{g_*}}+<dd^* \gamma+d\Phi(\gamma), d^*_{g_*}\theta>_{L^2_{g_*}}
\nonumber \\&& \hspace{5.5cm}=0.
\ea
 To show the opposite inclusion, consider $\gamma={\mathcal P}_{\Delta+d \circ \Phi}^{-1}(\alpha, \beta)$ for some $\beta \in
C^{l}_o(\Lambda^j T^*M)$ and $\alpha \in {\mc C}^{l-2, (l-2)/2}_o(\pi^*\Lambda^j M)$. Thus $\gamma, \alpha$ and $\beta$ satisfy the equation (\ref{closedheat}). Assume $l>3$. Taking the differential in the equation we deduce 
\ba
\begin{cases}
\f{\p}{\p t} d \gamma+dd^*d \gamma=0, \\
d\gamma(\cdot, 0)=0.
\end{cases}
\ea
But $dd^*d \gamma=\Delta d\gamma$. Hence we infer 
$d\gamma \equiv 0$. Indeed, this can be shown directly as follows.  Multiplying the above equation by $d\gamma$ and then integrating (first in space, then in time) lead to 
\ba
\|d\gamma(\cdot, t)\|_{L^2}^2+\int_0^t \|d^*d \gamma(\cdot, t)\|_{L^2}^2=\|d\gamma(\cdot, 0)\|_{L^2}=0,
\ea
where the $L^2$-norms are w.r.t.~$g_0$.
It follows that $\gamma \in {\mc C}^{l,l/2}_o(\pi^*\Lambda^j T^*M)$. 

The case $2<l<3$ requirs a different argument which also applies to the case $l>3$. Choose a complete set of $L^2$-orthonormal 
eigenforms $\gamma_k$  of degree $j$ for $\Delta$, such that each $\gamma_k$ is either harmonic, exact, or 
coexact. This is possible for the following reason.
Let $\phi$ be an eigenform with nonzero eigenvalue 
$\lambda$.  We write $\phi=h+d\psi+d^*\chi$, where 
$h$ is harmonic. There holds
\ba
dd^*d\psi+d^*dd^*\chi=\lambda h+\lambda d\psi
+\lambda d^* \chi.
\ea
It follows that
\ba
dd^*d\psi=\lambda d\psi, d^*d d^*\chi=\lambda 
d^*\chi, \lambda h=0.
\ea
Hence $\Delta d\psi=\lambda d\psi, \Delta d^*\chi
=\lambda d^* \chi.$  

Let $\phi_i$ be the exact 
forms among the $\gamma_k$ with $\Delta \phi_i=\lambda_i \phi_i$.  We have 
\ba
d\gamma(\cdot, t)=\sum_i a_i(t) \phi_i.
\ea
Multiplying the equation (\ref{closedheat}) with 
$d^*\phi_i$ and integrating lead to 
\ba
&&\f{d a_i}{dt}=\f{d}{dt}<d\gamma, \phi_i>=<\f{\p \gamma}{\p t}, d^* \phi_i>
=-<\Delta \gamma+d\Phi(\gamma), d^* \phi_i> \nonumber 
\\&&
=-<d^*d \gamma, d^* \phi_i>=-<d\gamma, dd^* \phi_i>
=-\lambda_i a_i.
\ea
Since $\lambda_i>0$ and $a_i(\cdot, 0)=0$, we infer 
$a_i=0$. Consequently, $\gamma \in  {\mc C}^{l,l/2}_o(\pi^*\Lambda^j T^*M)$.
\qed

\subsection{Linear parabolic theory for exact forms}
\hspace{1cm}\vspace{1mm}

We set for $l\ge 1$
\ba
&&{C}^{l}_{d}(\Lambda^j T^*M)=d({ C}^{l+1}(\Lambda^{j-1} T^*M)).
\ea
For a noninteger $l \ge 1$ we set
\ba
&&{\mc C}^{l,l/2}_{d}(\pi^*\Lambda^j T^*M)
=d({\mc C}^{l+1,(l+1)/2}(\pi^*\Lambda^j T^*M)).
\ea

Employing basic linear elliptic estimates one can easily show that ${C}^{l}_{d}(\Lambda^j T^*M)$ is a closed and hence Banach supspace of ${C}^{l}(\Lambda^{j} T^*M)$.
However, as it turns out, ${\mc C}^{l,l/2}_{d}(\pi^*\Lambda^j T^*M)$
is not a closed supspace of ${\mc C}^{l,l/2}(\pi^*\Lambda^j T^*M)$, and hence it is not a 
Banach space. The analytic reason for this is the lack of involvement of the time derivative
in its definition.

As a consequence of Theorem \ref{closed} we obtain the following result for exact forms.

\begin{theo} \label{exact} Let $l>2$ be a noninteger. Let $g_0, \Delta, \Phi_0, \Phi_1, \Phi$ and  $\delta_0$ be the same as in Theoerem \ref{full}.  
Assume (\ref{phi-small}). Then the parabolic map 
\ba
{\mc P}={\mathcal P}_{\Delta+d \circ \Phi}: {\mc C}^{l,l/2}_d(\pi^*\Lambda^j T^*M) \rightarrow 
{\mc C}^{l-2, (l-2)/2}_d(\pi^*\Lambda^j M) \times 
{\mc C}^{l}_d(\Lambda^j T^*M)
\ea
is an isomorphism. Moreover, the estimate 
(\ref{closednorm}) holds true with the same 
$C$. 
\end{theo}
\Pf As in the proof of Theorem \ref{closed}, it suffices to show 
\ba  \label{exactinclusion}
{\mc C}^{l,l/2}_d(\pi^*\Lambda^j T^*M)={\mathcal P}_{\Delta+d \circ \Phi}^{-1}({\mc C}^{l-2, (l-2)/2}_d(\pi^*\Lambda^j M) \times 
C^{l}_d(\Lambda^j T^*M))
\ea
for the extended $P$-map.  It is easy to see that the 
LHS of (\ref{exactinclusion}) is contained in its 
RHS. On the other hand, if $\alpha$ and $\beta$ are exact, integrating the equation (\ref{closedheat}) in time shows that $\gamma$ is also exact. Hence the RHS of 
(\ref{exactinclusion}) is also contained in the LHS of 
(\ref{exactinclusion}). \qed \\

\noindent {\bf Remark} Obviously, the analog of 
Theorem \ref{closed} for co-closed forms and the analog of Theorem \ref{exact} 
for co-exact forms hold true if $d\circ \Phi$ is replaced by $d^* \circ \Phi$. \\

\sect{Short time solutions of the gauge fixed Laplacian flow}
\vspace{2mm}

From now on $M$ stands for a compact 7-dimensional manifold which admits closed $G_2$ structures. As in the last section, we fix a background metric $g_*$ on $M$.
All function norms in this section are associated with 
$g_*$. But pointwise norms and other geometric quantities are associated with an initial $G_2$-structure $\sigma_0$ in some situations. This will be made clear in the discussions below.  

\subsection{Gauge fixing}
\hspace{1cm}\vspace{1mm}

To construct short time smooth solutions of the Laplacian flow, we employ as in [BX] the following DeTurck type gauge fixing of the Laplacian flow 
\ba \label{gaugefixing}
\f{\p \sigma}{\p t}=\Delta_{\sigma} \sigma +\mc{L}_{X(\sigma)} \sigma,
\ea
where $X(\sigma)$ is a vector field associated with $\sigma$ and $\mc{L}_{X(\sigma)}$ denotes the Lie derivative.  The game of this gauge fixing is to find a
suitable $X(\sigma)$ such that the operator $\Delta_{\sigma} \sigma +\mc{L}_{X(\sigma)} \sigma$ has maximal (strong) ellipticity. In [BX], a vector field is 
constructed from the induced metric and its Levi-Civita connection. Based on the new differential identities in 
Section 3,  
we introduce a new vector field which has a more transparent structure. 

Let a reference closed $G_2$-structure $\sigma_0$ be 
given. We set 
$\theta=\sigma-\sigma_0$ for a $G_2$-structure $\sigma$
and write 
\ba \label{theta}
\theta=f^0 \sigma_0+*_{\sigma_0} (f^1 \wedge \sigma_0)+
f^3
\ea
as in (\ref{decomp}).  
 We define 
\ba \label{X}
X_{\sigma_0}(\theta)=(\f{7}{3}df^0+2(d^7_7)_{\sigma_0} f^1)_{\sharp_{\sigma_0}}.
\ea
It is motivated by the identity (\ref{hodge}). Obviously, $X_{\sigma_0}(\theta)$ is defined 
for an arbitrary 3-form $\theta$ given by 
(\ref{theta}). \\

\noindent {\bf Definition 6.1}  The {\it $\sigma_0$-gauged Laplacian flow} is defined to be 
\ba \label{gauged}
\f{\p \sigma}{\p t}=\Delta_{\sigma} \sigma
+\mc{L}_{X_{\sigma_0}(\sigma-\sigma_0)} \sigma.
\ea

For closed $\sigma$ we have $\mc{L}_{X_{\sigma_0}(\sigma-\sigma_0)} \sigma=
d(X_{\sigma_0}(\sigma-\sigma_0)\lr \sigma)$. Hence the $\sigma_0$-gauged Laplacian flow for closed $\sigma$ can be written as follows
\ba
\f{\p \sigma}{\p t}=\Delta_{\sigma} \sigma+d (X_{\sigma_0}(\sigma-\sigma_0)\lr \sigma).
\ea

Next we relate the operator $\Delta_{\sigma}\sigma+d(X_{\sigma_0}(\sigma-\sigma_0)\lr \sigma)=\Delta_{\sigma} \sigma
+d(X_{\sigma_0}(\theta) \lr \sigma)$ to the Hodge Laplacian of $\sigma_0$.  We shall adopt the following notations: we use $\mathcal D$ to denote the linearization, i.e.~the directional derivative, of an operator, and write it as $D$ in the case of a pointwise operator without involving partial detivatives (thus a finite dimensional operator). For example, 
$({\mathcal D}_{\sigma} \Delta_{\sigma}\sigma)(\theta)
=\f{d}{ds} \Delta_{\sigma+s\theta}(\sigma+s\theta)|_{s=0}$ and 
$(D_{\sigma} *_{\sigma} \sigma)(\theta)=
\f{d}{ds} *_{\sigma+s\theta}(\sigma+s\theta)|_{s=0}.$ 

\begin{lem} \label{transform} Let a closed $G_2$ structure $\sigma_0$ be given. There holds for an arbitrary closed 
$G_2$-structure $\sigma$
\ba \label{delta}
\Delta_{\sigma} \sigma+d(X_{\sigma_0}(\theta)\lr \sigma)=-\Delta_{\sigma_0} \theta
-d(\Phi_{\sigma_0}(\theta))
\ea
with 
\ba \label{phi-formula}
&&\Phi_{\sigma_0}(\theta)=
A(\sigma_0, \sigma_0+\theta, \theta, \nabla_{\sigma_0} \theta)
+B(\sigma_0, \sigma_0+\theta, \tau_{\sigma_0}, \theta)-\tau_{\sigma_0},
\ea
where $\theta=\sigma-\sigma_0$ as above, and $A$ and $B$ are smooth in their first two arguments and linear in the other two arguments.  The functions $A$ and $B$ are pointwise functions, e.g. $B(\sigma_0, \sigma, \tau_{\sigma_0}, \theta)(p)
=B(\sigma_0(p), \sigma(p), \tau_{\sigma_0}(p), \theta(p))$. They are also universal, i.e. their formulas are independent of the point and the manifold. 
In other words, these formulas are induced from 
the case of the Eulidean space in terms of an inducing map.
\end{lem}
\Pf Because $\sigma$ and $\sigma_0$ are closed, we have $\Delta_{\sigma} \sigma=-d*_{\sigma} d *_{\sigma} \sigma$ and $\Delta_{\sigma_0} \sigma_0=-d *_{\sigma_0} d *_{\sigma_0} \sigma_0$. Hence we obtain
\ba
\Delta_{\sigma} \sigma-\Delta_{\sigma_0}\sigma_0
=-d(*_{\sigma}d*_{\sigma} \sigma-*_{\sigma_0} d *_{\sigma_0} \sigma_0).
\ea
There holds 
\ba
&&*_{\sigma}d*_{\sigma} \sigma-*_{\sigma_0} d *_{\sigma_0} \sigma_0=*_{\sigma_0} d (*_{\sigma} \sigma-*_{\sigma_0} \sigma_0)+(*_{\sigma}-*_{\sigma_0})d (*_{\sigma} \sigma-*_{\sigma_0} \sigma_0)\no +(*_{\sigma}-*_{\sigma_0})d *_{\sigma_0}\sigma_0.
\ea
We have $D_{\sigma_0}(*_{\sigma}\sigma)(\theta)=
*_{\sigma_0}(\f{4}{3}f^0
\sigma_0+*_{\sigma_0}(f^1 \wedge \sigma_0)-f^3)$, cf.~[J]. 
It follows that 
\ba \label{joyce}
*_{\sigma} \sigma-*_{\sigma_0} \sigma_0=*_{\sigma_0}(\f{4}{3}f^0
\sigma_0+*_{\sigma_0}(f^1 \wedge \sigma_0)-f^3)+q(\sigma_0, \sigma, \theta, \theta),
\ea
where $q$ is given by 
\ba
q(\sigma_0, \sigma, \theta, \theta)
=\int_0^1 \int_0^1 t D^2 (*_{\sigma} \sigma)|_{\sigma_0+st(\sigma-\sigma_0)} (\theta, \theta) dsdt
\ea
with $D^2$ denoting the second derivative operator. Thus $q$ is smooth in its first two arguments and 
linear in the other two arguments. 
Note that 
$q$ is a universal pointwise function and involves  no derivative of $\sigma_0$ or $\sigma$.

Now we infer from the above formulas
\ba
\Delta_{\sigma}\sigma &=&\Delta_{\sigma_0} \sigma_0
-d(H_{\sigma_0} \theta)+d((*_{\sigma_0}-*_{\sigma})d (*_{\sigma} \sigma-*_{\sigma_0} \sigma_0) +(*_{\sigma_0}-*_{\sigma})d *_{\sigma_0}\sigma_0)
\no-d*_{\sigma_0}d(q(\sigma_0, \sigma, \theta, \theta)),
\ea
where $H_{\sigma_0} \theta=*_{\sigma_0} d *_{\sigma_0}
(\f{4}{3}f^0
\sigma_0+*_{\sigma_0}(f^1 \wedge \sigma_0)-f^3)$ is the operator introduced in (\ref{H}). Consequently,
\ba
\Delta_{\sigma}\sigma+d(X_{\sigma_0}(\theta)\lr \sigma)&=&\Delta_{\sigma_0} \sigma_0-d(H_{\sigma_0} \theta)+d(\frac{7}{3}(df^0)_{\sharp_{\sigma_0}} \lr  \sigma
+2(d^7_7f^1)_{\sharp_{\sigma_0}} \lr \sigma)
\no +d((*_{\sigma_0}-*_{\sigma})d (*_{\sigma} \sigma-*_{\sigma_0} \sigma_0) +(*_{\sigma_0}-*_{\sigma})d *_{\sigma_0}\sigma_0)
-d*_{\sigma_0}d(q(\sigma_0, \sigma, \theta, \theta)),
\no
\ea

Applying (\ref{hodge}) with $\sigma_0$ playing the role of $\sigma$ we then 
deduce (\ref{delta}) with 
\ba
\Phi_{\sigma_0}(\theta)&=&-\tau_{\sigma_0}+(*_{\sigma}-*_{\sigma_0})d(*_{\sigma}\sigma-*_{\sigma_0}
\sigma_0)+(*_{\sigma_0}-*_{\sigma})*_{\sigma_0} \tau_{\sigma_0}
\nonumber \\ && +*_{\sigma_0} d(q(\sigma_0, \sigma, \theta, \theta))-\f{7}{3}(df^0)_{\sharp_{\sigma_0}} \lr \theta-2(d^7_7f^1)_{\sharp_{\sigma_0}} \lr \theta
\nonumber \\&&+
\f{2}{3}{\mathcal \xi_{\sigma_0}}\left(*_{\sigma_0}(*_{\sigma_0}\sigma_0 \wedge *_{\sigma_0}(f^1 \wedge *_{\sigma_0}\tau_{\sigma_0}))\right)-\f{7}{3}f^0 \tau_{\sigma_0}.
\ea
It is easy to see that $(*_{\sigma}-*_{\sigma_0})d(*_{\sigma}\sigma-*_{\sigma_0}
\sigma_0)-\f{7}{3}(df^0)_{\sharp_{\sigma_0}} \lr \theta-2(d^7_7f^1)_{\sharp_{\sigma_0}} \lr \theta
$ can be written in the form $A_1(\sigma_0, \sigma, \theta,$ 
$\nabla_{\sigma_0} \theta)$, the expression
$(*_{\sigma_0}-*_{\sigma})*_{\sigma_0} \tau_{\sigma_0}+
\f{2}{3}{\mathcal \xi_{\sigma_0}}\left(*_{\sigma_0}(*_{\sigma_0}\sigma_0 \wedge *_{\sigma_0}(f^1 \wedge *_{\sigma_0}\tau_{\sigma_0}))\right)-\f{7}{3}f^0 \tau_{\sigma_0}$ can be written in the form 
$B_1(\sigma_0, \sigma, \tau_{\sigma_0}, \theta)$, and $*_{\sigma_0}d(q(\sigma_0, \sigma, \theta, \theta))$ can be written in the form 
$A_2(\sigma_0, \sigma, \theta, \nabla_{\sigma_0}\theta)
+B_2(\sigma_0, \sigma, \nabla_{\sigma_0} \sigma_0, \theta, \theta)$. (Note that $d=\sum_i e^i \wedge \nabla_{e_i}$.)  By Lemma \ref{Klemma}, $\nabla_{\sigma_0} \sigma_0$ can be expressed in terms of $\tau_{\sigma_0}$.  Hence  $B_2(\sigma_0, \sigma,
\nabla_{\sigma_0} \sigma_0, \theta, \theta)$ can be 
rewritten in the form $\tilde B_2(\sigma_0, \sigma, 
\tau_{\sigma_0}, \theta)$.  Setting $A=A_1+A_2$ and 
$B=B_1+\tilde B_2$ we then arrive at (\ref{phi-formula}). It is easy to 
verify that $A$ and $B$ have the claimed properties.  
Moreover, the quantities $*_{\sigma}$, $*_{\sigma_0}$,  
$A$ and $B$ etc.~can all be given explicitly in terms of 
$\sigma$ and $\sigma_0$. This is in part a consequence of (\ref{metric}) and (\ref{volume}).  \qed \\

By Lemma \ref{transform}, the $\sigma_0$-gauged Laplacian flow can be written as follows
\ba \label{newflow}
\f{\p \sigma}{\p t}=-\Delta_{\sigma_0} \theta
-d (\Phi_{\sigma_0}(\theta)).
\ea

For convenince of presentation, we formulate a simple 
lemma, which is an easy consequence of Lemma \ref{simple}
and the nature of the function $A$. 

\begin{lem} \label{Alemma} There hold 
\ba \label{Abound}
|A(\sigma_0, \sigma, \theta, \gamma)|\le \f{C_0}{\sqrt{7}} |\theta| \cdot |\gamma| \text{ and  hence }
\|A(\sigma_0, \sigma, \theta, \cdot)\|_{C^0} \le 
C_0\|\theta\|_{C^0} 
\ea
for a universal positive constant $C_0\ge 1$, provided that $\|\theta\|_{C^0}\le \epsilon_0$, where $\epsilon_0$ is from Lemma \ref{simple}. 
\end{lem} 

\vspace{1mm}

\subsection{Short time solutions of the gauge fixed 
Laplacian flow} 
\hspace{1cm}\vspace{1mm}

Let a closed $G_2$-structure $\sigma_0$ be given. The smoothness requirement for $\sigma_0$ will be specified 
in each situation (or is clear from the context).  In this subsection we prove an existence and uniqueness theorem for the $\sigma_0$-gauged Laplacian flow with 
$C^{4+\mu}$ initial data. For simplicity of presentation, here we choose the background metric $g_*$ 
of the last subsection to be the induced metric $g_{\sigma_0}$.  Thus, all norms in this subsection are 
measured w.r.t.~$g_{\sigma_0}$.  The covariant derivative $\nabla$ means $\nabla_{\sigma_0}$, i.e.~it is associated with $g_{\sigma_0}$. Furthermore, $e_i$ denotes a local orthonormal frame for $g_{\sigma_0}$, 
and $e^i$ its dual. \\

\noindent {\bf Definition 6.2} For the convenience of presentation we introduce the following notation
\ba
\nu_{\sigma, \sigma_0}=\Delta_{\sigma} \sigma
+d(X_{\sigma_0}(\sigma-\sigma_0)\lr \sigma)=-\Delta_{\sigma_0}\theta-d(\Phi_{\sigma_0}(\theta))
\ea
for a closed $G_2$-structure $\sigma$, where $\theta=\sigma-\sigma_0$ as before.\\

The following lemma provides an elementary estimate for 
this quantity.

\begin{lem} \label{nulemma} Assume $\|\sigma-\sigma_0\|_{C^0}
\le \epsilon_0$, where 
$\epsilon_0$ is from Lemma \ref{simple}.
There holds
\ba
\|\nu_{\sigma, \sigma_0}\|_{C^0} \le  \eta_1(\|\sigma-\sigma_0\|_{C^2}, \|\nabla \tau_{\sigma_0}\|_{C^0}, \|\tau_{\sigma_0}\|_{C^0}),
\ea
where  $\eta_1$ is a universal continuous 
increasing (in each argument) positive function of its 
arguments with $\eta_1(0, 0, \cdot)=0$.
We also have for $0<\mu<1$
\ba
[\nu_{\sigma, \sigma_0}]_{\mu} \le 
C\eta_2(\|\sigma-\sigma_0\|_{C^{2+\mu}}, [\nabla \tau_{\sigma_0}]_{\mu}, \|\tau_{\sigma_0}\|_{C^{\mu}}),
\ea
where $C$ is a positive constant depending only on
$\|\sigma_0\|_{C^{1+\mu}}$, and $\eta_2$ is a 
universal positive function with the same properties as $\eta_1$.
Moreover, we have 
\ba
\|\nu_{\sigma, \sigma_0}\|_{C^{2+\mu}}
\le \bar C \eta_3(\|\sigma-\sigma_0\|_{C^{4+\mu}},
\|\nabla \tau_{\sigma_0}\|_{C^{2+\mu}}, 
\|\tau_{\sigma_0}\|_{C^{2+\mu}}),
\ea
where $\bar C$ is a positive constant depending only on
$\|\sigma_0\|_{C^{3+\mu}}$, and $\eta_3$ is a universal positive function with the same properties as $\eta_1$.
\end{lem} 
\Pf Obviously, there hold
$|d\tau_{\sigma_0}|\le C_0|\nabla \tau_{\sigma_0}|$ and $|\Delta_{\sigma_0} \theta|
\le C_0|\nabla^2 \theta|$
for a universal positive constant $C_0$. 
On the other hand, we have for the functions 
$A$ and $B$ in the formula (\ref{phi-formula})
for $d(\Phi_{\sigma_0}(\theta))$ 
\ba \label{A-formula}
d(A(\sigma_0, \sigma, \theta, \nabla \theta))&=&\sum_i e^i \wedge [(\nabla_{e_i})_1 A(\sigma_0, \sigma, \theta, \nabla \theta)+(\nabla_{e_i})_2
A(\sigma_0, \sigma, \theta, \nabla \theta)\nonumber \\&&+A(\sigma_0, \sigma, \nabla_{e_i} \theta, \nabla \theta)+A(\sigma_0, 
\sigma, \theta, \nabla_{e_i} \nabla \theta)]
\ea
and 
\ba \label{B-formula}
d(B(\sigma_0, \sigma, \tau_{\sigma_0}, \theta))
&=&\sum_i e^i \wedge [(\nabla_{e_i})_1
B(\sigma_0, \sigma, \tau_{\sigma_0}, \theta)+(\nabla_{e_i})_2
B(\sigma_0, \sigma, \tau_{\sigma_0},  \theta)\nonumber \\&&+B(\sigma_0, \sigma, \nabla_{e_i} \tau_{\sigma_0}, \theta)+B(\sigma_0, 
\sigma, \tau_{\sigma_0}, \nabla_{e_i} \theta)],
\nonumber \\
\ea
where $(\nabla_{e_i})_k, k=1,2$ means taking the covariant derivative with the $k$-th argument as the variable, while keeping the other arguments parallel. 
Note that by Lemma \ref{Klemma}, $\nabla \sigma_0$ can be expressed in terms of $\tau_{\sigma_0}$ . It follows that
\ba \label{middleestimate}
|d(\Phi_{\sigma_0}(\theta))|
&\le& C[|\nabla \theta|(|\theta| |\nabla \theta|+|\tau_{\sigma_0}|
|\theta|)+|\nabla \theta|^2+|\theta| 
|\nabla^2 \theta|
\nonumber \\&& +
|\theta|(|\tau_{\sigma_0}|^2+|\nabla \tau_{\sigma_0}|) 
+|\tau_{\sigma_0}| |\nabla\theta|_{\sigma_0}+|\nabla \tau_{\sigma_0}|],
\ea
where the positive constant $C$ depends only on  $\|\sigma\|_{C^0}$
and $\|g_{\sigma}^{-1}\|_{C^0}$, which can be estimated by using 
Lemma \ref{simple} and the assumption $\|\sigma-\sigma_0\|_{C^0}\le \epsilon_0$. 
Obviously, the first claim of the lemma follows from 
(\ref{middleestimate}). The second and third claims of the lemma follow from similar computations based on (\ref{A-formula}) and (\ref{B-formula}).\qed \\

\begin{theo} \label{gauge1} (existence of the $\sigma_0$-gauged Laplacian flow) Assume $\sigma_0 \in C^{2+\mu}$ for some $0<\mu<1$.  Let $\hat \delta_0=\hat \delta_0(\sigma_0)=
\hat \delta_0(2+\mu, g_{\sigma_0})$ be from Lemma \ref{linearlemma} 
below (for $l=2+\mu$), which depends only on $g_{\sigma_0}$ (in terms of its Riemannian norm
$|g_{\sigma_0}|_{C^{1+\mu}}$). 
 Let 
$\sigma_1$ be a closed $C^{2+\mu}$ 
$G_2$-structure on $M$ such that $\nu_{\sigma_0, \sigma_1} \in C^{2+\mu}$ and 
\ba \label{bound1}
\|\sigma_1-\sigma_0\|_{C^0}\le \f{1}{4}\min\{\epsilon_0, \hat \delta_0\}.
\ea
For each positive constant $K>0$ there is a positive 
constant $\rho(K, g_{\sigma_0})\le K$ depending only on $K$ and $g_{\sigma_0}$ (in terms of $|g_{\sigma_0}|_{C^{1+\mu}}$) with the following 
properties. Assume $0<T\le 1$,
\ba\label{nu-condition}
T\|\nu_{\sigma_1, \sigma_0}\|_{C^0}\le \f{1}{4}\min\{\epsilon_0, \hat \delta_0\},
\ea
\ba \label{A-condition1}
\|\tau_{\sigma_0}\|_{C^{1+\mu}}+\|\sigma_1-\sigma_0\|_{C^{2+\mu}}\le K,
\ea
and
\ba \label{A-condition2}
[\nu_{\sigma_0, \sigma_1}]_{\mu}+T^{(1-\mu)/2}\|\nu_{\sigma_1, \sigma_0}\|_{C^{2+\mu}}
\le \rho(K, g_{\sigma_0}).
\ea
Then there is a closed
$\mc{C}^{2+\mu, (2+\mu)/2}$-solution $\sigma=\sigma(t)$ of the $\sigma_0$-gauged 
Laplacian flow on $[0, T]$ with $\sigma(0)=\sigma_1$, 
such that 
\ba \label{bound0}
\|\sigma-\sigma_0\|_{C^0} \le \min\{\epsilon_0, 
\hat \delta_0\}
\ea
and
\ba \label{bound1}
\|\sigma-\sigma_0\|_{\mc{C}^{2+\mu, (2+\mu)/2}} \le 
5K+\f{1}{2}\min\{\epsilon_0,\hat \delta_0\}.
\ea

Let $l>2+\mu$ be a non-integer. If $\sigma_0 \in C^{l}$, then 
$\sigma\in \mc{C}^{l, l/2}$ for $t>0$. If in addition $\sigma_1 \in C^{l}$, then we have
$\sigma\in \mc{C}^{l, l/2}$ on $M \times [0, T]$. 

Finally, the solution depends smoothly on $\sigma_1$ and $\sigma_0$.
\end{theo}

\begin{theo} \label{gauge2} ($\mc{C}^{4+\mu, (4+\mu)/2}$ estimates) Assume $\sigma_0 \in C^{4+\mu}$ and 
$\sigma_1 \in C^{4+\mu}$, and everything 
as in Theorem \ref{gauge1}, except that  
$\hat \delta_0=\hat \delta_0(\sigma_0)=
\hat \delta_0(4+\mu, g_{\sigma_0})$ (from Lemma \ref{linearlemma}
with $l=4+\mu$), 
which depends only on $g_{\sigma_0}$ (in terms of $|g_{\sigma_0}|_{C^{3+\mu}}$).
In addition, assume  
\ba
\|\tau_{\sigma_0}\|_{C^{3+\mu}}+\|\sigma_1-\sigma_0\|_{C^{3+\mu}}\le K.
\ea
 Then all the conclusions of Theorem \ref{gauge1} hold true. Moreover,
there holds 
\ba \label{bound2}
\|\sigma-\sigma_0\|_{\mc{C}^{4+\mu, (4+\mu)/2}} 
\le C(K, g_{\sigma_0})
(\|\sigma_1-\sigma_0\|_{C^{4+\mu}}+\|d\tau_{\sigma_0}\|_{C^{2+\mu}}).
\ea
for a positive constant $C(K, g_{\sigma_0})$ depending only on 
$K$ and $g_{\sigma_0}$ (in terms of $|g|_{C^{3+\mu}}$).
\end{theo}

\begin{theo} \label{gauge3} (global uniqueness)
Let $\sigma_0$ be a $C^2$ $G_2$-structure on $M$.
Let $\sigma=\sigma(t)$ and $\tilde \sigma=\tilde \sigma(t)$ be two $\mc{C}^{2,1}$ solutions of the 
$\sigma_0$-gauged Laplacian flow on a common interval 
$[0, T]$, such that $\sigma(0)=\tilde \sigma(0)$.
Then $\sigma \equiv \tilde \sigma$ on $[0, T]$. 
(See Appendix for the definition of $\mc{C}^{2,1}$.) 
\end{theo}

\noindent {\bf Remark} For the purpose of 
obtaining  existence and uniqueness of a short time solution of the Laplacian flow for a given initial $G_2$-structure $\sigma_0$, it suffices to consider the case $\sigma_0=\sigma_1$, cf.~the next section. In this case, the conditions in Theorem \ref{gauge1} and Theorem
\ref{gauge2} are obviously simplified.  Besides the independent interest of the gauge fixed Laplacian flow, we need to consider the general case of $\sigma_0$ for the purpose of obtaining regularity of solutions of the Laplacian flow and 
long time existence and convergence of the Laplacian flow starting near a torsion-free $G_2$-structure, cf.~Sections 7, 8 and 9.\\

We need some preparations for the proofs of these 
theorems. Let $C^l_o(\Lambda^3_+T^*M)$ denote the set of sections in 
$C^l_o(\Lambda^3 T^*M)$ with values in $\Lambda^3_+T^*M$, and $C^{l,l/2}_o(\pi^*(\Lambda^3_+T^*M))$ denote the set of sections in $C^{l,l/2}_o(\Lambda^3T^*M)$ with values in $\pi^*(\Lambda^3_+T^*M)$. By Lemma \ref{simple}, $C^l_o(\Lambda^3_+T^*M)$ is a domain in $C^l_o(\Lambda^3 T^*M)$, and $C^{l,l/2}_o(\pi^*(\Lambda^3_+T^*M))$ is a domain in $C^{l,l/2}_o(\Lambda^3T^*M)$. Consider $P(\sigma)=\f{\p \sigma}{\p t}-\Delta_{\sigma} \sigma-d (X_{\sigma_0}(\sigma-\sigma_0)\lr \sigma)$, the $P$-operator of 
$-\Delta_{\sigma}\sigma-d(X_{\sigma_0}(\sigma-\sigma_0)\lr \sigma)$:
\ba
P: C^{l,l/2}_o(\pi^*(\Lambda^3_+T^*M)) \rightarrow 
C^{l-2,(l-2)/2}_o(\pi^*(\Lambda^3T^*M)) 
\ea
and the corresponding $P$-map
\ba
{\mathcal P}: C^{l,l/2}_o(\pi^*(\Lambda^3_+T^*M)) \rightarrow 
C^{l-2,(l-2)/2}_o(\pi^*(\Lambda^3T^*M)) \times C^l_o(\Lambda^3 T*M).
\ea
It is obviously a smooth map.

\begin{lem} \label{linearlemma}  Let $\sigma_0\in 
C^{l}_o(\Lambda^3_+T^*M)$ for a noninteger $l>2$. Let $\delta_0=\delta_0(l,g_{\sigma_0})$ be the positive constant from Theorem \ref{full} with $g_*=g_0=g_{\sigma_0}$. (In this case, the dependence 
of $\delta_0$ on $\|g_0\|_{C^0}, 
\|g_0^{-1}\|_{C^0}, l$ and $g_*$ is reduced to 
the dependence on $l$ and $g_{\sigma_0}$, which is in terms of  $|g_{\sigma_0}|_{C^{l-1}}$.) Set
$\hat \delta_0=\hat \delta_0(l, g_{\sigma_0})= C_0^{-1}\delta_0$, where $C_0$ is from Lemma 
\ref{Alemma}.
  Let $\sigma \in C^{l, l/2}(\pi^*(\Lambda^3T^*M))$ such that $\|\sigma-\sigma_0\|_{C^0}
\le \min\{\epsilon_0, \hat \delta_0\}$  with $\epsilon_0$ being from Lemma \ref{simple}. (By Lemma \ref{simple}, we then have $\sigma\in C^{l, l/2}(\pi^*(\Lambda^3_+T^*M))$.) Then the linearization of $\mathcal P$ at $\sigma$:
\ba
{\mathcal D}_{\sigma}{\mathcal P}:
C^{l,l/2}_o(\pi^*(\Lambda^3T^*M)) \rightarrow 
C^{l-2,(l-2)/2}_o(\pi^*(\Lambda^3T^*M)) \times C^l(\Lambda^3 T*M)
\ea
is an isomorphism. Moreover, there hold
\ba \label{linearestimate}
\|{\mathcal D}_{\sigma}{\mathcal P}\| \le C 
\text{ and }  \|({\mathcal D}_{\sigma}{\mathcal P})^{-1}\| \le C 
\ea
for a positive constant $C=C(l, T, 
\|\tau_{\sigma_0}\|_{C^{l-1}}, 
\|\sigma-\sigma_0\|_{C^{l-1, (l-1)/2}}, g_{\sigma_0})$,
where the dependence on $g_{\sigma_0}$ is in terms 
of its Riemannian norm $|g_{\sigma_0}|_{C^{l-1}}$.
\end{lem}
\Pf There holds 
\ba \label{full-linear}
{\mathcal D}_{\sigma} {\mathcal P}(\gamma)=
({\mathcal D}_{\sigma} P, \gamma(\cdot, 0))
\ea
for all $\gamma$.
By 
(\ref{delta}) we have 
\ba \label{d}
{\mathcal D}_{\sigma} P(\gamma)=
\f{\p \gamma}{\p t}+\Delta_{\sigma_0} \gamma
+d({\mathcal D}_{\sigma} \Phi_{\sigma_0}(\gamma)).
\ea
 It follows from (\ref{phi-formula}) that
\ba \label{dd}
&&({\mathcal D}_{\sigma} \Phi_{\sigma_0})(\gamma)=
\Phi_0(\gamma)+\Phi_1(\nabla_{\sigma_0} \gamma),
\ea
where
\ba \label{ddd}
&&\Phi_0=D_{\sigma}A(\sigma_0, \cdot, \theta, \nabla \theta)+A(\sigma_0, \sigma, \cdot, \nabla \theta)
+D_{\sigma}B(\sigma_0, \cdot, \tau_{\sigma_0}, \theta)+B(\sigma_0, \sigma, \tau_{\sigma_0}, \cdot),
\nonumber \\
&&\Phi_1=A(\sigma_0, \sigma, \theta,\cdot),
\ea
and $\theta=\sigma-\sigma_0$ as before. 
By the assumption $\|\sigma-\sigma_0\|_{C^0}\le \epsilon_0$ and Lemma \ref{Alemma} there holds
$\|\Phi_1\|_{C^0} \le C_0\|\theta\|_{C^0}.$
Applying Theorem \ref{closed} and employing the nature 
of the functions $A$ and $B$ (cf.~(\ref{A-formula}) and 
(\ref{B-formula})) we arrive at the desired conclusions. \qed \\

\begin{lem} \label{hessianlemma} Let $\sigma_0 \in 
C^{l}_o(\Lambda^3_+T^*M)$ and $\sigma \in C^{l,l/2}_o(\pi^*(\Lambda^3_+T^*M)$ 
for a noninteger $l>2$. 
Assume $\|\sigma-\sigma_0\|_{\sigma_0} \le \epsilon_0$.  Then the second derivative 
operator of $\mc P$ at $\sigma$ 
\ba
&&{\mc D}^2_{\sigma} {\mc P}: C^{l,l/2}_o(\pi^*(\Lambda^3T^*M)) \times 
C^{l,l/2}_o(\pi^*(\Lambda^3T^*M)) \rightarrow 
\no \hspace{1.2cm} C^{l-2,(l-2)/2}_o(\pi^*(\Lambda^3T^*M)) \times C^l_o(\Lambda^3 T*M)
\ea
 satisfies the bound 
\ba \label{hessianbound}
\|{\mc D}^2_{\sigma} {\mc P}\| \le 
C
\ea
for a positive constant  $C=C(\|\tau_{\sigma_0}\|_{C^{l-1}},
\|\sigma-\sigma_0\|_{C^{l,l/2}})$.
\end{lem}
\Pf We have
\ba
{\mc D}^2_{\sigma} {\mc P}(\gamma, \gamma')=({\mc D}_{\sigma}^2 P(\gamma, \gamma'),0)
\ea
and
\ba
{\mc D}^2_{\sigma} P(\gamma, \gamma')=
d({\mc D}^2_{\sigma}\Phi_{\sigma_0}(\gamma, \gamma')).
\ea
By (\ref{dd}) we have
\ba
({\mc D}^2_{\sigma}\Phi_{\sigma_0})(\gamma, \gamma')
=(D_{\sigma} \Phi_0)(\gamma', \gamma)+(D_{\sigma}
\Phi_1)(\gamma', \gamma).
\ea
By the formulas (\ref{ddd}) for $\Phi_0$ and $\Phi_1$ and the nature of $A$ and $B$ (cf.~(\ref{A-formula}) and 
(\ref{B-formula})) we then 
deduce the bound (\ref{hessianbound}). \qed \\

\noindent {\it Proof of Theorem \ref{gauge1} } Consider the above $P$-map ${\mc P}$ with $l=2+\mu$. 
Let $\hat \delta_0=\hat \delta_0(2+\mu, g_{\sigma_0})$ be given by  Lemma \ref{linearlemma}
for $l=2+\mu$, as stated in the theorem. 
 Set 
\ba
\tilde \sigma_1 &=&\sigma_1+t\nu_{\sigma_1, \sigma_0} 
\ea
 and 
\ba
\tilde \theta=\tilde \sigma_1-\sigma_0.
\ea

For convenience, we set $\nu=\nu_{\sigma_1, \sigma_0}$.
Then we have 
\ba \label{check}
\|\tilde \theta\|_{C^0} \le \|\sigma_1-\sigma_0\|_{C^0}+T
\|\nu\|_{C^0} \le
\min\{\f{\epsilon_0}{2}, \f{\hat \delta_0}{2}\}
\ea
on account of (\ref{bound1}).  
Moreover, one readily deduces by using the definition of 
the ${\mc C}^{2+\mu, (2+\mu)/2}$ norm 
\ba \label{tilde-s}
\|\tilde \theta\|_{C^{2+\mu, (2+\mu)/2}} &\le&  \|\sigma_1-\sigma_0\|_{C^{2+\mu}}+T \|\nu\|_{C^{2+\mu}}+\|\nu\|_{C^{\mu}}+T^{\f{1-\mu}{2}}\|\nabla \nu\|_{C^0}+T^{\f{2-\mu}{2}}\|\nabla^2 \nu\|_{C^0} \nonumber \\
&\le& \|\sigma_1-\sigma_0\|_{C^{2+\mu}}+4\|\nu\|_{C^{2+\mu}}.
\ea

 By Lemma \ref{linearlemma} and the assumptions of the 
 theorem,
$\mc P$ is an isomorphism and satisfies the estimates
 \ba \label{iso-estimate}
 \|{\mathcal D}_{\tilde \sigma_1} {\mathcal P}\|
 \le C \text{ and } \|({\mathcal D}_{\tilde \sigma_1} {\mathcal P})^{-1}\|
 \le C
 \ea
for a positive constant $C$ depending only on 
$K$ and $g_{\sigma_0}$. (In this proof, the directly indicated dependence of various constants on  $g_{\sigma_0}$ are all in terms of 
$|g_{\sigma_0}|_{C^{1+\mu}}$.)

By Lemma \ref{hessianlemma} and (\ref{iso-estimate}) we can apply the inverse function theorem
for mappings of Banach spaces to deduce that
$\mathcal P$ is a smooth diffeomorphism from 
an open neighborhood $U_{\tilde \sigma_1}$ of 
$\tilde \sigma_1$ in $\mc{C}^{2+\mu}$ onto the open ball 
$B_{r_1}({\mathcal P}(\tilde \sigma_1))$ in 
$\mc{C}^{\mu, \mu/2}$, for a positive number 
$r_1$ having the same dependences as the above $C$.
Moreover, $U_{\tilde \sigma_1} \subset 
B_{r_2}(\tilde \sigma_1)$ for a positive number 
$r_2$ with the same dependences as $C$. In addition, we choose $r_1$ and $r_2$ such that 
$r_2\le \f{1}{2}\min\{\epsilon_0, \hat \delta_0\}$.  On the other hand, we have ${\mc P}(\tilde \sigma_1)=(P(\tilde \sigma_1), \sigma_1)$ and
\ba
P(\tilde \sigma_1)&=&\f{\p}{\p t} \tilde \theta
+[\Delta_{\sigma_0} \tilde \theta+d(\Phi_{\sigma_0}(\tilde \theta))] \nonumber \\ &=&[-\Delta_{\sigma_0}\theta-d(\Phi_{\sigma_0}(\theta)]+[\Delta_{\sigma_0} \tilde \theta+
d(\Phi_{\sigma_0}(\tilde \theta))] \nonumber \\
&=& t \Delta_{\sigma_0} \nu
+d(\Phi_{\sigma_0}(\theta+t\nu)-
\Phi_{\sigma_0}(\theta)).
\ea
Employing the formula (\ref{phi-formula}) and 
calculating as in the proof of Lemma 
\ref{nulemma} we deduce
\ba
\|{\mc P}(\tilde \sigma_1)-(0, \sigma_1)\|_{C^{\mu, \mu/2}} \le C_1 \|t\nu\|_{\mc{C}^{2+\mu, (2+\mu)/2}},
\ea
where $C_1$ depends only on  $\|\sigma_0\|_{C^{1+\mu}},
\|\theta\|_{C^{2+\mu}}, \|t\nu\|_{\mc{C}^{1+\mu, (1+\mu)/2}}$ and $\|\tau_{\sigma_0}\|_{C^{1+\mu}}$.
Since $\nabla \sigma_0$ can be expressed in terms of 
$\tau_{\sigma_0}$ (by Lemma \ref{Klemma}), 
we have $\|\sigma_0\|_{C^{1+\mu}}\le C_2(1+\|\tau_{\sigma_0}\|_{C^{\mu}})$ for a universal positive constant $C_2$. 
  There hold  
\ba
\|t\nu\|_{\mc{C}^{1+\mu, (1+\mu)/2}}
\le \max\{T, T^{(1-\mu)/2}, T^{(2-\mu)/2}\}
\|\nu\|_{C^{1+\mu}}\le T^{(1-\mu)/2}\|\nu\|_{C^{1+\mu}}
\ea
and
\ba
\|t\nu\|_{\mc{C}^{2+\mu, (2+\mu)/2}}
\le [\nu]_{\mu}
+\max\{T, T^{(1-\mu)/2, (2-\mu)/2}\}
\|\nu\|_{C^{2+\mu}} \le [\nu]_{\mu}
+T^{(1-\mu)/2}
\|\nu\|_{C^{2+\mu}} .
\ea
We infer 
\ba
\|{\mc P}(\tilde \sigma_1)-(0, \sigma_1)\|_{C^{\mu, \mu/2}}
\le  C_3([\nu]_{\mu}
+T^{(1-\mu)/2}
\|\nu\|_{C^{2+\mu}})
\ea
for a positive constant $C_3$ depending only on 
$K$ and $g_{\sigma_0}$.
We set $\rho(K, g_{\sigma_0})=\min\{2^{-1}C_3^{-1} r_1, K\}$. 
Then the condition  (\ref{A-condition2}) implies that  
$(0, \sigma_1) \in B_{r_1}({\tilde P}(\tilde \sigma_1))$.  Consequently, we obtain a solution $\sigma=
{\mathcal P}|_{U_{\tilde \sigma_0}}^{-1}(0, \sigma_1)$ of the $\sigma_0$-gauged Laplacian flow on $[0, T]$ with 
$\sigma(0)=\sigma_1$. There holds 
\ba
\|\sigma-\sigma_0\|_{\mc{C}^{2+\mu, (2+\mu)/2}}
\le \|\sigma-\tilde \sigma_1\|_{\mc{C}^{2+\mu, (2+\mu)/2}}
+\|\tilde \sigma_1-\sigma_0\|_{\mc{C}^{2+\mu, (2+\mu)/2}}
\le 5K+\f{1}{2}\min\{\epsilon_0,\hat \delta_0\}.
\ea
On the other hand, we have
\ba \label{thetaC}
\|\sigma-\sigma_0\|_{C^0}&\le& \|\sigma-\tilde \sigma_1\|_{C^0}
+\|\sigma_1-\sigma_0\|_{C^0}+T\|\nu\|_{C^0}
\nonumber \\ &\le& r_2+\min\{\f{\epsilon_0}{2}, 
\f{\hat \delta_0}{2}\}\le \min\{\epsilon_0, 
\hat \delta_0\}. 
\ea

The $\mc{C}^{l, l/2}$ regularity for $t>0$ and the 
$\mc{C}^{4+\mu, (4+\mu)}$ and $\mc{C}^{l,l/2}$ regularities up to $t=0$ (under the corresponding given assumptions) follow from the standard 
regularity arguments in coordinate charts, cf. [LSU] and [Y2].

The smooth dependence of $\sigma$ on $\sigma_0$ and 
$\sigma_1$ follows from the formula 
$\sigma=
{\mathcal P}|_{U_{\tilde \sigma_0}}^{-1}(0, \sigma_1)$
and the inverse function theorem. \qed \\

\noindent {\it Proof of  Theorem \ref{gauge2}} Since 
$\hat \delta_0(4, g_{\sigma_0}) \le \hat \delta_0(2, g_{\sigma_0})$ (cf.~Theorem \ref{full}), 
the conclusions of Part I hold true.  
To derive the estimate (\ref{bound2}) we view the 
gauge fixed Laplacian flow (\ref{newflow}) with the 
solution $\sigma$ as a linear equation in the form (\ref{closedheat}), i.e.
\ba \label{linearview}
\f{\p \theta}{\p t}+\Delta \theta
+d(\Phi_0(\theta)+\Phi_1(\nabla \theta))=d\tau_{\sigma_0},
\ea
 with $\theta=\sigma-\sigma_0$,  $\Delta=\Delta_{\sigma_0}$, $\Phi_0(\theta)=B(\sigma_0, \sigma, \tau_{\sigma_0}, \gamma)$ and $\Phi_1(\nabla \gamma)=A(\sigma_0, \sigma, \theta, 
\nabla \gamma).$
By (\ref{thetaC}) we deduce 
\ba
\|\Phi_1\|_{C^0}\le C_0 \hat \delta_0 \le \delta_0(4+\mu,
g_{\sigma_0})\le \delta_0(3+\mu, g_{\sigma_)}),
\ea
where $C_0$ is from Lemma \ref{Alemma}.
Hence we can apply Theorem \ref{closed} (or Theorem \ref{full}) to deduce 
\ba
\|\theta\|_{\mc{C}^{3+\mu, (3+\mu)/2}}
\le C_1(\|\sigma_1-\sigma_0\|_{C^{3+\mu}}
+\|d\tau_{\sigma_0}\|_{C^{1+\mu}}),
\ea
where $C_1$ depends only on $\|\tau_{\sigma_0}\|_{C^{2+\mu}}, 
\|\theta\|_{C^{2+\mu, (2+\mu)/2}}$ and  $|g_{\sigma_0}|_{C^{2+\mu}}$.  Combining 
this with the estimate (\ref{bound1}) we deduce 
an estimate for $\|\theta\|_{\mc{C}^{3+\mu, (3+\mu)/2}}$.
Then we repeat the argument for $4+\mu$ instead of 
$3+\mu$ and arrive at the estimate (\ref{bound2}). 
\\

\noindent {\it Ptoof of Theorem \ref{gauge3}} Let $\sigma$ and $\tilde \sigma$ be two
$C^{2,1}$ solutions of the $\sigma_0$-gauged Laplacian 
flow  on 
an interval $[0, T]$, such that 
$\tilde \sigma(0)=\sigma_1$.  Set $\theta=\sigma-\sigma_1, \tilde \theta=
\tilde \sigma-\sigma_1$ and $\gamma=\tilde \sigma-\sigma$. By Lemma \ref{transform} and the 
Bochner-Weitzenb\"{o}ck formula 
(\ref{bochner}) we have 
\ba
\f{\p \gamma}{\p t}=-\nabla^* \nabla \gamma
-{\mathcal R}(\gamma)
-d(\Phi_{\sigma_0}(\tilde \sigma-\sigma_0)-
\Phi_{\sigma_0}(\sigma-\sigma_0)).
\ea
Multiplying this equation by $\gamma$ and then integrating lead to 
\ba \label{top}
\f{d}{dt} \int_M |\gamma|^2 
+\int_M |\nabla \gamma|^2=-\int_M {\mathcal R}(\gamma) \cdot \gamma- \int_M (\Phi_{\sigma_0}(\tilde \sigma-\sigma_0)-\Phi_{\sigma_0}(\sigma-\sigma_0)) \cdot d^* \gamma.
\ea
By (\ref{phi-formula}) we deduce 
\ba
&&\Phi_{\sigma_0}(\tilde \sigma-\sigma_0)-
\Phi_{\sigma_0}(\sigma-\sigma_0)=A(\sigma_0,
\sigma, \theta, \nabla \gamma)+
(A(\sigma_0, \tilde \sigma, \tilde \theta, \nabla
\tilde \theta)-A(\sigma_0, \sigma, \theta, \nabla \tilde \theta)) \nonumber \\
&&\hspace{3.5cm}+(B(\sigma_0, \tilde \sigma, \tau_{\sigma_0},
\tilde \theta)-B(\sigma_0,\sigma, \tau_{\sigma_0},
\theta)).
\ea
It follows that
\ba \label{middle}
|\Phi_{\sigma_0}(\tilde \sigma-\sigma_0)-
\Phi_{\sigma_0}(\sigma-\sigma_0)| \le C_1 |\gamma|
+C_2 |\theta| \cdot |\nabla \gamma|
\ea
on $M \times [0, T]$, with positive numbers $C_1$ and $C_2$.  Combining (\ref{top}) with 
(\ref{middle}) we then deduce 
\ba \label{bottom}
\f{d}{dt}\int_M |\gamma|^2
+\int_M |\nabla \gamma|^2 \le \int_M(C_3|\gamma|^2+
C_4 |\theta|^2 |\nabla \gamma|^2)
\ea
for positive numbers $C_3$ and $C_4$. 
There holds  for $0<T_1\le T$
\ba \label{theta-estimate}
\max_{[0, T_1]}|\theta|\le C_5 T_1, 
\ea 
where $C_5=\max_{M \times [0, T]}|\f{\p \sigma}{\p t}|$. 
Now we assume $T_1 \le C_4^{-1/2}  C_5^{-1}$. Then (\ref{bottom}) yields 
\ba
\f{d}{dt}\int_M |\gamma|^2 \le C_3\int_M |\gamma|^2
\ea
on $[0, T_1]$,
which implies $\gamma\equiv 0$ on $[0, T_1]$ because $\gamma=0$ at $t=0$.  Then we repeat the above argument 
with the time $0$ replaced by $T_1$, and with 
the new definitions $\theta=\sigma-\sigma(T_1)$ and 
$\tilde \theta=\tilde \sigma-\sigma(T_1)$. After finitely many such steps we then conclude 
$\gamma=0$ on $[0, T]$. 
\qed \\

\sect{Short time solutions of the Laplacian flow}
\vspace{2mm}

In this section we apply the results from the last section to prove Theorem \ref{existence}. In Subsection 6.1 we prove the existence part, and in Subsection 6.2
we prove the uniqueness part. 

\subsection{Existence}
\hspace{1cm}\vspace{1mm}

Let $l>4$ be a non-integer, and $\sigma_0\in C^{l}$ 
and $\sigma_1\in C^l$ be closed $G_2$-structures on $M$.
Let $\sigma=\sigma(t)$ be a $\mc{C}^{l, l/2}$ solution of the $\sigma_0$-gauged Laplacian flow on an interval $[0, T]$ with the initial value $\sigma_1$.  We consider the ODE
\ba \label{ODE}
\f{d}{dt}\phi(\cdot,t)=-X_{\sigma_0}(\sigma(t)-\sigma_0)(\phi(\cdot, t))
\ea
with the initial condition
\ba \label{initial}
\phi(\cdot, 0)=Id,
\ea
where $X_{\sigma_0}$ is given by (\ref{X}).
By a difference quotient argument one easily shows that 
the solution $\phi$ inherits the spacial regularity of $X_{\sigma_0}(\sigma-\sigma_0)$, and improves a $C^k$ regularity of it in the time direction 
to $C^{k+1}$ regularity.
Since $X_{\sigma_0}(\sigma-\sigma_0) \in {\mc C}^{l-1, (l-1)/2}$, we deduce that $\phi \in {\mc C}^{l-1, l/2}$.

\begin{lem} \label{goback}  $ \hat \sigma(t)=\phi(\cdot, t)^* \sigma(t)$
is a ${\mc C}^{l-2, (l-1)/2}$
solution of the Laplacian flow on $[0, T]$ with 
$\sigma(0)=\sigma_1$.  
\end{lem}
\Pf By the above regularity of $\phi$ we obviously have 
$\hat \sigma\in {\mc C}^{l-2, (l-1)/2}$. (The derivative 
of $\phi$ involved in $\phi^* \sigma$ causes the drop of regularity.)  Now we have with $X_{\sigma_0}=X_{\sigma_0}(\sigma-\sigma_0)$
\ba
\f{\p \hat \sigma}{\p t} &=&
\phi(\cdot, t)^* ({\mc L}_{-X_{\sigma_0}}\sigma)
+\phi(\cdot, t)^*(\Delta_{\sigma} \sigma+{\mc L}_{X_{\sigma_0}} \sigma)\nonumber \\
&=&\Delta_{\phi(\cdot, t)^*\sigma} \phi(\cdot, t)^* \sigma \nonumber \\
&=&\Delta_{\hat \sigma} \hat \sigma.
\ea
It is obvious that $\hat \sigma(0)=\sigma_1$.
 \qed \\ 

\subsection{Uniqueness}
\hspace{1cm}\vspace{1mm}

Let $\sigma$ be a 
function with closed $G_2$-structures as values on 
a time interval $[0, T]$ such that $ \sigma(0)=\sigma_1$. 
Analogous to the situation in [BX], we consider the 
following nonlinear evolution equation for diffeomorphims
\ba \label{diff-equation}
\f{\p \phi}{\p t}=Z_{\sigma_0, \sigma}(\phi)\circ \phi
\ea
with the initial condition $\phi(\cdot, 0)=Id$, where
\ba
Z_{\sigma_0, \sigma}(\phi)=-X_{\sigma_0}((\phi(\cdot, t)^{-1})^*\sigma-\sigma_0).
\ea

\begin{theo} \label{diff-theorem} 
 Assume $\sigma_0\in C^{l}$ and $\sigma \in {\mc C}^{l, l/2}$ for a non-integer $l>2$. Then there is a unique
${\mc C}^{l+1, (l+1)/2}$ solution of the evolution equation 
(\ref{diff-equation}) on a time interval 
$[0, T_1]$ ($T_1>0$) with $\phi(\cdot, 0)=Id.$
\end{theo}
\Pf We compute the linearization of $Z_{\sigma_0, \sigma}$ at a given diffeomorphism $\phi:M\rightarrow M$ in the direction of a vector field $Y$ along $\phi$.
(Thus $Y(p) \in T_{\phi(p)}M$ for each $p\in M$.) 
Let $\phi_s$ be a family of diffeomorphisms of $M$ such that $\phi_0=\phi$ and 
\ba
\frac{d}{ds}\phi_s|_{s=0}=Y. 
\ea
(For example, $\phi_s=exp_{\phi}(sY)$ for small $s$, where $exp$ is the exponential map of a Riemannian metric on $M$.)
There hold
\ba
\frac{d}{ds}\phi_s^{-1}|_{s=0}=-(\phi^{-1})_*Y
\ea
and 
\ba
\frac{d}{ds}(\phi_s^{-1})^*\sigma|_{s=0}=-\mathcal{L}_Y((\phi^{-1})^*\sigma).
\ea
Set $\bar \sigma=(\phi^{-1})^*\sigma$, which is a closed $G_2$-structure. Then we have
\ba \label{flowlinear}
{\mathcal D}_{\phi} Z_{\sigma_0, \sigma}(Y)&\equiv&  -\frac{d}{ds}X_{\sigma_0}((\phi_s^{-1})^*\sigma-\sigma_0)|_{s=0}
\nonumber \\
&=& X_{\sigma_0}(\mathcal{L}_Y(\bar \sigma)) \nonumber \\ &=&
X_{\sigma_0}(d(Y\lr \bar \sigma))
\nonumber \\ &=&
X_{\sigma_0}(d(Y\lr \sigma_0))+X_{\sigma_0}(d(Y\lr (\bar \sigma-\sigma_0)).
\ea
As is easy to verify, there holds $Y\lrcorner \sigma_0= *(Y^{\flat}\wedge *\sigma_0)$ (cf.~[B2]), where $*=*_{\sigma_0}$ and $Y^\flat$ is the 1-form dual to $X$ w.r.t.~$g_{\sigma_0}$.
By Lemma \ref{2-7}  we have  
\ba \label{2-7-0}
d*(Y^\flat\wedge*\sigma_0)=-\frac{3}{7}d^7_1Y^\flat\sigma_0-\frac{1}{2}*(d^7_7Y^\flat\wedge\sigma_0)+d^7_{27}Y^\flat+\zeta(Y^\flat)
\ea
with $\zeta=\zeta_{\sigma_0}$. By (\ref{X}), Lemma \ref{hodge-1form}, (\ref{2-7-0}) and the Bochner-Weitzenb\"{o}ck 
formula we then infer
\ba \label{XX}
X_{\sigma_0}(d(Y\lrcorner \sigma_0))&=&(-d^1_7d^7_1Y^\flat-d^7_7d^7_7Y^\flat)_\#+
X_{\sigma_0}(\zeta(Y^\flat)) \nonumber \\
&=&-(\Delta_{\sigma_0} Y^{\flat})_\# 
+\f{1}{3}*d*\xi_{\sigma_0}\left(*(*\sigma_0 \wedge *(Y^\flat \wedge *\tau_{\sigma_0}))\right)+X_{\sigma_0}(\zeta(Y^\flat)) \nonumber \\
&=&-\nabla^* \nabla Y-{\mathcal{R}}(Y^\flat)_\#
+\f{1}{3}*d*\xi_{\sigma_0}\left(*(*\sigma_0 \wedge *(Y^\flat \wedge *\tau_{\sigma_0}))\right)+ X_{\sigma_0}(\zeta(Y^\flat)). 
\ea                         
 On the other hand, we have
\ba
X_{\sigma_0}(d(Y\lrcorner (\bar \sigma-\sigma_0))
&=&A_2(\sigma_0, \bar \sigma-\sigma_0, \nabla^2 Y)+A_1(\sigma_0, \nabla (\bar \sigma-\sigma_0), \nabla Y)
\no+A_0(\sigma_0, \nabla^2 (\bar \sigma-\sigma_0), Y),
\ea
where $A_0, A_1$ and $A_2$ are universal pointwise functions, 
smooth in their first arguments, and linear in their 
second and third arguments. We arrive at 
\ba
{\mathcal D}_{\phi} Z_{\sigma_0, \sigma}(Y)=-\nabla^*
\nabla Y+A_2(\sigma_0, \bar \sigma-\sigma_0, \nabla^2 Y)
+W_{\sigma_0, \sigma, \phi}(Y, \nabla Y),
\ea
where the first order linear differential operator $W_{\sigma_0, \sigma, \phi}(Y, \nabla Y)$ is the sum of $A_0, A_1$ and the lower order terms on the far right hand side of (\ref{XX}).  
Obviously, ${\mathcal D}_{\phi} Z_{\sigma_0, \sigma}$ is strongly elliptic when $\|\bar \sigma-\sigma_0\|_{C^0}$ is small enough, which is the case for small time because of the fact $\bar \sigma(0)=\sigma_0$. 
Now we can apply the general result on evolutions of mappings in [Y2] to deduce the desired existence and uniqueness of short time solutions. (The basic mechanism for the said general result is similar to Theorem 
\ref{gauge1} and its proof.) \qed \\

Now we are ready to prove Theorem \ref{existence}. (Note that only the existence part of Theorem \ref{diff-theorem} is needed.) \\

\noindent {\it Proof of Theorem \ref{existence}} 
Let $\sigma_1$ be a ${C}^{4+\mu}$ closed $G_2$-structures 
on $M$. We choose a $C^{\infty}$ closed $G_2$-structure $\sigma_0$ sufficiently close to $\sigma_1$ in 
$C^{3+\mu}$, such that the conditions of Theorem \ref{gauge1} and Theorem \ref{gauge2} are satisfied for $K=1$ and a suitable $T>0$. Let $\sigma$ denote the unique ${\mc C}^{4+\mu, (4+\mu)/2}$ solution of the $\sigma_0$-gauged 
Laplacian flow on $[0, T]$ with the initial value $\sigma_1$ as given by Theorem \ref{gauge1}.  Applying 
Lemma \ref{goback} we then obtain 
a closed ${\mc C}^{2+\mu, (3+\mu)/2}$ solution 
$\hat \sigma(t)=\phi(\cdot, t)^* \sigma(t)$ of the Laplacian flow with the initial value 
$\sigma_1$ on the time interval $[0, T]$ given by Theorem \ref{gauge1}. The claimed ${\mc C}^{l-2, 
(l-1)/2}$ regularity follows from Theorem \ref{gauge1} 
and Lemma \ref{goback}.  

For a given $0<\epsilon<T$. Let $\psi$ be the solution of (\ref{ODE}) on $[\epsilon, T]$ with 
the initial value $Id$ at $t=\epsilon$. 
By Theorem \ref{gauge1}, $\sigma$ is smooth for $t>0$. Hence
$\psi$ is smooth.  Consequently, the pullback $\psi(\cdot, t)^*\sigma(t)$ is a smooth solution of the Laplacian flow on the time interval $[\epsilon, T]$. 
Obviously, it equals $\phi(\cdot, t)^* \sigma(t)$ for 
a family of diffeomorphisms $\phi(\cdot, t)$ of class 
${\mc C}^{3+\mu, (4+\mu)/2}$.  

Renaming $\hat \sigma$ we then obtain a desired solution of the Laplacian flow $\sigma=\sigma(t)$.

Next we show the uniqueness. Let $\gamma_1=\gamma_1(t)$ and $\gamma_2=\gamma_2(t)$ be 
two ${\mc C}^{2+\mu,(2+\mu)/2}$ solutions of the Laplacian flow on 
a common interval $[0, T]$ for some $T>0$, such that 
$\gamma_1(0)=\gamma_2(0)$.  We set $\sigma_0=\gamma_1(0)=\gamma_2(0)$. For $i=1,2$, let $\phi_{\gamma_i}$ be 
the ${\mc C}^{3+\mu, (3+\mu)/2}$ solution of the 
equation (\ref{diff-equation}) on an interval 
$[0, T_i]\subset [0, T]$, with $\gamma_i$ playing the role of $\sigma$ and with the initial value 
$Id$, as provided by Theorem \ref{diff-theorem}.  We set 
for each $i$
\ba \label{hatgamma}
\hat \gamma_i(t)=(\phi_{\gamma_i}(\cdot, t)^{-1})^* \gamma_i(t),
\ea
 which is of class $\mc{C}^{2+\mu, (2+\mu)/2} \subset {\mc C}^{2, 1}$.  Then we have for $i=1,2$
\ba
\f{\p \hat \gamma_i}{\p t}&=&(\phi_i(\cdot, t)^{-1})^*{\mc L}_{(\phi(\cdot, t)^{-1})_*X_{\sigma_0}(\hat \gamma_i-\sigma_0)}
\gamma_i+ (\phi_{\gamma_i}(\cdot, t)^{-1})^* \Delta_{\gamma_i}
\gamma_i \nonumber \\
&=& {\mc L}_{X_{\sigma_0}(\hat\gamma_i-\sigma_0)} \hat \gamma_i+\Delta_{\hat\gamma_i} \hat\gamma_i.
\ea
Thus, for each $i$, $\hat \gamma_i$ is a ${\mc C}^{2,1}$ solution of the $\sigma_0$-gauged Laplacian flow (\ref{gauged}) on 
$[0, T_i]$. Obviously, we also have
$\hat \gamma_i(0)=\sigma_0$ for each $i$. Hence  Theorem \ref{gauge3} implies that $\hat \gamma_1$ and
$\hat \gamma_2$ agree on $[0, T_0]$ for  $T_0=\min\{T_1, T_2\}$.  Consequently, for each $i$, $\phi_{\gamma_i}$ satisfies on $[0, T_0]$ the same ODE
\ba
\f{\p \phi}{\p t}=-X_{\sigma_0}(\hat \gamma-\sigma_0)
\circ \phi,
\ea
where $\hat \gamma$ stands for $\hat \gamma_1=\hat \gamma_2$.
Since $\phi_{\gamma_1}$ and $\phi_{\gamma_2}$ have the same initial
value, 
we deduce that $\phi_{\gamma_1} \equiv \phi_{\gamma_2}$
on $[0, T_0]$. By (\ref{hatgamma}) we then infer that  
$\gamma_1 \equiv \gamma_2$ on $[0, T_0]$. 

Next we repeat the above argument with $T_0$ being the new time origin. This way we can extend the interval on which $\gamma_1$ and $\gamma_2$ agree. By  a simple continuity argument we then conlcude that 
$\gamma_1 \equiv \gamma_2$ on $[0, T]$. \qed\\

\sect{Long time existence and convergence of the gauge
fixed Laplacian flow}
\vspace{2mm}

\subsection{A Sobolev-type  inequality} 
\hspace{1cm}\vspace{1mm}

Consider  a compact manifold $N$ of dimension $n\ge 3$ equipped with a Riemannian metric $g$. The Sobolev 
constant $C_S(N, g)$ (for the exponent $2$) is defined to be the smallest positive number for which the following Sobolev inequality holds 
\ba \label{sobolev}
\int_N |f|^{\f{2n}{n-2}} dvol \le C_{S}(N, g)\int_N |\nabla f|^2 dvol+V_g(N)^{-\f{2}{n}} \int_N f^2 dvol
\ea
for all $f\in C^1(N)$, where $V_g(N)$ denotes the volume of $(N, g)$.  We have the following $L^1$ version of Moser type 
maximum principle. 

\begin{theo} \label{moser}
 Let $T>0$ and $f$ be a nonnegative Lipschitz continuous function on $M \times [0,T]$ satisfying
\ba \label{linearequation1}
\frac{\partial f}{\partial t} \leq -\Delta f + bf
\ea
on $N \times [0, T]$ in the sense of distributions, where $\Delta$ denotes the Hodge Laplacian on functions and $b$ is a nonnegative constant.
 Then
we have for each $p \in N$ and $0<t \le T$
\ba \label{moser-estimate}
\max_{M \times [t, T]} |f|\le t^{-\f{n+2}{4}}C_n
\left(\max\{b, \f{n}{4}(1+\f{n}{2})^2\}\right)^2 
T^{\f{n+2}{2}}\left( \max\{C_{S}(N,g), TV_g(N)^{-\f{2}{n}}\}\right)^{\frac{n}{2}}
 \int_{M\times [0, T]} |f|,
\ea
where $C_n$ is a positive constant depending only on 
$n$.
\end{theo}

\Pf This is the global formulation of a corresponding 
local version in [Y5] (see also [Sa]). First we have
\ba \label{moser2}
\max_{M \times [t, T]}|f|&\leq& (1+\frac{2}{n})^{\frac{c_n}{2}} \left( \max\{C_{S}(N,g), TV_g(N)^{-\f{2}{n}}\}\right)^{\frac{n}{4}}\left(
2b+\frac{n}{2}(1+\frac{n}{2})^2 \cdot \frac{1}{t}\right) ^{\f{n+2}{4}} \nonumber \\
&& \cdot
\Bigl(\int^T_0\int\limits_{N} f^{2}(\cdot, s)dvol ds\Bigr)^{\frac{1}{2}},
\ea
for all $0<t\le T$,
where $c_n=\sum_0^{\infty} 2k(1+\frac{2}{n})^{-k}.$
This is the global formulation of a corresponding local version in [Y4]. We can adapt its proof in [Y4]. The cut-off function $\eta$ in that proof 
is not needed here, hence we can take $\eta \equiv 1$.
The local Sobolev inequality used there is 
replaced by (\ref{sobolev}). Then the arguments there 
lead to (\ref{moser}) straightforwardly

Applying (\ref{moser2}) to $0<t'<t\le T$ with $t'$ playing the role of the time origin we infer
\ba
\max_{M \times [t, T]} |f| &\le& C\bar T^{\f{n}{4}}
(1+\f{1}{t-t'})^{\f{n+2}{4}} \left(\int_{M\times [t', T]} 
|f|^2\right)^{1/2} \nonumber \\
&\le& C\bar T^{\f{n}{4}}
(1+\f{1}{t-t'})^{\f{n+2}{4}} (\max_{M \times [t', T]} |f|^{1/2}) \cdot \left(\int_{M\times [0, T]} |f|\right)^{1/2},
\ea
where $C=(1+\f{2}{n})^{c_n/2}$ and 
$\bar T=\max\{C_{S}(N,g), TV_g(N)^{-\f{2}{n}}\}$. We may assume that $\max_{M\times [t, T]} |f|$ is positive. (If it is zero, then the estimate (\ref{moser}) holds true trivially.) Then we deduce for $i\ge 0$ 
\ba
\f{\max_{M \times [t, T]} |f|^{2^{-i}}}{\max_{M \times [t', T]}
|f|^{2^{-i-1}}} \le \left(C\bar T^{\f{n}{4}}\right)^{2^{-i}}
(1+\f{1}{t-t'})^{\f{n+2}{4} \cdot 2^{-i}}\left(\int_{M\times [0, T]} |f|\right)^{2^{-i-1}}.
\ea
For a sequence of positive times $t=t_0>t_1>t_{k+1}, k\ge 0$ we  infer from it
\ba
\f{\max_{M \times [t, T]} |f|}{\max_{M \times [t_{k+1}, T]}
|f|^{2^{-k-1}}} \le \Pi_{0\le i \le k} \left(C \bar T^{\f{n}{4}}\right)^{2^{-i}}
(1+\f{1}{t_i-t_{i+1}})^{\f{n+2}{4} \cdot 2^{-i}}\left(\int_{M\times [0, T]} |f|\right)^{2^{-i-1}}.
\ea
Choosing $t_i=t(1-\sum_{1\le j \le i} 2^{-j})$ for 
$i\ge 1$ we have $1+\f{1}{t_i-t_{i+1}}=1+t^{-1}2^{i+1}
\le t^{-1}(T+2^{i+1})$. Hence we deduce
\ba
\f{\max_{M \times [t, T]} |f|}{\max_{M \times [t_{k+1}, T]}
|f|^{2^{-k-1}}} \le \Pi_{0\le i \le k} \left(C\bar T^{\f{n}{4}}\right)^{2^{-i}}
[t^{-1}(T+2^{i+1})]^{\f{n+2}{4} \cdot 2^{-i}}\left(\int_{M\times [0, T]} |f|\right)^{2^{-i-1}}.
\ea
Letting $k \rightarrow \infty$ we obtain
\ba
\max_{M \times [t, T]} |f|\le t^{-\f{n+2}{4}}C^2 \bar T^{\f{n}{2}}(\Pi_{0\le i \le \infty} 
(T+2^{i+1})^{\f{n+2}{4} \cdot 2^{-i}}) \int_{M\times [0, T]} |f|.
\ea
Replacing $T+2^{i+1}$ by $(T+1)2^{i+1}$ we then  arrive at 
(\ref{moser}) with 
\ba
C_n=4\left(1+\f{2}{n}\right)^{\sum_1^{\infty} 2i(1+\frac{2}{n})^{-i}} 2^{(n+2) \sum_0^{\infty}(i+1)2^{-i-2}}.
\ea
\qed\\

\subsection{Long time existence and convergence of the gauge fixed Laplacian flow: 
the statement and preliminaries}
\hspace{1cm}\vspace{1mm}

Consider a $G_2$ structure $\sigma_0$. Let $\lambda_0=\lambda_0(\sigma_0)$ denote the 
first eigenvalue of the Hodge Laplacian 
$\Delta_{\sigma_0}$ on exact 3-forms.  It is obviously positive, because a harmonic form which is also exact must be trivial. There holds
\ba
\int_M |d^*_{\sigma_0} \gamma|^2  
\ge \lambda_0 \int_M |\gamma|^2
\ea
for all exact $C^1$ 3-forms $\gamma$, or more generally, 
$W^{1,2}$ 3-forms.   Indeed, this is a consequence of the decomposition $\gamma=\sum_i a_i \gamma_i$, where 
the $\gamma_i$ are the exact forms among a complete set of $L^2$ orthonormal 
eigenforms of the Hodge Laplacian, see the proof of 
Theorem \ref{closed}. 

We have the following preliminary result.

\begin{lem} \label{starting} Let $0<\mu<1$, $K>0$, and $\hat \sigma_0$ be a $C^{4+\mu}$ torison-free $G_2$-structure on 
$M$.  Then there are 
positive constants $\lambda_0$ and $\rho$ depending only on 
$\hat \sigma_0$ and $K$, and $\hat \delta_0^*, T_0$ and $c$ depending only 
on $\hat \sigma_0, \mu$ and $K$, with the following properties. Let $\sigma_0$ and $\sigma_1$ be 
$C^{4+\mu}$ $G_2$-structure on $M$ such that 
\ba \label{con1}
\|\sigma_0-\hat \sigma_0\|_{C^0, \hat \sigma_0} \le \epsilon_0,
\ea
\ba \label{con2}
\|\sigma_0-\hat \sigma_0\|_{C^{4+\mu}, \hat \sigma_0}
\le K,
\ea 
\ba \label{con3}
\|\sigma_1-\sigma_0\|_{C^0, \sigma_0}\le \f{1}{4}\min\{\epsilon_0, \hat \delta_0(\sigma_0)\},
\ea
\ba \label{con4}
\|\sigma_1-\sigma_0\|_{C^{4+\mu}, \sigma_0} \le K,
\ea
and
\ba \label{con5}
\|\sigma_1-\sigma_0\|_{C^{2+\mu}, \sigma_0}\le \rho,
\ea
where $\hat \delta_0(\sigma_0)$ is from Theorem \ref{gauge2}.
Then there hold $\lambda_0(\sigma_0)\ge \lambda_0$ and 
$\hat \delta_0(\sigma_0)\ge \hat \delta_0^*$. 
  On the other hand, there is a unique $\mc{C}^{4+\mu, (4+\mu)/2}$ solution $\sigma=\sigma(t)$ of the 
$\sigma_0$-fixed Laplacian flow on $M \times [0, T_0]$
with $\sigma(0)=\sigma_1$.  There hold
\ba \label{estimate1}
\|\sigma-\sigma_0\|_{C^0}\le \min\{\epsilon_0, 
\hat \delta_0(\sigma_0)\}
\ea
and
\ba \label{estimate2}
\|\sigma-\sigma_0\|_{\mc{C}^{4+\mu, (4+\mu)/2}}\le c.
\ea
\end{lem} 
\Pf The bounds for $\lambda_0(\sigma_0)$ and 
$\hat \delta_0(\sigma_0)$ follow from simple compactness
arguments.  The unique existence of the solution $\sigma=\sigma(t)$ 
on a uniform time interval $[0, T_0]$ and the estimates
(\ref{estimate1}) and (\ref{estimate2}) follow from 
Theorem \ref{gauge1} and Theorem \ref{gauge2}. Note that, by Lemma \ref{nulemma}, 
$\|\nu_{\sigma_1, \sigma_0}\|_{C^{2+\mu}}$ can be estimated in terms of $\hat \sigma_0, K$ and $\mu$, 
and $[\nu_{\sigma_1, \sigma_0}]_{\mu}$ can be estimated 
in terms of $\|\sigma_1-\sigma_0\|_{C^{2+\mu}}$ multiplied by a positive constant depending only on 
$\hat \sigma_0$ and $K$.  Hence the condition (\ref{A-condition2}) in Theorem \ref{gauge1} follows from 
(\ref{con5}) and a suitable choice of $T_0$. 

We sketch an alternative argument for obtaining the above uniform existence and estimates.  
Since $\sigma_0$ is torsion-free, there holds $\Delta_{\sigma_0}\sigma_0=0$. On the other hand, we obviously have $\Phi_{\sigma_0}(\sigma_0)=0$. Hence 
$\sigma_0$ is a constant valued solution of the 
$\sigma_0$-gauged Laplacian flow on $[0, \infty)$.
Now we choose e.g. $T_0=1$. (In this argument we can 
choose the value of $T_0$ first, and then determine other constants.) Applying the inverse function theorem at $\sigma_0$ (restricted to the time interval 
$[0, T_0]$) we then obtain the desired existence and 
estimates, provided that $\|\sigma_1-\sigma_0\|_{\mc{C}^{4+\mu}}$ is sufficiently small.  (So we obtain this way a somewhat weaker result than the above one.) \qed \\


Now we formulate the long time existence and convergence 
theorem for the gauge fixed Laplacian flow. Its proof will be presented in the next subsections.

\begin{theo} \label{gaugedconverg} Let $0<\mu<1$, $K>0$, and $\hat \sigma_0$ be a $C^{4+\mu}$  $G_2$-structure on 
$M$. Then there are positive constants $c$ and $\varepsilon_0$ depending only on $\hat \sigma_0,\mu$ 
and $K$ with the following properties.  Let $\sigma_0$ and $\sigma_1$ be two cohomologous closed $C^{4+\mu}$ 
$G_2$-structures on $M$ satisfying (\ref{con1}),
(\ref{con2}), (\ref{con3}), (\ref{con4}) and (\ref{con5}). Assume that $\sigma_0$ is torsion-free. In addition, assume that
\ba \label{con6}
\int_M |\sigma_1-\sigma_0|^2 \le 
\varepsilon_0,
\ea 
where the metric $g_{\sigma_0}$ is used for the norm
and the volume form. (The notation for the volume form is 
omitted. The reader is also advised to be aware that 
this $\varepsilon_0$ is different from $\epsilon_0$ of 
Lemma \ref{simple}.)
Then the $\sigma_0$-gauged Laplacian flow 
\ba \label{gaugedflow1}
\f{\p \sigma}{\p t}=\Delta_{\sigma} \sigma
+d(X_{\sigma_0}(\sigma-\sigma_0)\lr \sigma)
\ea
with initial value $\sigma_1$ has a unique 
${\mathcal C}^{4+\mu, (4+\mu)/2}$ solution 
$\sigma=\sigma(t)$ on $[0, \infty)$ which  converges in $\mc{C}^{4+\mu, (4+\mu)/2}$ to $\sigma_0$
at exponential rate as $t\rightarrow \infty$. 

If $\sigma_0 \in \mc{C}^{l}$ for a non-integer 
$l>4$, then $\sigma(t)-\sigma_0$ converges in 
$\mc{C}^{l+1, (l+1)/2}$ to $0$ at exponential rate.
\end{theo}
\vspace{1mm}

\subsection{$L^2$-decay}
\hspace{1cm}
\vspace{1mm}

Let $0<\mu<1$ and $K>0$ be given. Consider $\hat \sigma_0, \sigma_0$ and $\sigma_1$ satisfying the conditions of Theorem \ref{gaugedconverg}.
Let $\sigma=\sigma(t)$ be the solution of the 
$\sigma_0$-fixed Laplacian flow on $[0, T_0]$ with 
$\sigma(0)=\sigma_1$ as given by Lemma \ref{starting}.
We proceed to prove that $\sigma$ extends to a solution of the $\sigma_0$-gauged Laplacian flow on $[0, \infty)$ and converges to $\sigma_0$ as $t\rightarrow \infty$.

Henceforth we employ the metric $g_{\sigma_0}$ for 
all geometric meaurements and operations. 

\begin{lem} \label{8lemma} There holds 
\ba
\label{88}
\|\sigma-\sigma_0\|_{C^0(M\times [\f{T_0}{8}, T_0])}^2
\le C\int_{M} |\sigma_1-\sigma_0|^2
\ea
for a positive constant $C=C(\hat \sigma_0, K)$.
\end{lem}

The proof of this lemma will be given below. \\

\noindent {\bf Definition 7.1} For $0<\epsilon\le \min\{\epsilon_0, \delta_0^*\}$ let $I_{\epsilon}$ denote the set of $T\ge T_0$ such that $\sigma$ extends to a $C^{4+\mu, (4+\mu)/2}$ solution of the $\sigma_0$-gauged 
Laplacian flow on $[0, T]$, with the following three properties
\ba \label{b1}
\|\sigma-\sigma_0\|_{C^0(M\times [0, T])}\le \min\{\epsilon_0, \hat \delta_0(\sigma_0)\},
\ea
\ba \label{b2}
\|\sigma-\sigma_0\|_{C^0(M\times [\f{1}{8}{T_0}, T])}
\le \epsilon,
\ea
and 
\ba \label{b3}
\|\sigma-\sigma_0\|_{\mc{C}^{4+\mu, (4+\mu)/2}(M \times [t-\f{7T_0}{8}, t])}\le 2c
\ea
for all $T_0\le t \le T$, where $c=c(\hat \sigma_0, \mu, K)$ is from Lemma
\ref{starting}. \\
\vspace{1mm}

\noindent {\bf Remark} Alternatively, we can replace
the condition (\ref{b3}) by the following one
\ba \label{b33}
\|\sigma(t)-\sigma_0\|_{C^{4+\mu}}\le 2c
\ea
for all $T_0\le t \le T$.  Then the proof below also goes through 
with some modifications. \\ 

Set $\theta=\sigma-\sigma_0$ and $\theta_0=\sigma_1-\sigma_0$. Note that (\ref{b2})
and Lemma \ref{starting} imply the following estimate
\ba \label{b4}
\max_{M\times [0, T]} \{|\nabla \theta|,
|\nabla^2 \theta|, |\nabla^3 \theta|, 
|\nabla^4 \theta|, |\f{\p \theta}{\p t}|,
|\f{\p^2 \theta}{\p t^2}|, |\nabla \f{\p \theta}{\p t}|,
|\nabla^2 \f{\p \theta}{\p t}|\} \le 2c.
\ea
 We derive various exponential decay estimates for $\theta$, starting with the $L^2$-dcay.

Since $\sigma_0$ is torsion-free,
we deduce from (\ref{newflow}) and (\ref{phi-formula})
\ba \label{newflow1}
\f{\p \theta}{\p t}=-\Delta \theta
-d (\Phi_{\sigma_0}(\theta))=-dd^* \theta
-d(\Phi_{\sigma_0} (\theta)) 
\ea
with
\ba \label{Aformula}
\Phi_{\sigma_0}(\theta)=A(\sigma_0, \sigma, 
\theta, \nabla\theta).
\ea
(Note that e.g. $\Delta=\Delta_{\sigma_0}$ and 
$\nabla=\nabla_{\sigma_0}$.)
Integrating (\ref{newflow1})
yields 
\ba
\theta=\sigma_1-\sigma_0-d \int_0^t (d^* 
\theta-\Phi_{\sigma_0}(\theta)).
\ea
Since $\sigma_1$ and $\sigma_0$ are cohomologous, it follows that $\theta$ is exact.

We write $\theta(t)=\theta(\cdot, t)$ and often abbreviate it to $\theta$.

\begin{lem} \label{L2decaylemma} There is a positive constant $\epsilon_1=\epsilon_1(\hat \sigma_0, K)$ with the following properties.  Let $T\in I_{\epsilon}$ with $\epsilon\le
\epsilon_1$. 
Then there holds for each $t\in [0, T]$
\ba \label{L2}
\int_M |\theta(t)|^2 \le Ce^{-\lambda_0 t}\int_M |\theta_0|^2.
\ea
for a positive constant $C=C(\hat \sigma_0, K)$.
\end{lem}
\Pf By (\ref{newflow1}), (\ref{Aformula}) and Lemma 
\ref{Alemma} we have 
\ba \label{dt}
&&\f{d}{dt}\int_M |\theta|^2
=2\int_M \theta \cdot (-\Delta_{\sigma_0} \theta-
d(\Phi_{\sigma_0}(\theta))
=-2\int_M |d^* \theta|^2-\int d^* \theta \cdot 
\Phi_{\sigma_0}(\theta)\no
\le -2\int_M |d^*\theta|^2
+C_1\max_t |\theta| \int_M |\nabla \theta|^2,
\ea
where $C_1$ is a universal positive constant. 
Since $d\theta=0$,  we can apply Bochner-Weitzenb\"{o}ck formula and the bound (\ref{b4}) (for controlling 
the curvature) to deduce 
\ba \label{boch}
\int_M |\nabla \theta|^2 \le\int |d^* \theta|^2
+C_2 \int_M |\theta|^2
\ea
with $C_2=C_2(\hat \sigma_0, K)$.
We set
\ba
\epsilon_1=\min\{\f{1}{C_1}, \f{\lambda_0}{C_2C_1}\}.
\ea
Then we deduce from (\ref{dt}) and (\ref{boch}), on account of the bound (\ref{b2}) and the assumption $\epsilon\le \epsilon_1$ 
\ba
\f{d}{dt}\int_M |\theta|^2 &\le& -(2-C_1\epsilon)\int |d^*\theta|^2
+C_1C_2\epsilon\int_M |\theta|^2 \nonumber \\
&\le& -((2-C_1\epsilon)\lambda_0-C_1C_2\epsilon) \int_M |\theta|^2 \nonumber \\
&\le& -\lambda_0 \int_M |\theta|^2,
\ea 
as long as $T_0/8\le t \le T$.  Consequently, we have
for $t \in [T_0/8, T]$
\ba \label{halfway}
\int_M |\theta|^2 \le e^{-\lambda_0 (t-T_0/8)}
\int_M |\theta(T_0/8)|^2.
\ea

To handle the time interval $[0, T_0/8]$  we argue as follows. By the computation in (\ref{dt}) and the bound 
(\ref{b4}) we have 
\ba
\f{d}{dt} \int_M |\theta|^2 \le -2\int_M |d^*\theta|^2
+C_3\int_M |\theta| \cdot |\nabla \theta|
\ea
for $C_3=C_3(\hat \sigma_0, K)$. Employing this inequality, (\ref{boch}) and the Cauchy-Schwarz 
inequality we then deduce for all $0\in [0, T]$
\ba
\f{d}{dt} \int_M |\theta|^2 \le
-\int_M |d^*\theta|^2+C_4 \int_M |\theta|^2
\ea
for $C_4=C_4(\hat \sigma_0, K)$. It follows that 
\ba \label{smalltime}
\int_M |\theta(t)|^2 \le C_5 \int_M |\theta_0|^2
\ea
for $C_5=C_5(\hat \sigma_0, K)$ and $0\le t \le T_0$. 
Combining (\ref{halfway}) and (\ref{smalltime}) we then arrive at (\ref{L2}).
\qed\\

\noindent {\bf Remark} We'll derive decay estimates for other quantities from the above $L^2$-decayof $\theta$. Alternatively, one can also adapt the above arguments to handle other quantities, as they all can be handled in terms of exact forms. However, that approach involves 
additional or stronger conditions for the initial date, which is not satisfactory. 

\subsection{$C^0$-decay and gradient $C^0$-decay}
\hspace{1cm}
\vspace{1mm}

Next we derive decay estimates for $\|\theta(t)\|_{C^0}$
and $\|\nabla \theta(t)\|_{C^0}$. 

\begin{lem} \label{C0decaylemma}  Let $T\in I_{\epsilon}$ with $\epsilon\le
\epsilon_1$. 
Then there holds for each $t\in [T_0/8, T]$
\ba \label{C0}
\|\theta(t)\|_{C^0}^2 \le C e^{-\lambda_0 t}
\int_M |\theta_0|^2.
\ea
Moreover, there holds for each $t\in (0, T]$
\ba \label{dL2}
\int_{t^*}^t \int_M |\nabla \theta|^2 \le 
C e^{-\lambda_0 t} \int_M |\theta_0|^2,
\ea
where $t^*=\max\{t-T_0, 0\}$.
\end{lem} 
\Pf Applying Bochner-Weitzenb\"{o}ck formula we deduce 
\ba \label{evolution}
\f{\p \theta}{\p t}&=&-\nabla^* \nabla \theta
-{\mathcal R}(\theta)-d (\Phi_{\sigma_0}(\theta))
\nonumber \\
&=&-\nabla^* \nabla \theta
-{\mathcal R}(\theta)-d(A(\sigma_0, \sigma, \theta, \nabla \theta)).
\ea
There holds
\ba \label{dA}
d(A(\sigma_0, \sigma, \theta, \nabla \theta)) &=&
\sum_i e^i \wedge \nabla_{e_i} A(\sigma_0, \sigma, 
\theta, \nabla \theta) \nonumber \\ 
&=&\sum_i e^i \wedge (\nabla_{e_i})_2 
A(\sigma_0, \sigma, \theta, \nabla \theta)
+\sum_i e^i \wedge A(\sigma_0, \sigma, \nabla_{e_i} \theta, \nabla \theta)
\nonumber \\&&+ \sum_i e^i \wedge 
A(\sigma_0, \sigma, \theta, \nabla_{e_i} \nabla \theta),
\ea
where $e_i$ stands for a local orthonormal frame and $e^i$ its dual, and $(\nabla_{e_i})_2$ means to take the covariant 
derivative of $A(\sigma_0, \sigma, \theta, \nabla \theta)$ with the second argument $\sigma$ as the variable, while keeping the other arguments parallel.
Then we infer, on account of the bounds (\ref{b1}) 
and (\ref{b4}) 
\ba \label{newb0}
\f{\p }{\p t}|\theta|^2 &\le& -\Delta |\theta|^2
-2|\nabla \theta|^2
+C_6|\theta|^2 + C_6(|\nabla \theta| |\theta|^2
+|\theta| |\nabla \theta|^2+|\theta|^2 |\nabla^2 \theta|),
\ea
where $C_6$ depends only on $\hat \sigma_0, \mu$
and $K$.  There holds $C_6 |\theta| 
|\nabla \theta|^2\le |\nabla \theta|^2+\f{1}{4}
C_6^2|\theta|^2 |\nabla \theta|^2$. 
Applying this and  (\ref{b4}) we then infer
\ba \label{newb}
\f{\p }{\p t}|\theta|^2 &\le& -\Delta |\theta|^2-|\nabla \theta|^2
+C_7|\theta|^2
\ea
for $t\in [T_0/8, T]$, with $C_7
=C_7(\hat \sigma_0, \mu, K)$.  
 Applying Theorem \ref{moser} to $|\theta|^2$ over the interval 
$[t-T_0/8, t]$ (with $t-T_0/8$ as the new time origin)
and appealing to the bound (\ref{b4}) 
we then obtain for $t\in [T_0/8, T]$
\ba \label{4-2}
|\theta(t)|^2 &\le& C_{8} \int_{t-T_0/8}^{t} \int_M |\theta|^2
\ea
with $C_{8}=C_{8}(\hat \sigma_0, K)$.  Combining (\ref{4-2}) and 
Lemma \ref{L2decaylemma} we then arrive at (\ref{C0}). 

Integrating (\ref{newb}) we 
infer (\ref{dL2}).  
\qed \\

\noindent {\bf Remark} 1) The differential inequality 
(\ref{evolution}) is not strong enough to lead to   
\ba \label{similar}
\f{\p}{\p t}|\theta| \le -\Delta |\theta|+C|\theta|
\ea
for a positive constant $C$. This is because of the second order part contained in the term $-d(A(\sigma_0,
\sigma, \theta, \nabla \theta))$. A differential inequality like (\ref{similar}) would allow one to 
obtain the estimate (\ref{C0}) in terms of the $L^2$-version (\ref{moser2}) of the maximum principle, 
which is weaker than the $L^1$-version (\ref{moser}).
\\ 
2) If we apply the $L^2$-version (\ref{moser2}) instead of the $L^1$-version (\ref{moser}) to (\ref{newb}), then
we would obtain the following estimate 
\ba\label{C0}
\|\theta(t)\|_{C^0}^4 \le C e^{-\lambda_0 t}
\int_M |\theta_0|^2.
\ea
This estimate is enough for deriving the convergence of the gauge fixed Laplacian flow and then the convergence of the Laplacian flow. However, it only leads to 
a H\"{o}lder continuity of the limit map of the Laplacian flow. \\

The above remarks also apply to the similar situations below. Now it is convenient to present the proof of Lemma 
\ref{8lemma}. \\

\noindent {\it Proof of Lemma \ref{8lemma}} \,\,Here we deal with the solution $\sigma=\sigma(t)$ on $[0, T_0]$ given by Lemma \ref{starting}.  First, 
arguing as in the proof of Lemma \ref{L2decaylemma}, 
using (\ref{estimate1}) and (\ref{estimate2}) instead of 
(\ref{b1}) and (\ref{b4}), 
we deduce (\ref{smalltime}), with a new $C_5$, which of the same nature as before.  Second, arguing as in the proof of Lemma \ref{C0decaylemma} and applying (\ref{estimate1}) and (\ref{estimate2}) instead of (\ref{b1}) and (\ref{b4}) we also deduce 
(\ref{newb}), with a different $C_8$ which is of the same nature as before. It is clear that we can then apply Theorem \ref{moser} as in the proof of Lemma 
\ref{C0decaylemma} to obtain (\ref{88}). \qed \\

Next we derive a $C^0$-decay estimate for $\nabla \theta$.

\begin{lem} \label{dC0lemma} Let $T \in I_{\epsilon}$ with $\epsilon\le \epsilon_1$. Then there holds
\ba \label{dC0}
\max_t |\nabla \theta|^2 \le Ce^{-\lambda_0t} \int_M |\theta_0|^2
\ea
for $t\in [T_0/8, T]$ and $C=C(\hat \sigma_0,
\mu, K)$.
\end{lem}
\Pf We take the covariant derivative in the equation 
(\ref{evolution}) to obtain
\ba
\f{\p \gamma}{\p t}=-\nabla^* \nabla \gamma
-{\mathcal R}_2 \gamma-\nabla {\mathcal R}_1 \theta
-\nabla d(A(\sigma_0, \sigma, \theta, \gamma))
\ea
for $\gamma=\nabla \theta$, where $\mc{R}_1$ and $\mc{R}_2$ are some linear actions of the curvature operator.  
Employing (\ref{dA}) we obtain a similar formula for 
$\nabla d(A(\sigma_0, \sigma, \theta, \nabla \theta))$.
On account of the bounds (\ref{b1}) and (\ref{b4}) we then deduce
\ba
|\nabla d(A(\sigma_0, \sigma, \theta, \gamma))|
&\le& C_{9}|\theta|(|\gamma|+|\nabla \gamma|+
|\nabla^2 \gamma|)
+C_{9} |\gamma|(|\gamma|+|\nabla \gamma|)
\nonumber \\
&\le& C_{10}(|\gamma|+|\theta|) 
\ea
for positive constants $C_{9}$ and $C_{10}$ depending only 
on $\hat \sigma_0, \mu$ and $K$. 
It follows that 
\ba \label{nabla}
\f{\p}{\p t} |\gamma|^2 &\le& -\Delta |\gamma|^2
-2|\nabla \gamma|^2+C_{11}(|\gamma|^2+|\theta|^2)
\ea
with $C_{11}=C_{11}(\hat \sigma_0, \mu, K)$. 
The extra term $C_{11}|\theta|^2$ in this  differential inequality can be handled by various means. One way is to
combine (\ref{newb}) and (\ref{nabla}) 
to deduce
\ba \label{combine1}
\f{\p}{\p t}(|\theta|^2+|\gamma|^2)
&\le& -\Delta(|\theta|^2+|\gamma|^2)
-(|\gamma|^2+|\nabla^2 \theta|^2)+(C_8+C_{11})(|\theta|^2+|\gamma|^2).
\ea
Applying Theorem \ref{moser}, (\ref{combine1})  and 
the integral estimate (\ref{dL2}) as before we deduce 
(\ref{dC0}).\qed \\

\subsection{$L^2$-decay and $C^0$-decay for $\f{\p \theta}{\p t}$}
\hspace{1cm}
\vspace{1mm}

Let $T \in I_{\epsilon}$ with $\epsilon\le \epsilon_1$. Integrating (\ref{nabla}) and (\ref{combine1}) and employing the previous $L^2$-decay  estimates for 
$\theta$ and $\nabla \theta$ we deduce 
\ba
\int_{T_0/16}^{t} \int_M |\nabla^2 \theta|^2 \le C_{12}e^{-\lambda_0t} \int_M |\theta_0|^2
\ea
 for $t \le [T_0/16, T]$ and $C_{12}=C_{12}(\hat \sigma_0, \mu, K)$. (Here we employ a simple cut-off function of $t$ similar to the one used in the proof of Lemma \ref{interiorlemma}.) Employing this estimate, the above estimates and the evolution equation (\ref{newflow1}) we 
then infer
\ba \label{tL2}
\int_{T_0/16}^{t} \int_M |\f{\p \theta}{\p t}|^2 \le Ce^{-\lambda_0t} \int_M |\theta_0|^2
\ea
for $t \in [T_0/16, T]$ and $C=C(\hat \sigma_0, 
\mu, K)$.\\

\begin{lem} \label{tC0decaylemma} Let $T\in I_{\epsilon}$ with $\epsilon\le \epsilon_1$. Then there holds 
\ba \label{tC0}
\max_t |\f{\p \theta}{\p t}|^2 \le Ce^{-\lambda_0t} \int_M |\theta_0|^2
\ea
for $T_0/8\le t \le T$ and $C=C(\hat \sigma_0, \mu, K)$.
\end{lem}
\Pf  Set $\vartheta\equiv \f{\p \theta}{\p t}$.  
Taking the time derivative in the evolution equation 
(\ref{newflow1})  we deduce
\ba
\f{\p \vartheta}{\p t}&=&-\nabla^* \nabla 
\vartheta-{\mathcal R} \vartheta
-d(\f{\p}{\p t} \Phi_{\sigma_0}(\theta)).
\ea
There holds 
\ba \label{t-der1}
\f{\p}{\p t} \Phi_{\sigma_0}(\theta)=
D_2A(\sigma_0, \sigma, \theta, \gamma)(\vartheta)
+A(\sigma_0, \sigma, \vartheta, \gamma)+
A(\sigma_0, \sigma, \theta, \nabla \vartheta)).
\ea
Following the pattern of computations in (\ref{dA})  we deduce 
\ba \label{t-der2}
|d(\f{\p}{\p t} \Phi_{\sigma_0}(\theta))|
&\le& C_{13}|\theta|(|\gamma| |\vartheta|
+|\nabla \gamma||\vartheta|
+|\gamma| |\nabla \vartheta|+|\nabla \vartheta|
+|\nabla^2 \vartheta|)+C_{13}|\gamma|(|\gamma| 
|\vartheta|+|\vartheta|+|\nabla \vartheta|)\nonumber \\
&&+C_{13}
(|\vartheta||\nabla \gamma|+|\nabla \theta| |\nabla 
\vartheta|)
\ea
for $C_{13}=C_{13}(\hat \sigma_0, \mu, K)$.
Employing (\ref{t-der1}), (\ref{t-der2}) and the bound (\ref{b4}) we then deduce 
\ba
\f{\p}{\p t}|\vartheta|^2\le -\Delta |\vartheta|^2
-|\nabla \vartheta|^2+C_{14}|\vartheta|(|\vartheta|
+|\gamma|+|\theta|)
\ea
for $C_{14}=C_{14}(\hat \sigma_0, \mu, K)$. 
Combining this with (\ref{combine1}) we then infer 
\ba \label{combine2}
\f{\p}{\p t}(|\vartheta|^2+|\gamma|^2+
|\theta|^2)&\le& -\Delta (|\vartheta|^2+|\gamma|^2+
|\theta|^2)
-(|\nabla \vartheta|^2+|\nabla \gamma|^2+|\nabla \theta|^2)\nonumber \\ &&+C_{15}(|\vartheta|^2+|\gamma|^2+|\theta|^2)
\ea
for $C_{15}=C_{15}(\hat \sigma_0, \mu, K)$.
Applying this differential inequality, Theorem \ref{moser}, the integral estimate (\ref{tL2}) 
and the above $L^2$-decay estimates for $\theta$ and $\gamma$  we 
arrive at (\ref{tC0}). \qed\\

\subsection{$\mc{C}^{4+\mu, (4+\mu)/2}$-decay}
\hspace{1cm}\vspace{1mm}

\begin{lem} \label{4lemma} Let $T\in
I_{\epsilon}$ with $\epsilon\le \epsilon_1$. Assume $T_0 \le t \le T$. Then there holds
\ba \label{4decay}
\|\theta\|^2_{{\mathcal C}^{4+\mu, (4+\mu)/2}(M\times [t-\f{3}{4}T_0,
t])} &\le& Ce^{-\lambda_0 t}\int_M |\theta_0|^2,
\ea
with $C=C(\hat \sigma_0, \mu, K)$. 
\end{lem}
\Pf We view the evolution equation (\ref{newflow1}) with the given solution $\theta$ as a linear equation in the form of (\ref{closedheat}), similar to (\ref{linearview}).
Thus we have 
\ba
\f{\p \theta}{\p t}+\Delta \theta+d(\Phi_0(\theta)+\Phi_1(\nabla \theta))=0
\ea
with $\Phi_0\equiv 0$ and $\Phi_1(\nabla \theta)=
A(\sigma_0, \sigma, \theta, \nabla \theta)$. By Lemma 
\ref{Alemma} and the bound (\ref{b1})  there holds 
\ba
\|\Phi_1\|_{C^0(M\times [0, T])}\le C_0  \epsilon
\le C_0 \bar \delta_0^* \le \delta_0(\sigma_0).
\ea
Now consider  $T_0\le t \le T$. By the bounds (\ref{b1}) and (\ref{b3}) we have 
\ba
\|\Phi_1\|_{\mc{C}^{3+\mu, (3+\mu)/2}(M \times [t-\f{7T_0}{8}, t])}\le C_{16}
\ea
for $C_{16}=C_{16}(\hat \sigma_0, K)$. 
Applying Theorem
\ref{interiortheorem} with $l=4+\mu, m=\mu, \epsilon=\f{1}{8}T_0$ and 
$t-\f{7}{8}T_0$ as the new time origin we then arrive at 
\ba
\|\theta\|_{{\mathcal C}^{4+\mu, (4+\mu)/2}(M\times [t-\f{3T_0}{4},
t])} &\le&
C_{17}\|\theta\|_{{\mathcal C}^{\mu, \mu/2}(M\times [t-\f{7}{8}T_0,
t])}
\ea
for $C_{17}=C_{17}(\hat \sigma_0, K)$. 
Combining this with Lemma \ref{C0decaylemma}, 
\ref{dC0lemma} and \ref{tC0decaylemma} we then 
arrive at 
(\ref{4decay}). \qed \\

\noindent {\it Proof of Theorem \ref{gaugedconverg}}

We define the number $\varepsilon_0$ in the theorem as follows 
\ba \label{e0}
\varepsilon_0=(4C)^{-1}\min\{\epsilon_1^2,
c^2, K^2,
\epsilon_0^2, (\hat\delta_0^*)^2, \rho^2\},
\ea
where $C$ is the larger of the $C$ from Lemma \ref{8lemma} and the $C$ from Lemma \ref{4lemma}, and $\rho=\rho(\hat \sigma_0, K)$ is from Lemma \ref{starting}. \\

\noindent {\it Claim 1} The set $I_{\epsilon_1}$ is
nonempty, indeed $T_0\in I_{\epsilon_1}$. \\

Indeed, the estimates (\ref{estimate1}) and (\ref{estimate2}) imply the conditions (\ref{b1}) and 
(\ref{b3}) for $T=T_0$, and 
Lemma \ref{8lemma}, the assumption (\ref{con6}) and 
(\ref{e0}) imply the condition (\ref{con2}) for 
$\epsilon=\epsilon_1$ and $T=T_0$. It follows that
$T_0 \in I_{\epsilon_1}$. \\

\noindent {\it Claim 2} The set $I_{\epsilon_1}$ is closed.  \\

This follows from elementary convergence 
and continuity arguments based on (\ref{con3}). \\

\noindent {\it Claim 3} The set $I_{\epsilon_1}$ is 
open in $[T_0, \infty)$. \\

To prove this claim, assume $T\in I_{\epsilon_1}$. 
By Lemma \ref{4lemma}, the assumption (\ref{con6})
and (\ref{e0}) we have for $T_0\le t \le T$
\ba \label{tt}
\|\theta\|_{{\mathcal C}^{4+\mu, (4+\mu)/2}(M\times [t-\f{3T_0}{4},
t])} \le \f{1}{2}\min\{\epsilon_1,
c, K,
\epsilon_0, \hat\delta_0^*, \rho\}.
\ea
Applying Lemma \ref{starting} to the initial 
$G_2$-structure $\sigma(T-\f{1}{2}T_0)$ with $T-\f{1}{2}T_0$ as the time origin we then obtain a $\mc{C}^{4+\mu, (4+\mu)/2}$-solution of the $\sigma_0$-gauged Laplacian 
flow on $[T-\f{1}{2}T_0, T+\f{1}{2}T_0]$. By  its uniqueness 
property, it agrees with $\sigma(t)$ on $[T-\f{1}{2}T_0,
T]$. Hence it extends $\sigma(t)$ to a $\mc{C}^{4+\mu,
(4+\mu)/2}$ solution of the $\sigma_0$-gauged 
Laplacian flow on $[0, T+\f{1}{2}T_0]$.  By 
(\ref{tt}) and continuity we have 
\ba
\label{ttt}
\|\theta\|_{{\mathcal C}^{4+\mu, (4+\mu)/2}(M\times [t-\f{3T_0}{4},
t])} \le c
\ea
for all $t\in [T_0, T']$, whenever $T'>T$ and 
$T'-T$ is sufficiently small. For such a $T'$ and 
a $t$, there are two possible cases to consider. One is that $t-\f{7T_0}{8} \ge T_0$, the other is that 
$t-\f{7T_0}{8}<T_0$. In th first case, we write the time interval $[t-\f{7T_0}{8}, t]$ as the union of 
$[t-\f{7T_0}{8}, t-\f{3T_0}{8}]$ and $[t-\f{3T_0}{8}, t]$.  Then we can apply the estimates (\ref{ttt}) 
to the both subintervals.  By the triangular 
inequality we then deduce
\ba
\label{tttt}
\|\theta\|_{{\mathcal C}^{4+\mu, (4+\mu)/2}(M\times [t-\f{7T_0}{8},
t])} \le 2c.
\ea
In the latter case, we write $[t-\f{7T_0}{8}, t]$ as the union of $[t-\f{7T_0}{8}, T_0]$ and $[T_0, t]$. Then we can apply the estimate (\ref{ttt}) to the second subinterval, while apply the estimate (\ref{estimate2}) in Lemma 
\ref{starting} to the first subinterval. By the triangular inequality we again arrive at (\ref{tttt}).   We conclude that  
$T'$ belongs to $I_{\epsilon_1}$, whenever $T'>T$ and 
$T'-T$ is sufficiently small.   
It follows that $I_{\epsilon_1}$ is open. 

Combining the above three claims we infer that 
$I_{\epsilon_1}=[T_0, \infty)$. Hence the solution 
$\sigma(t)$ has been extended to $[0, \infty)$.
By Lemma \ref{4lemma}, $\theta$ converges  
in $\mc{C}^{4+\mu, (4+\mu)/2}$ to zero at exponential rate as $t\rightarrow \infty$. 

Finally we assume that $\sigma_0 \in \mc{C}^{l}$ for 
a non-integer $l>4$. Then $\sigma-\sigma_0\in \mc{C}^{l+1, (l+1)/2}$ by the regularity property provided by 
Theorem \ref{gauge1}. By the above $\mc{C}^{4+\mu,
(4+\mu)/2}$ convergence we have 
\ba
\|\Phi_1\|_{C^0(M\times [t-\f{7T_0}{8}, t])}\le 
\delta_0(l+1, \sigma_0),
\ea
whenever $t$ is large enough. Hence we can argue as 
in the proof of Lemma \ref{4lemma} to obtain for large 
$t$ an estimate for $\|\theta\|_{\mc{C}^{l+1, (l+1)/2}(M\times 
[t-1, t])}$ similar to (\ref{4decay}). This is precisely
the desired convergence of $\theta$ in $\mc{C}^{l+1, (l+1)/2}$
to zero at exponential rate. (In particular, since $\sigma_0 \in C^{4+\mu, (4+mu)/2}$, $\sigma-\sigma_0$ 
converges in $\mc{C}^{5+\mu, (5+\mu)/2}$ to 
$0$ at 
exponential rate.) \qed \\ 

\sect{Long time existence and convergence of the Laplacian flow}
\vspace{2mm}

\subsection{Local structure of the moduli space of torsion-free $G_2$-structures} 
\hspace{1cm}\vspace{1mm}

To establish the long time existence and convergence of the Laplacian flow starting near a torsion-free 
$G_2$-structure, we'll need some results on the local structure of the moduli space of torsion-free 
$G_2$-structures. On the other hand, the said convergence of the Laplacian flow also reveals  an interesting dynamic property of this moduli space.
Consider a non-integer $l>1$. (In the following presentation, $l$ is allowed to be $\infty$ with the convention $\infty+1=\infty$, except for the norms and their associated objects.) Let $\mc{T}_l$ denote the space of torsion-free $C^l$ $G_2$-structures on $M$, and $\text{Diff}\,^{l+1}_0(M)$ be the group of $C^{l+1}$ diffeomorphisms of $M$ which are $C^{l+1}$-isotopic to the identity map.  Obviously, this group acts on $\mc{T}_l$. The quotient $\mc{T}_l/\text{Diff}\,^{l+1}_0(M)$ is the moduli space of torsion-free $C^l$ $G_2$-structures. Let $\pi_l:
\mc{T}_l \rightarrow  \mc{T}_l/\text{Diff}\,^{l+1}_0(M)$ be the projection. Occationally, we abbreviate $\mc{T}_{\infty}$, $\text{Diff}\,^{\infty}_0(M)$ and $\pi_{\infty}$ to $\mc{T}, \text{Diff}\,_0(M)$ and 
$\pi$ respectively.

Let $\sigma_0$ be a given torsion-free $C^l$  $G_2$-structure  on $M$. In 
this subsection, all the geometric operations and measurements are  w.r.t.~$g_{\sigma_0}$.
For $0\le j\le 7$ let ${\mathcal H}^j \equiv {\mathcal H}_{\sigma_0}^j(M)\subset 
C^l_o(\Lambda^3T^*M)$ denote the 
space of harmonic $j$-forms on $M$ w.r.t. $\sigma_0$, which represents the DeRham
cohomology group $H^j(M, \mathbb{R})$. \\

\noindent {\bf Definition 9.1}
For $l, r>0$ and $\gamma \in C_0^l(\Lambda^3 T^*M)$,
let $B^l_r(\gamma)$ denote the open ball of center 
$\gamma$ and radius $r$ in $C_o^l(\Lambda^3 T^*M)$.
We set $\mc{B}_r^l(\gamma)=B^l_{r}(\gamma) \cap \mc{H}^3$ for $\gamma \in \mc{H}^3$. 

\begin{theo} \label{moduli} Let $0<\mu<1$, $2+\mu\le l \le \infty$ (a non-integer), and let $\sigma_0$ be  a given torsion-free $C^l$ $G_2$-structure on $M$. 
Then there are a positive number $r_0 \equiv r_0( \sigma_0, \mu)\le \epsilon_0$ (with $\epsilon_0$ from Lemma \ref{simple}) depending only on $\sigma_0$ and $\mu$, and a smooth embedding $\Xi_{\sigma_0}: \mc{B}^{2+\mu}_{r_0}(\sigma_0)
 \rightarrow C^l_0(\Lambda^3T^*M)$  whose image consists of 
torsion-free $G_2$-structures, such that $\Xi_{\sigma_0}(\gamma)$ is cohomologous to $\gamma$ for all $\gamma \in
\mc{B}^{2+\mu}_{r_0}$. (Since $r_0\le \epsilon_0$, $\mc{B}^{2+\mu}_{r_0}$ consists of $G_2$-structures.) Moreover, $\Xi_{\sigma_0}(\mc{B}^{2+\mu}_{r_0})$ provides a local slice of the space of torsion-free 
$C^l$ $G_2$-structures under the action of 
${\text{Diff}\,}_0^{l+1}(M)$. As a consequence, 
the collection of $(\mc{B}^{2+\mu}_{r_0}(\sigma_0), \Xi_{\sigma_0})$ for all $\sigma_0 \in \mc{T}_l$ provides a natural smooth structure on $\mc{T}_l/\text{Diff}\,_0^l(M)$.  

We also have for all $h\in 
\mc{B}^{2+\mu}_{r_0}(\sigma_0)$
\ba \label{h}
\|\Xi_{\sigma_0}(h)-\sigma_0\|_{C^l} \le C\|h-\sigma_0\|_{C^l}
\ea
with a positive constant $C$ depending only on  $\sigma_0$ and $l$.
\end{theo}

This result is a refinement of Joyce's result \cite{Joyce} on local moduli of torsion-free $G_2$-structures, and can be viewed as the elliptic version of Theorem \ref{stability}.
For its proof we refer to [XY2].

We consider 
for each $0\le j\le 7$ the projection map $H_{l}: C^l_o(\Lambda^jT^*M)
\rightarrow {\mathcal H}^j$ which sends 
each closed
$C^{l}$ $j$-form $\gamma$ to the unique 
harmonic form (w.r.t.~$\sigma_0$) in the cohomology class of $\gamma$, $[\gamma]_l=
\{\gamma+d\beta: \beta \text{ is a } C^{l+1} 
(j-1)-\text{form on } M\}$. The following result is a simple consequence of basic elliptic regularity.

\begin{lem} \label{ll} There holds $H_l=H_{l'}$
on $C^{l}_o(\Lambda^jT^*M)$ for $l'\le l$. 
\end{lem}

Set $\hat {\mathcal B}_{r_0}^{2+\mu}(\sigma_0)=
H_{2+\mu}^{-1}(({\mathcal B}_{r_0}^{2+\mu}(\sigma_0))$
and $ \hat {\mathcal B}_{r_0}^{2+\mu, l}(\sigma_0)=
H_l^{-1}({\mathcal B}_{r_0}^{2+\mu}(\sigma_0))=
\hat {\mathcal B}_{r_0}^{2+\mu}(\sigma_0)\cap C^l$.
We define the smooth projection maps
\ba
\hat \Xi_{\sigma_0}=\Xi_{\sigma_0} \circ H_{2+\mu}: \hat {\mathcal B}_{r_0}^{2+\mu}(\hat \sigma_0)
\rightarrow \mc{T}_{2+\mu}
\ea
and 
\ba
\Pi_{\sigma_0}=\pi_{2+\mu} \circ \hat \Xi_{\sigma_0}: \hat {\mathcal B}_{r_0}^{2+\mu}(\hat \sigma_0)
\rightarrow \mc{T}_{2+\mu}/\text{Diff}\,^{2+\mu}_0(M).
\ea
Restricting them to $\hat {\mathcal B}_{r_0}^{2+\mu, l}(\hat \sigma_0)$ we then obtain the smooth projection maps 
\ba
\hat \Xi_{\sigma_0}: \hat {\mathcal B}_{r_0}^{2+\mu, l}(\hat \sigma_0)
\rightarrow \mc{T}_l
\ea
and 
\ba
\Pi_{\sigma_0}: \hat {\mathcal B}_{r_0}^{2+\mu,l}(\hat \sigma_0)
\rightarrow \mc{T}_l/\text{Diff}\,^l_0(M).
\ea
Note that the former equals $\Xi_{\sigma_0} \circ H_{l}$, while the latter equals $\pi_l \circ \hat \Xi_{\sigma_0}$.

\begin{lem} \label{harmonic-bound} There holds for each $1\le j \le 7$ and all $\gamma \in C^l_o(\Lambda^jT^*M)$
\ba \label{harmonic}
\max\{\|H_l(\gamma)-\sigma_0\|_{C^l}, \|H_l(\gamma)-\gamma\|_{C^l}\} \le C\|\gamma-\sigma_0\|_{C^l}
\ea
for a positive constant $C$ depending only on $\sigma_0$ and $l$. 
\end{lem} 
\Pf Let $\beta \in C^{l+1}(\Lambda^2 T^*M)$  be the unique solution of the equation 
\ba
\Delta \beta=-d^* \gamma
\ea
subject to the $L^2$-orthogonality condition $\beta \perp_{L^2} \mathcal{H}^2$.
Set $h=\gamma+d\beta$. Then there holds 
$(d^*+d)h=0$, and hence $h \in
\mathcal{H}^3$. It follows that $h=H_l(\gamma)$. By basic elliptic 
estimates we have 
\ba
\|\beta\|_{C^{l+1}} \le C_1\|d^*\gamma\|_{C^{l-1}}
\ea
for a positive constant $C_1$ and hence 
\ba
\|H_l(\gamma)-\gamma\|_{C^l} \le C_2\|d^*\gamma\|_{C^{l-1}}
\ea
for a positive constant $C_2$, where $C_1$ and $C_2$ depend only on $\sigma_0$ and $l$. Next we observe $d^*\gamma=d^*(\gamma-\sigma_0)$ and hence 
$\|d^*\gamma\|_{C^{l-1}}\le 
C_3\|\gamma-\sigma_0\|_{C^l}$
for a universal positive constant $C_3$. The first estimate in (\ref{harmonic}) follows. The second estimate follows from the first by the triangular inequality. \qed \\

\begin{lem} \label{moduli-bound} 
There holds for all
$\gamma \in \hat {\mc{B}}^{2+\mu,l}_{r_0}(\sigma_0)$ 
\ba \label{Xi}
\max\{\|\hat \Xi_{\sigma_0}(\gamma)-\gamma\|_{C^l},
\|\hat \Xi_{\sigma_0}(\gamma)-\sigma_0\|_{C^l}\}
\le C\|\gamma-\sigma_0\|_{C^l}
\ea
for a positive constant $C$ depending only on $\sigma_0$ and $l$.
\end{lem}
\Pf  The desired estimates follow from Theorem \ref{moduli}, Lemma \ref{harmonic-bound} and the triangular inequality. \qed \\

\subsection{Convergence I} 
\hspace{1cm}
\vspace{1mm}

In this subsection, all the geometric operations and measurements are w.r.t.~$\sigma_0$ as given below. 

\begin{theo} \label{converg1} Let $0<\mu<1$, $K>0$ and 
$\hat \sigma_0$ a given $C^{4+\mu}$ structure on $M$.
 Assume that 
$\sigma_0$ and $\sigma_1$ satisfy the conditions in Theorem \ref{gaugedconverg}. Then there is a unique $C^{2+\mu, (3+\mu)/2}$ solution 
$\sigma=\sigma(t)$ of the Laplacian flow on $[0, \infty)$ which takes the initial
value $\sigma_1$. It converges in $\mc{C}^{2+\mu, (2+\mu)}$ at exponential rate to a torsion-free $C^{2+\mu}$ $G_2$-structure $\sigma_{\infty}$ on $M$ which is 
$C^{3+\mu}$ isotopic to $\sigma_0$.  Moreover, there holds
\ba
\|\sigma_{\infty}-\sigma_0\|_{C^{2+\mu}}\le C\|\sigma_1-\sigma_0\|_{C^{4+\mu}} 
\ea
for a positive constant $C$ depending only on 
$\hat \sigma_0, K$ and $\mu$.
If $\sigma_0 \in C^l$ for an non-integer $l>4$, then 
$\sigma_{\infty} \in C^{l-2}$ and is $C^{l-1}$ isotopic 
to $\sigma_0$, and $\sigma$ converges in 
$\mc{C}^{l-2, (l-3)/2}$ to $\sigma_{\infty}$ at exponential rate.
\end{theo}

\Pf By Theorem \ref{gaugedconverg}, we have a unique $\mc{C}^{4+\mu, (4+\mu)/2}$ solution $ \sigma=\sigma(t)$ of the $\sigma_0$-gauged Laplacian flow on $[0, \infty)$ which takes the initial value $\sigma_1$ and converges in $\mc{C}^{4+\mu, (4+\mu)/2}$ at exponential rate to $\sigma_0$. As in Section 7, we consider the ODE
(\ref{ODE}) on $M\times [0, \infty)$ with the present $\sigma$ and the initial 
condition (\ref{initial}). Since $\theta=\sigma-\sigma_0$ converges to $0$ in $\mc{C}^{4+\mu,
(4+\mu)/2}$ at exponential rate, $X_{\sigma_0}(\sigma-\sigma_0)$
converges to $0$ in $\mc{C}^{3+\mu, (3+\mu)/2}$ at exponential rate. Consequently, the $\mc{C}^{3+\mu, (4+\mu)/2}$ solution $\phi$
of (\ref{ODE}) exists for all time and converges in $\mc{C}^{3+\mu, (4+\mu)/2}$ at exponential rate to a 
$C^{3+\mu}$ map $\phi_{\infty}$.  Taking derivative
in (\ref{ODE}) we infer 
\ba
\nabla_{\f{d}{dt}}d\phi=\nabla_{d\phi}X_{\sigma_0}
\ea
with $X_{\sigma_0}=X_{\sigma_0}(\sigma-\sigma_0)$.
Consequently, there holds for each $p\in M$ and $v\in T_pM$
\ba
\f{d}{dt}|d\phi(v)|^2&=& 2d\phi(v) \cdot \nabla_{d\phi(v)}X_{\sigma_0} \nonumber \\
&\ge&-2 |d\phi(v)|^2 |\nabla X_{\sigma_0}(\phi(p))|.
\ea
Because $|d\phi(v)|=|v|$ at $t=0$, integration then yields 
\ba
|d\phi(v)| \ge |v|e^{-\int_0^t\|\nabla  X_{\sigma_0}\|_{C^0}}
\ea
at time $t$. We conclude
\ba
|d\phi_{\infty}(v)|\ge |v|e^{-\int_0^{\infty}\|\nabla  X_{\sigma_0}\|_{C^0}}
\ea
and hence $d\phi_{\infty}$ is an isomorphism everywhere.
Consequently, $\phi_{\infty}$ is a local diffeomorphism. 
Since it is homotopic to the identity map, it is a diffeomorphism.

Now we set $\hat \sigma(\cdot, t)=\phi(\cdot, t)^*\sigma(\cdot, t)$. By Lemma \ref{goback}, 
$\hat \sigma$ is a $\mc{C}^{2+\mu, (3+\mu)/2}$ solution of 
the Laplacian flow on $[0, \infty)$ with the initial value $\sigma_1$. By the above reasoning and the convergence of $\sigma$ to $\sigma_0$, it converges 
in $\mc{C}^{2+\mu, (3+\mu)/2}$ at exponential rate 
to $\hat \sigma_{\infty}=\phi_{\infty}^*\sigma_0$.
This $G_2$-structure is obviously torsion-free because $\sigma_0$ is so. This also follows from the Laplacian flow equation and the convergence of 
$\f{\p \hat \sigma}{\p t}$ to $0$.  

If $\sigma_0\in C^l$ for $l>4$, then the corresponding 
convergence result in Theorem \ref{gaugedconverg} 
and the above reasoning imply that $\hat \sigma$ converges in $\mc{C}^{l-2, (l-1)/2}$ to $\hat \sigma_{\infty}$ at exponential rate. Moreover, 
$\phi$ converges in $C^{l-1, l/2}$ to $\phi_{\infty}$.
Hence $\sigma_{\infty} \in C^{l-2}$ and is $C^{l-1}$
isotopic to $\sigma_0$. \qed \\

\subsection{Convergence II}
\hspace{1cm}
\vspace{1mm}

\noindent {\bf Definition 9.2} Let $0<\mu<1$, 
$4+\mu \le l\le \infty$ (a non-integer), and  $\gamma\in C^{4+\mu}_o(\Lambda^3T^*M)$.   For $K>0$ and $r>0$ we set 
\ba
B^l_{K,r}(\gamma)=\{\gamma' \in B^l_K(\gamma):
\|\gamma'-\gamma\|_{C^{2+\mu}}\le r, 
\|\gamma'-\gamma\|_{C^{4+\mu}} \le K\}.
\ea
It will be called a {\it $K$-strong (or simply,  strong) $C^{2+\mu}$ ball (or neighborhood) }in the 
$C^{l}$ space.

\begin{theo} \label{converg2} Let $0<\mu<1, K>0$ and $\hat 
\sigma_0$ a given $C^{4+\mu}$ torsion-free $G_2$-structure on 
$M$. There is a positive constant $\hat r_0\le \epsilon_0$ depending only on $K, \mu$ and $\hat \sigma_0$ with 
the following properties. 
Let $\sigma_1 \in B^l_{K,\hat r_0}(\hat \sigma_0)$ for a 
non-integer $l\ge 4+\mu$. 
Then there is a unique $C^{l-2, (l-1)/2}$ solution of the Laplacian flow on $[0, \infty)$ with the initial 
value $\sigma_1$. It converges in $C^{l-2, (l-1)/2}$
at exponential rate to a torsion free $C^{l-2}$ 
$G_2$-structure $\sigma_{\infty}$, which is $C^{l-1}$ isotopic to $\hat \Xi_{\hat \sigma_0}(\sigma_1)$.  
\end{theo} 

\Pf We apply the results in Subsection 9.1 
with $\hat \sigma_0$ playing the role of $\sigma_0$ there. Let $\sigma_1 \in B^l_{K, \hat r_0}(\hat \sigma_0)$ for some $\hat r_0>0$.  By Lemma \ref{harmonic-bound} there hold
\ba
\|H_{l}(\sigma_1)-\hat \sigma_0\|_{C^{4+\mu}, \hat \sigma_0}
\le C_1K
\ea
and
\ba
\|H_l(\sigma_1)-\hat \sigma_0\|_{C^{2+\mu}, \hat \sigma_0}\le C_1 \hat r_0,
\ea
where $C_1$ depends only on $\hat \sigma_0$ and $\mu$.
Assume $\hat r_0\le C_1^{-1}r_0$. Then $H_l(\sigma_1) \in 
\mc{B}^{2+\mu,l}_{r_0}(\hat \sigma_0)$. It follows that
$B^l_{K,\hat r_0}(\hat \sigma_0) \subset \hat 
{\mc{B}}^{2+\mu,l}_{r_0}(\hat \sigma_0)$. Now we set 
$\sigma_0=\hat \Xi_{\hat \sigma_0}(\sigma_1)$. By Lemma \ref{moduli-bound} we have
\ba
\|\sigma_1-\sigma_0\|_{C^{4+\mu}, \hat \sigma_0}
\le C_2K, \,\,\, \|\sigma_0-\hat \sigma_0\|_{C^{4+\mu}, \hat \sigma_0} \le C_2K
\ea
and
\ba
\|\sigma_1-\sigma_0\|_{C^{2+\mu}, \hat \sigma_0}
\le C_2\hat r_0, \,\,\, \|\sigma_0-\hat \sigma_0\|_{C^{2+\mu}, \hat \sigma_0} \le C_2 \hat r_0
\ea
for a positive constant $C_2$ depending only on $\hat \sigma_0$ and $\mu$. We assume $\hat r_0\le C_2^{-1} \epsilon_0$
(with $\epsilon_0$ from Lemma \ref{simple}). Then the
ratios between the norms measured in $\hat \sigma_0$ and the norms measured in terms of $\sigma_0$ (in both directions) are bounded by positive constants depending only on $\hat \sigma_0$ and $\mu$. Hence we have
\ba
\|\sigma_1-\sigma_0\|_{C^{4+\mu}, \sigma_0} \le C_3K
\ea
and 
\ba
\|\sigma_1-\sigma_0\|_{C^{2+\mu}, \sigma_0} \le C_3 \hat r_0
\ea
for a positive constant $C_3$ depending only on 
$\hat \sigma_0$ and $\mu$. 

Now we replace $K$ in Theorem \ref{gaugedconverg} by 
$\max\{C_2, C_3\}K$
and obtain the corresponding $\rho$ there. Then it is 
clear that we can define $\hat r_0$ according to 
the above 
two conditions and the 
conditions in Theorem \ref{gaugedconverg} and Lemma \ref{starting}. Then we can apply Theorem \ref{converg1} to deduce the desired long time existence and convergence of the Laplacian flow with the initial value $\sigma_1$. The claimed 
isotopy property of the limit also follows from the same theorem. \qed \\

\sect{The limit map of the Laplacian flow}

Set $\mc{F}(\sigma_0, \sigma_1)=\sigma_{\infty}$, where $\sigma_{\infty}$ is the limit given in Theorem \ref{converg1}. Note that by the uniqueness part of 
Theorem \ref{existence}, this map is actually independent of $\sigma_0$. We'll keep $\sigma_0$ as an argument for the following reasons. First, in the 
proof below, we'll employ several quantities which  depend on both $\sigma_1$ and $\sigma_0$. So it is natural to treat everything in the framework of two 
arguments $\sigma_0$ and $\sigma_1$. Second, since we construct our arguments without using the uniqueness part of Theorem \ref{existence}, they have a broader scope of possible applications.

\begin{theo} \label{smooth} Let $\hat \sigma_0, \sigma_0$ and $\sigma_1$ be as in Theorem \ref{converg1}. Then the map $\mc{F}(\sigma_0, \sigma_1)$ is a Lipschitz continuous  function on $\sigma_0$ and $\sigma_1$ w.r.t.~ 
$C^{4+\mu}$-norm on $\sigma_0$, $C^{4+\mu}$-norm on 
$\sigma_1$, and $C^{2+\mu}$-norm on $\mc{F}(\sigma_0, 
\sigma_1)$. In general, $\mc{F}(\sigma_0, \sigma_1)$ is Lipschitz continuous  w.r.t.~$C^{l}$-norm on 
$\sigma_0$, $C^l$-norm on $\sigma_1$, and $C^{l-2}$-norm on $\mc{F}(\sigma_0, \sigma_1)$, provided that $\sigma_0 \in C^{l}$ and $\sigma_1\in C^l$ for $l\ge 4+\mu$. 
\end{theo}

\Pf For two initial $G_2$-structures $\sigma_1$
and $\bar \sigma_1$, and two torsion-free reference $G_2$-structures $\sigma_0$ and $\bar \sigma_0$ as in the situation of Theorem \ref{converg1}, we consider the corresponding solutions $\sigma=\sigma(t)$ of the $\sigma_0$-gauged Laplacian flow, and $\bar \sigma= 
\bar \sigma(t)$ of the $\bar \sigma_0$-gauged Laplacian flow.  Set $\theta=\sigma-\sigma_0$, $\bar \theta=\bar \sigma-\bar \sigma_0$, and $\gamma=\bar \theta-\theta$.
We first derive estimates for $\gamma$. There holds
\ba \label{gamma-equation}
\f{\p \gamma}{\p t}&=&-\Delta_{\sigma_0} \gamma
-d(\Phi_{\sigma_0}(\gamma))
+(\Delta_{\sigma_0}-\Delta_{\bar \sigma_0})\bar \theta-d(\Phi_{\bar \sigma_0}(\bar \theta)-\Phi_{\sigma_0}(\bar \theta)).
\ea
We handle the two difference terms on the RHS of 
(\ref{gamma-equation}) in terms of integration. For example, there holds
\ba 
d(\Phi_{\bar \sigma_0}(\bar \theta)-\Phi_{\sigma_0}(\bar \theta))
=d\left((\int_0^1  \f{d}{ds}\Phi_{\sigma_0+s(\bar \sigma_0-\sigma_0)}ds) 
(\gamma)\right).
\ea
Since $\gamma$ is exact, we can apply the arguments in the proof of Theorem \ref{gaugedconverg} to obtain 
decay estimates for  $\gamma$. Here we employ the exponential decay of $\theta$ and $\bar \theta$ to handle the integration terms resulting from the  three difference terms on the RHS of (\ref{gamma-equation}). 
(Note that the non-homogeneous terms arising from 
the last two  terms can be handled by 
elementary integration techniques). We deduce for 
$2<l\le 4+\mu$ 
\ba
\|\gamma\|_{\mc{C}^{l, l/2}(M\times [t-1, t])}^2 \le C_1 e^{-\lambda_0 t}
(\|\bar \sigma_0-\sigma_0\|_{C^l}^2+\|\bar \sigma_1-\sigma_1\|_{C^l}^2)
\ea
with a positive constant $C_1$ depending only on $\hat \sigma_0, \mu$ and $K$. All the constants below in this proof have this same dependence.

Next let $\phi$ be the solution of the ODE (\ref{ODE}) 
corresponding to the solution $\sigma$ (with $\phi(0)=Id$), and let $\bar \phi$ be the solution of the ODE (\ref{ODE}) 
corresponding to the solution $\bar \sigma$, i.e.
\ba
\f{d \bar \phi}{dt}=-X_{\bar \sigma_0}(\bar \theta)(\bar \phi),
\ea
also with $\bar\phi(0)=Id$. Then $\phi^*\sigma$ is the solution of the Laplacian flow with the initial value $\sigma_1$, and $\bar \phi^*\bar \sigma$ is the solution of the Laplacian flow with the initial value $\bar \sigma_1$ as studied in the proof of Theorem \ref{converg1}.  Now we embed $M$ into a 
Euclidean space and set $\psi=\bar \phi-\phi$. We deduce
\ba \label{map-difference}
\f{d \psi}{dt}&=&-(X_{\bar \sigma_0}(\bar \theta)(\phi+\psi)-X_{\sigma_0}(\theta)(\phi))
\nonumber \\
&=& -(X_{\sigma_0}(\theta)(\phi+\psi)-
X_{\sigma_0}(\theta)(\phi))
-(X_{\bar \sigma_0}(\bar \theta)(\bar \phi)-X_{\sigma_0}(\bar \theta)(\bar \phi))
-X_{\sigma_0}(\gamma)(\bar \phi)\nonumber \\
\ea
and $\psi(0)=0$. As above, we can handle the two difference terms in the bottom line of (\ref{map-difference}) 
by intgeration. For example, there holds
\ba
(X_{\bar \sigma_0}(\bar \theta)(\bar \phi)-X_{\sigma_0}(\bar \theta)(\bar \phi))=
\int_0^1 \f{d}{ds}X_{\sigma_0+s(\bar \sigma_0-\sigma_0)}(\bar \theta)
(\bar \phi) ds.
\ea
(We can assume that $\|\bar \sigma_0-\sigma_0\|_{C^0}\le \epsilon_0$. Then $\sigma_0+s(\bar \sigma_0-\sigma_0)$ are $G_2$-structures and as smooth as $\sigma_0$ and $\bar \sigma_0$ for $0\le s \le 1$.)
When treating the first one, we need to make sure to use quantities defined on $M$. For each $p\in M$ and $t\ge 0$ 
choose a shortest geodesic $c(t), 0\le t \le d(p, q)$ w.r.t.~$\sigma_0$. Then we can integrate along $c(t)$
to get a desired formula for the first difference term. 
Employing these formulas and the exponential decay estimates for $\theta, \bar \theta$ and $\gamma$, and integrating (\ref{map-difference}), 
we deduce for all $t>0$
\ba \label{psi-estimate1}
\|\psi(\cdot, t)\|_{C^0} \le C_2(\|\bar \sigma_0
-\sigma_0\|_{C^{1+\mu}}+\|\bar \sigma_1-\sigma_1\|_{C^{1+\mu}})
\ea
with a positive constant $C_2$. 
Taking derivatives in 
(\ref{map-difference}) and arguing in similar fashions
we then obtain for all $t>0$
\ba \label{psi-estimate2}
\|\psi(\cdot, t)\|_{C^{3+\mu}} \le C_3(\|\bar \sigma_0
-\sigma_0\|_{C^{4+\mu}}+\|\bar \sigma_1-\sigma_1\|_{C^{4+\mu}})
\ea
with a positive constant $C_3$. 

Now we combine the above estimates to deduce for all $t\ge 1$
\ba
\|\bar \phi^*\bar \sigma(\cdot, t)-\phi^*\sigma(\cdot, t)\|_{C^{2+\mu}} &\le& \|\phi^*\gamma(\cdot, t)\|_{C^{2+\mu}}
+\|\psi^* \bar \sigma(\cdot, t)\|_{C^{2+\mu}}
\nonumber \\
&\le& C_4(\|\bar \sigma_0-\sigma_0\|_{C^{4+\mu}}
+\|\bar \sigma_1-\sigma_1\|_{C^{4+\mu}}) 
\ea
with a positive constant $C_4$. Taking the limit as 
$t\rightarrow \infty$ we then deduce 
\ba
\|\mc{F}(\bar \sigma_0, \bar \sigma_1)-
\mc{F}(\sigma_0, \sigma_1)\|_{C^{2+\mu}}
\le C_4 (\|\bar \sigma_0-\sigma_0\|_{C^{4+\mu}}
+\|\bar \sigma_1-\sigma_1\|_{C^{4+\mu}}). 
\ea
The general case of $l$ is similar.\qed \\

The definition domain of $\mc{F}$ is a domain in a Banach space, as the following lemma displays. 

\begin{lem} Let $l>0$ be an non-integer. Set $\mc{X}^l=\{(\gamma_0, \gamma_1):
\gamma_0\in C^l_o(\Lambda^3T^*M), 
\gamma_1 \in C^l_o(\Lambda^3T^*M), 
\gamma_1-\gamma_0 \in dC^{l+1}(\Lambda^3T^*M)\}$ and 
$\mc{Y}^l=\mc{X}^l \cap (C^l_o(\Lambda^3_+T^*M)\times 
C^l_o(\Lambda^3_+T^*M))$. Then $\mc{X}^l$ is 
a closed subspace of $C^l_o(\Lambda^3T^*M) \times C^l_o(\Lambda^3T^*M)$, and 
$\mc{Y}^l$ is a domain of $\mc{X}^l$. 
\end{lem}

\Pf First observe that $dC^{l+1}(\Lambda^3 T^*M)$ is a closed subspace of $C^l_o(\Lambda^3 T^*M)$. To show this, consider a sequence $\beta_k \in C^{l+1}(\Lambda^3 T^*M)$ such that $d\beta_k \rightarrow f$ in 
$C^l_o(\Lambda^3 T^*M)$. We solve the equation 
$\Delta \gamma_k=d\beta_k$ with $\gamma_k \perp 
\mc{H}^3$. Then we have $\Delta(\gamma_k-\gamma_{k'})=
d\beta_k-d\beta_{k'}$. It follows that $\|\gamma_k-
\gamma_{k'}\|_{C^{l+2}} \le C\|d\beta_k-d\beta_{k'}\|_{C^l}$. Hence $\gamma_k \rightarrow \gamma$ in $C^{l+2}$ for some $\gamma$. 
But $\Delta \gamma_k=d\beta_k$ implies $dd^* \gamma_k=d\beta_k$. Hence we infer $dd^* \gamma=f$, 
which implies $f\in dC^{l+1}(\Lambda^3T^*M)$. 

Obviously, the above closedness implies the desired 
closedness of $\mc{X}^l$.  By Lemma \ref{simple}, 
$\mc{Y}^l$ is a domain of $\mc{X}^l$. \qed \\

\begin{lem} \label{submanifold} Let $l>2$ be a non-integer. Then $\mc{T}^l$ is a smooth Banach submanifold of the Banach space 
$C^l_o(\Lambda^3 T^*M)$.
\end{lem}

We refer to [XY2] for the proof of this lemma.

\begin{theo} \label{smooth1} For given $\hat \sigma_0, \mu$ and $K$ as in Theorem \ref{gaugedconverg}, let ${\mc{U}}(\hat \sigma_0, \mu, K)$ denote the neighborhood in $\mc{Y}^{4+\mu}$ defined by the conditions in that theorem.  Then the map $\mc{F}: \mc{U}(\hat \sigma_0,
\mu, K) \rightarrow \mc{T}^{2+\mu}$ is smooth.  Moreover, for each $l>4+\mu$, the restriction $\mc{F}: \mc{U}(\hat \sigma_0,
\mu, K)\cap C^l \rightarrow \mc{T}^{l-2}$ is smooth.
\end{theo} 
\Pf This is a lengthy proof, which we break into three parts. In the first part, we decompose the difference form $\gamma$ in the proof of Theorem \ref{converg1} and derive the 
associated estimates. In the second part, we decompose the difference map $\psi$ in that proof and derive the associated estimates. In the last part, we draw the final conclusions. \\
1) We employ the notations in the proof of Theorem 
\ref{converg1} and set ${\bf p}=(\sigma_0, \sigma_1), {\bf q}=(\bar \sigma_0, \bar \sigma_1).$
 First observe for the equation 
(\ref{gamma-equation}) 
\ba
d(\Phi_{\sigma_0}(\bar \theta)-\Phi_{\sigma_0}(\theta))
=L_0\gamma+Q_0(\gamma, \gamma)
\ea
with
\ba
L_0 \gamma &=& d(D_{\sigma} A(\sigma_0,
\sigma, \theta, 
\nabla_{\sigma_0} \theta)(\gamma)+A(\sigma_0, \sigma, \gamma, \nabla_{\sigma_0}\theta)+ 
A(\sigma_0, \sigma, \theta,
\nabla_{\sigma_0} \gamma))
\ea
and
\ba
Q_0(\gamma, \gamma)
&=& d\int_0^1 tD^2_{\sigma} A(\sigma_0,
\sigma_0+\theta+st\gamma, \theta+st\gamma, 
\nabla_{\sigma_0} \theta+st \nabla_{\sigma_0} \gamma)(\gamma, \gamma)dsdt.
\ea
Similarly, we can write the sum of the second and third terms on the RHS of (\ref{gamma-equation}) as the sum of a linearized 
term and a quadratic term:
\ba
(\Delta_{\sigma_0}-\Delta_{\bar \sigma_0})\bar \theta-d(\Phi_{\bar \sigma_0}(\bar \theta)-\Phi_{\sigma_0}(\bar \theta))
=L_1(\bar \sigma_0-\sigma_0)+Q_1(\bar \sigma_0-\sigma_0, \bar \sigma-\sigma_0)+Q_2(\bar \sigma_0-\sigma_0,
\gamma),
\ea
where $L_1$ is independent of the quantities with bar.
It follows that 
\ba
\f{\p \gamma}{\p t}=-\Delta_{\sigma_0} \gamma
+L_0\gamma+L_1(\bar \sigma_0-\sigma_0)
+Q_0(\gamma, \gamma)+Q_1(\bar \sigma_0-\sigma_0,
\bar \sigma_0-\sigma_0)+Q_2(\bar \sigma_0-\sigma_0,
\gamma).
\ea
Note that $L_0, L_1, Q_0, Q_1$ and $Q_2$ are time-dependent, and converge to zero at exponential rate in suitable norms as $t\rightarrow \infty$. Now we consider the equation 
\ba \label{gamma1-equation}
\f{\p \gamma_1}{\p t}=-\Delta_{\sigma_0} \gamma_1
+L_0\gamma_1+L_1(\bar \sigma_0-\sigma_0)
\ea
with the initial condition $\gamma_1(0)=\gamma(0)=
(\bar \sigma_1-\sigma_1)-(\bar \sigma_0-\sigma_0)=(\bar \sigma_1-\bar \sigma_0)-(\sigma_1-\sigma_0)$
and the equation 
\ba \label{gamma2-equation}
\f{\p \gamma_2}{\p t}=-\Delta_{\sigma_0} \gamma_2
+L_0\gamma_2
+Q_0(\gamma, \gamma)+Q_1(\bar \sigma_0-\sigma_0,
\bar \sigma_0-\sigma_0)+Q_2(\bar \sigma_0-\sigma_0,
\gamma).
\ea
with the initial condition $\gamma_2(0)=0$. 
We can write $L_0$ in the following form
\ba
L_0 \gamma'=d(\Phi_0(\gamma')+\Phi_1(\nabla_{\sigma_0} \gamma'))
\ea
with $\Phi_0(\gamma')=D_{\sigma} A(\sigma_0,
\sigma, \theta, 
\nabla_{\sigma_0} \theta)(\gamma')+A(\sigma_0, \sigma, \gamma', \nabla_{\sigma_0}\theta)$ and $\Phi_1(\nabla_{\sigma_0}\gamma')= 
A(\sigma_0, \sigma, \theta,
\nabla_{\sigma_0} \gamma')$. As in the proof of Lemma \ref{linearlemma} we have $\|\Phi_1\|_{C^0, 
\sigma_0} \le C_0\|\theta\|_{C^0, \sigma_0}$. By the 
estimate for $\theta$ we can then apply Theorem \ref{closed} (or Theorem \ref{full}) to obtain a unique $\mc{C}^{4+\mu, (4+\mu)/2}$ solution $\gamma_1$ and a unique $\mc{C}^{4+\mu, (4+\mu)/2}$ solution $\gamma_2$
on $[0, \infty)$.   We also obtain the following estimates for all $t\ge 1$
\ba \label{gamma1}
\|\gamma_1\|_{\mc{C}^{4+\mu}(M\times [t-1, t])}
\le C_1e^{-\f{1}{2}\lambda_0t}(\|\bar \sigma_0-\sigma_0\|_{C^{4+\mu}}+\|\bar \sigma_1-
\sigma_1\|_{C^{4+\mu}})
\ea
and
\ba \label{gamma2}
\|\gamma_2\|_{\mc{C}^{4+\mu}(M\times [t-1, t])}
\le C_1e^{-\f{1}{2}\lambda_0t}(\|\bar \sigma_0-\sigma_0\|_{C^{4+\mu}}^2+\|\bar \sigma_1-
\sigma_1\|_{C^{4+\mu}}^2)
\ea
with a positive constant $C_1$ depending only on 
$\hat \sigma_0, \mu$ and $K$.  Obviously, there holds $\gamma=\gamma_1+\gamma_2$.

For $\sigma_0'\in C^{3+\mu}(\Lambda^3 T^*M)$ 
$\sigma_1'\in C^{4+\mu}(\Lambda^3 T^*M)$ such that 
$\sigma_1'-\sigma_0'$ is exact, we consider the following general version of (\ref{gamma1-equation})
\ba \label{Gamma1-equation}
\f{\p \gamma_1'}{\p t}=-\Delta_{\sigma_0} \gamma_1'
+L_0\gamma_1'+L_1\sigma_0'
\ea
with the initial condition $\gamma_1'(0)=
\sigma_1'-\sigma_0'$. Set ${\bf p}'=(\sigma_0', 
\sigma_1')$. Let $\Gamma_{1, \bf p}({\bf p}')$ denote the unique $\mc{C}^{4+\mu, (4+\mu)/2}$ solution on 
$M \times [0, \infty)$.  We have the following generalization of (\ref{gamma1})
\ba \label{Gamma1}
\|\Gamma_{1, {\bf p}}({\bf p}')\|_{\mc{C}^{4+\mu}(M\times [t-1, t])}
\le C_1e^{-\f{1}{2}\lambda_0t}(\|\sigma_0'\|_{C^{4+\mu}}+
\|\sigma_1'\|_{C^{4+\mu}}).
\ea 
We obviously have 
\ba \label{gamma-Gamma}
\gamma_1=\Gamma_{1, {\bf p}}({\bf q}-{\bf p}).
\ea

\noindent 2) Next we consider the equation (\ref{map-difference}). By the estimate (\ref{psi-estimate1}), we can achieve the following
by assuming $\|\bar \sigma_0-\sigma_0\|_{C^{1+\mu}}
+\|\bar \sigma_1-\sigma_1\|_{C^{1+\mu}}$ to be small enough: for each $p\in M$ and 
$t\ge 0$,  the distance between $\bar \phi(p, t)$ 
and $\phi(p, t)$ is less than half of the injectivity radius of $M$ (w.r.t.~$\sigma_0$).  Then we can handle the difference terms in (\ref{map-difference}) by unique 
shortest geodesics. The resulting quantities then retain the previous regularity and estimates. 
This way, we decompose the far right hand side of (\ref{map-difference}) into a linearized part and a quadratic part and deduce 
\ba
\f{d \psi}{d t}&=&\hat L_0\psi+\hat L_1\gamma
+\hat L_2 (\bar \sigma_0-\sigma_0)
+\hat Q_0(\psi, \psi)+\hat Q_1(\gamma, \gamma)
+\hat Q_2(\bar \sigma_0-\sigma, \bar \sigma_0-
\sigma_0)
+\hat Q_3(\psi, \gamma)\nonumber \\&&+\hat Q_4(\psi, \bar \sigma_0-
\sigma_0)+\hat Q_5(\gamma, \bar \sigma_0-\sigma_0),
\ea
where $\hat L_0, \hat L_1$ and $\hat L_2$ are independent of the quanitites with bar. Note that 
the involved operators $\hat L_0, \hat L_1, \hat L_2,$
$\hat Q_0$ etc.~ are all time-dependent and decay exponentially in suitable norms. 
We further write $\hat L_2\gamma=\hat L_2 \gamma_1+
\hat L_2\gamma_2$. Then we have $\psi=\psi_1+\psi_2$, where $\psi_1$ is the unique solution of the ODE
\ba \label{psi1}
\f{d \psi_1}{d t}=\hat L_0\psi_1+\hat L_1\gamma_1
+\hat L_2 (\hat \sigma_0-\sigma_0)
\ea
with the initial condition $\psi_1(0)=0$, and 
$\psi_2$ is the unique solution of the ODE 
\ba \label{psi2}
\f{d \psi_2}{d t}&=&\hat L_0\psi_2+\hat L_1\gamma_2
+\hat Q_0(\psi, \psi)+\hat Q_1(\gamma, \gamma)
+\hat Q_2(\bar \sigma_0-\sigma, \bar \sigma_0-
\sigma_0)+\hat Q_3(\psi, \gamma)+\hat Q_4(\psi, \bar \sigma_0-
\sigma_0) \nonumber \\&&+\hat Q_5(\gamma, \bar \sigma_0-\sigma_0),
\ea
with the initial condition $\psi_2(0)=0$. Employing 
the decay estimates for all the involved quantities 
we obtain the limits $\psi_1^{\infty}$ and $\psi_2^{\infty}$ of $\psi_1$ and $\psi_2$ respectively as $t\rightarrow \infty$, which satisfy
\ba \label{psi1-estimate}
\|\psi_1^{\infty}\|_{C^{3+\mu}}\le C_2 \|\bar \sigma_0-\sigma_0\|_{C^{4+\mu}}+||\bar \sigma_1-
\sigma_1\|_{C^{4+\mu}})
\ea
and 
\ba \label{psi2-estimate}
\|\psi_2^{\infty}\|_{C^{3+\mu}}\le C_2 \|\bar \sigma_0-\sigma_0\|_{C^{4+\mu}}^2+||\bar \sigma_1-
\sigma_1\|_{C^{4+\mu}}^2).
\ea
There holds $\psi_{\infty}=
\psi_{1, \infty}+\psi_{2, \infty}$. On the other hand, we have the following generalization of (\ref{psi1})
(analogous to (\ref{Gamma1-equation}))
\ba
\f{d}{dt} \Psi_{1, {\bf p}}({\bf p})=\hat L_0\Psi_{1, {\bf p}}({\bf p})
+\hat L_1\Gamma_{1, {\bf p}}({\bf p})
+\hat L_2 (\sigma_0')
\ea
with the initial condition $\Psi_{1, {\bf p}}({\bf p})
=0$, its limit $\Psi_{1, {\bf p}}^{\infty}({\bf p})$ as 
$t\rightarrow \infty$ and the estimate  
\ba \label{Psi1-estimate}
\|\Psi^{\infty}_{1, {\bf p}}({\bf p})\|_{C^{3+\mu}}\le C_2 \|\bar \sigma_0-\sigma_0\|_{C^{4+\mu}}+||\bar \sigma_1-
\sigma_1\|_{C^{4+\mu}}).
\ea
Thus $\Psi^{\infty}_{1, {\bf p}}$ is a bounded linear operator. There holds $\psi_1^{\infty}=\Psi^{\infty}_{1, {\bf p}}(\bar \sigma_0-\sigma_0, \bar \sigma_1-\sigma_1)$. \\
3) Now we calculate 
\ba
\mc{F}({\bf q})-\mc{F}({\bf p})&=&\mc{F}(\bar \sigma_0, \bar \sigma_1)-\mc{F}(\sigma_0, \sigma_1) =
\bar \phi_{\infty}^* \bar \sigma_0-\phi_{\infty}^*
\sigma_0 =\bar\phi_{\infty}^* 
(\bar \sigma_0-\sigma_0)+\psi_{\infty}^* \sigma_0
\nonumber \\
&=& \phi_{\infty}^* (\bar \sigma_0-\sigma_0)
+\psi_{\infty}^*(\bar \sigma_0-\sigma_0)
+\psi_{1,\infty}^*
\sigma_0+\psi_{2,\infty}^* \sigma_0.
\ea
By the estimates (\ref{psi2-estimate}) and (\ref{psi1-estimate}) we infer that $\mc{F}$ is differentiable at $\bf p$.  Moreover, we have 
\ba
\mc{D}_{\bf p} \mc{F}=\phi_{\infty}^* \pi_0
+\Psi^{\infty}_{1,{\bf p}}(\cdot)^*
\sigma_0,
\ea
where $\pi_0(\gamma', \gamma'')=\gamma'$.  Adapting the above arguments to handle the difference $\mc{D}_{\bf q} \mc{F}-\mc{D}_{\bf p} \mc{F}$, we deduce that 
$\mc{D}_{\bf p} \mc{F}$ is Lipschitz continuous. it follows that $\mc{F}$ is $C^1$ as a map into $C^{2+\mu}_o(\Lambda^3T^*M)$. By Lemma \ref{submanifold},
it is also $C^1$ as a map into $\mc{T}^{2+\mu}$.  

The above scheme can easily be extended to higher order derivatives of $\mc{F}$, and we derive that 
$\mc{F}$ is $C^{\infty}$.  Applying the $C^{l, l/2}$ estimates we then obtain the claimed smoothness of the $C^l$ restriction of $\mc{F}$. \qed \\

Finally we set $\mc{F}(\sigma_1)=\sigma_{\infty}$, where 
$\sigma_{\infty}$ is the limit of the Laplacian flow 
given in Theorem \ref{converg2}.  The following theorem contains Theorem \ref{limitmap} 
in Introduction as a special case.

\begin{theo} Let $0<\mu<1, K>0, \hat 
\sigma_0$ a given $C^{4+\mu}$ $G_2$-structure on 
$M$, and $\hat r_0$ be given by Theorem \ref{converg2}.
Let $4+\mu \le l \le \infty$ be a non-integer.
Then the limit map $\mc{F}: B^l_{K, \hat r_0}(\hat \sigma_0) \rightarrow \mc{T}^{l-2}$ is smooth. Moreover, there holds
\ba 
\pi_{l-2} \circ \mc{F}=\Pi_.
\ea
(In the case $l=\infty$, we have the convention $\infty-2=\infty$.)
\end{theo}

\Pf By the proof of Theorem \ref{converg2} there holds $\mc{F}(\sigma_1)=\mc{F}(\hat \Xi_{\hat \sigma_0}(\sigma_1), \sigma_1)$.  Hence the claimed 
smoothness follows from Theorem \ref{smooth1} 
and the smoothness of $\hat \Xi_{\hat \sigma_0}$. 
By Theorem \ref{converg2}, $\mc{F}(\sigma_1)$ is 
$C^{l-1}$ isotopic to $\hat \Xi_{\hat \sigma_0}(\sigma_1)$.  It follows that 
$\pi_{l-2}(\mc{F}(\sigma_1))=\pi_{l-2}(\hat \Xi_{\hat \sigma_0}(\sigma_1))=\Pi(\sigma_1)$.  \qed \\

\section*{Appendix: Space-time Function Spaces }

Let $M$ be a compact manifold of dimension 
$n\ge 1$.   Let a background Riemannian metric $g_*$ on $M$ be given. 
We assume that it has the required smoothness in 
each individual situation below. The norms defined below depend on the choice of $g_*$, but we an easily relate 
the norms w.r.t.~one background metric to those w.r.t.~another background metric.

Each tensor bundle $E$ associated with the tangent bundle $TM$ is equipped with the natural metric induced 
from $g_*$ and the natural connection $\nabla$ induced from the Levi-Civita connection (still called the Levi-Civita connection). We'll use these metric and connection in the definitions below.  
 
In this Appendix, we define various H\"older spaces of $E$-valued functions (i.e.~sections of $E$) used in this article. In particular, we define spacetime H\"{o}lder 
spaces which play a crucial role in our parabolic theory.
We basically follow the definitions given in [Y3]. 
Note that it is only for convenience of presentation that we restrict to tensor bundles associated with 
the tangent bundle. Our theory extends straightforwardly to a general vector bundle over $M$, which is equipped with a metric and a metric-compatible connection. 

From now on we fix a tensor bundle $E$. 

\subsection*{$C^{l}$-spaces}

Let $k\geq 0$ be an integer. We define the space $C^k(E)$ to be the space of continuous sections $\zeta$ of $E$ that have up to $k$-th order continuous  covariant derivatives, and define the 
$C^k$ norm as follows
\ba
\parallel  \zeta \parallel_{C^k(E)}=\sum_{i=0}^k\sup_M|\nabla^i \zeta|.
\ea
Equipped with this norm, the space $C^k(E)$ is a Banach 
space. \\

\noindent {\bf Remark} We write this norm as 
$\| \zeta\|_{C^k(E), g_*}$, if we need to indicate the 
background metric $g_*$. We replace the subscript $g_*$ by $\sigma$ if $g_*$ is the induced metric of a $G_2$-structure $\sigma$, i.e.~$g_*=g_{\sigma}$. Similar notations are also used 
for the other norms in this paper. \\

Next let $0<\mu<1$. We define the H\"older semi-norm $[\zeta]_\mu$ of a section $\zeta$ of $E$:
\ba
[\zeta]_\mu=\sup_{p, q\in M, 
0<d(p, q)\le 1} \sup_{\gamma} \frac{|P_\gamma(\zeta(p))-\zeta(q)|}{d(p,q)^\mu},
\ea
where $\gamma$ runs through all piecewise $C^1$-curves in $M$ going from $p$ to $q$ and having length not exceeding $2d(p, q)$, and $P_\gamma$ denotes the parallel transport along $\gamma$. (Alternatively, we can restrict to geodesics $\gamma$. Then we obtain an 
equivalent seminorm.) Note that the condition $[\zeta]_\mu<\infty$ can be interpreted as a fractional 
differentiability. 

Let $l=k+\mu$ for an integer $k\ge 0$. The H\"older space $C^{l}(E)$ consists of sections $\zeta$ of $C^k(E)$ with $[\nabla^k \zeta]_\mu<\infty$.  The norm $\parallel \zeta \parallel_{C^{l}}$ is defined as 
\ba
\parallel \zeta\parallel_{C^{l}}=\parallel \zeta\parallel_{C^k}+[\nabla^k \zeta]_\mu.
\ea
Equipped with this norm, $C^{l}(E)$ is a Banach space.

\subsection*{$\mc{C}^{m,l}$-spaces}
Let $I$ be a bounded closed interval of 
$\mathbb R$ with coordinate $t$. We abbreviate  for the derivative $\frac{\partial}{\partial t}$ to $\partial_t$. Let $\pi$ be the projection of $M\times I$ onto $M$ and $\pi^*E$ be the pull-back of $E$ to $M\times I$.

 For integers $m\geq 0, l\geq 0$, we define $\mc{C}^{m,l}(\pi^*E)$ to be the space of sections $\zeta$ of $\pi^*E$ which have continuous partial derivatives of the form $\partial_t^j\nabla^i \zeta$ with $i+2j\leq m$ and $ j\leq l$. When dealing with parabolic equations of the type of 
the heat equation, it is natural to count one time derivative as two space derivatives. This is the underlying reason for the above factor $2$ in front of $j$. This factor also appears below for the same reason.

The norm $\parallel\cdot\parallel_{C^{m,l}(\pi^*E)}$ is defined as follows
\ba
\parallel \zeta \parallel_{C^{m,l}(\pi^*E)}=\sup_{M\times I}\sum_{i+2j\leq m, j\leq l}|\partial^j_t\nabla^i \zeta|.
\ea

\noindent {\bf Remark} 
We often abbreviate the above notation to $\|\zeta\|_{C^{m,l}}$. We also write it as 
$\|\zeta\|_{C^{m,l}(M\times I)}$ if we need to emphasize the base domain.   Similar abbreviations and notations are also used for other norms or spaces in this paper.\\
 
It is easy to show that equipped with the above norm, $\mc{C}^{m,l}(\pi^*E)$ is a Banach space.
 
\subsection*{$\mc{C}^{l,m/2}$-spaces}
We now introduce ``fractional" differentiability in both the time and space directions.  For $0<\mu<1$, we define the $\mu$-H\"older semi-norm in the time direction
\ba
[\zeta]_{\mu, M\times I, I}=\sup_{p\in M, 0\leq t_2-t_1\leq 1, t_1,t_2\in I}\frac{|\zeta(p,t_2)-\zeta(p,t_1)|}{|t_2-t_1|^\mu}
\ea
and the $\mu$-H\"older semi-norm in the space direction  
\ba
[\zeta]_{\mu, M\times I, M}=\sup_{t\in I}[\zeta(\cdot, t)]_{\mu}.
\ea

Now let $l$ and $m$ be nonnegative non-integers with $2l\geq m$.  We define the space $\mc{C}^{l,m/2}(\pi^*E)$ as the space of sections in $\mc{C}^{[l],[m/2]}(\pi^*E)$ with finite $\mc{C}^{l,m/2}$-norm, which is defined as follows 
\ba \label{lm-norm}
\parallel \zeta \parallel_{\mc{C}^{l,m/2}(\pi^*E)}=\sum_{i+2j\leq [l], j\leq [m/2]}\max_{M\times I}|\partial^j_t\nabla^i \zeta|+<\zeta>^{(l)}_{M\times I, M}+<\zeta>^{(m/2)}_{M\times I, M},
\ea
with the $(l)$-H\"older semi-norm in the space direction 
\[<\zeta>^{(l)}_{M\times I, M}=\sum_{i+2j=[l], j\leq [m/2]}[\partial^j_t\nabla^i \zeta]_{l-[l],M\times I, M},\]
and the $(m/2)$-H\"older semi-norm in the time direction
\[<\zeta>^{(m/2)}_{M\times I, I}=\sum_{0<m-i-2j<2,i\leq [l]}[\partial^j_t\nabla^i \zeta]_{(m-i-2j)/2, M\times I, I}.\]
It is easy to show that, equipped with the norm 
(\ref{lm-norm}), $\mc{C}^{l,m/2}$ is a Banach space. 

Of particular importance is the case $l=m$, i.e.~the spaces $\mc{C}^{l,l/2}$.  They are the natural spaces for formulating a priori estimates for solutions of parabolic equations, see Theorem \ref{full}.  The formula for the $\mc{C}^{l,l/2}$-norm takes a slightly simpler form:
\ba
\parallel \zeta \parallel_{\mc{C}^{l,l/2}(\pi^*E)}=\sum_{i+2j\leq [l]} \max_{M\times I}|\partial^j_t\nabla^i \zeta|+<\zeta>^{(l)}_{M\times I, M}+<\zeta>^{(l/2)}_{M\times I, I},
\ea
 with the $(l)$-H\"older semi-norm in the space direction 
\ba
<\zeta>^{(l)}_{M\times I, M}=\sum_{i+2j=[l]}[\partial_t^j\nabla^i \zeta]_{l-[l],M\times I, M}
\ea
and the $(l/2)$-H\"older semi-norm in the time direction 
\ba
<\zeta>^{(l/2)}_{M\times I, I}=\sum_{0<l-i-2j<2}[\partial^j_t\nabla^i \zeta]_{\frac{l-i-2j}{2},M\times I, I}.
\ea
For example, we have for $0<\mu<1$
\ba
\|\zeta\|_{\mc{C}^{\mu, \mu/2}}=\max_{M\times I} 
|\zeta|+[\zeta]_{\mu,M\times I, M}+[\zeta]_{\f{\mu}{2},
M\times I, I},
\ea
\ba
\|\zeta\|_{\mc{C}^{1+\mu, (1+\mu)/2}}
=\max_{M\times I} (|\zeta|+|\nabla \zeta|)
+[\nabla \zeta]_{\mu,M\times I, M}+
[\zeta]_{\f{1+\mu}{2}, M\times I, I}+[\nabla \zeta]_{\f{\mu}{2},
M\times I, I},
\ea
and
\ba
\|\zeta\|_{\mc{C}^{2+\mu, (2+\mu)/2}} &=&
\max_{M\times I}(|\zeta|+
|\p_t \zeta|+|\nabla \zeta|+
|\nabla^2 \zeta|)
+([\partial_t \zeta]_{\mu,M\times I, M}
+[\nabla^2 \zeta]_{\mu,M\times I, M}) \nonumber \\
&&+
([\partial_t \zeta]_{\frac{\mu}{2},M\times I, I}
+[\nabla \zeta]_{\frac{1+\mu}{2},M\times I, I}
+[\nabla^2 \zeta]_{\frac{\mu}{2},M\times I, I}).
\ea
   
Finally we present another separate definition which is used in the formulation of some results in this paper.\\

\noindent {\bf Definition 11.1} 
We define the inverse tensor $(g_2)_{g_1}^{-1}$ of a Riemannian metric $g_2$ w.r.t.~another Riemannian metric $g_1$ as follows.
There holds $g_2(v_1, v_2)=g_1(Av_1, v_2)$ for 
a section $A$ of $T^*M \otimes TM=Hom(TM, TM)$. By the positive definiteness of $g_2$, $A$ is invertible at 
each point. We set 
\ba
(g_2)^{-1}_{g_1}(v_1, v_2)=g_1(A^{-1}v_1, v_2).
\ea
Then we define 
\ba
\|g_2^{-1}\|_{C^l, g_1}=\|(g_2)^{-1}_{g_1}\|_{C^l, g_1}.
\ea
We write it as $\|g_2^{-1}\|_{C^l}$, if the metric 
$g_1$ is clear from the context. Note that the eigenvalues of $g_2^{-1}$ w.r.t.~$g_1$ are the reciprocals of the eigenvalues of $g_2$ w.r.t.~$g_1$.

\end{document}